\documentclass{amsart}
\usepackage{amsthm,amsmath,amssymb}
\usepackage{stmaryrd}           
\usepackage{mathrsfs}           
\usepackage{hyperref}
\usepackage[all,2cell]{xy}
\newcommand{\M}{\ensuremath{\mathscr{M}}}
\newcommand{\N}{\ensuremath{\mathscr{N}}}
\newcommand{\p}{\ensuremath{\mathscr{P}}}
\newcommand{\sS}{\ensuremath{s\mathcal{S}}}

\newcommand{\I}{\ensuremath{\mathscr{I}}}
\newcommand{\V}{\ensuremath{\mathscr{V}}}
\newcommand{\K}{\ensuremath{\mathscr{K}}}
\newcommand{\B}{\ensuremath{\mathscr{B}}}
\newcommand{\C}{\ensuremath{\mathscr{C}}}
\newcommand{\D}{\ensuremath{\mathscr{D}}}
\newcommand{\E}{\ensuremath{\mathscr{E}}}
\newcommand{\Oscr}{\ensuremath{\mathscr{O}}}
\newcommand{\R}{\ensuremath{\mathscr{R}}}
\newcommand{\W}{\ensuremath{\mathscr{W}}}
\newcommand{\Q}{\ensuremath{\mathscr{Q}}}
\newcommand{\DD}{\ensuremath{\mathbf{\Delta}}}
\newcommand{\ep}{\ensuremath{\varepsilon}}
\newcommand{\tep}{\ensuremath{\widetilde{\varepsilon}}}
\newcommand{\gm}{\gamma}
\newcommand{\Id}{\ensuremath{\operatorname{Id}}}
\newcommand{\ob}{\ensuremath{\operatorname{ob}}}
\newcommand{\map}{\operatorname{Map}}

\newcommand{\tva}{two-variable adjunction}
\newcommand{\tvva}{two-variable \V-adjunction}
\newcommand{\tvhva}{two-variable $\Ho\V$-adjunction}
\newcommand{\csmhc}{closed symmetric monoidal homotopical category}
\newcommand{\vheovhc}{\V-homotopical equivalence of \V-homotopical categories}

\newcommand{\colim}{\operatorname{colim}}
\newcommand{\hocolim}{\operatorname{hocolim}}
\newcommand{\uhc}{\operatorname{uhocolim}}
\newcommand{\holim}{\operatorname{holim}}
\newcommand{\uhl}{\operatorname{uholim}}
\newcommand{\coeq}{\operatorname{coeq}}
\newcommand{\eqlzr}{\operatorname{eq}}

\newcommand{\lan}{\operatorname{Lan}}

\newcommand{\tholan}{\mathbf{L}\operatorname{Lan}}
\newcommand{\ran}{\operatorname{Ran}}

\newcommand{\thoran}{\mathbf{R}\operatorname{Ran}}

\newcommand{\Ho}{\operatorname{Ho}}

\newcommand{\ten}{\otimes}

\newcommand{\thoten}{ \ensuremath{\overset{\mathbf{L}}{\ten}} }
\newcommand{\xodot}{\,\overline{\odot}\,}
\newcommand{\hoodot}{ \ensuremath{\overset{\mathbb{L}}{\odot}} }
\newcommand{\thoodot}{ \ensuremath{\overset{\mathbf{L}}{\odot}} }
\providecommand{\oast}{\varoast} 
\newcommand{\xoast}{\,\overline{\oast}\,}
\newcommand{\hooast}{ \ensuremath{\overset{\mathbb{L}}{\oast}} }
\newcommand{\thooast}{ \ensuremath{\overset{\mathbf{L}}{\oast}} }
\newcommand{\smsh}{\wedge}
\newcommand{\hosmash}{ \ensuremath{\overset{\mathbb{L}}{\wedge}} }

\newcommand{\cten}[1]{\{#1\}}
\newcommand{\hocten}[1]{\mathbb{R}\cten{#1}}
\newcommand{\thocten}[1]{\mathbf{R}\cten{#1}}

\newcommand{\homr}{\operatorname{hom}_r}
\newcommand{\hohomr}{\ensuremath{\mathbb{R}\homr}}
\newcommand{\thohomr}{\ensuremath{\mathbf{R}\homr}}
\newcommand{\homl}{\operatorname{hom}_\ell}
\newcommand{\hohoml}{\ensuremath{\mathbb{R}\homl}}
\newcommand{\thohoml}{\ensuremath{\mathbf{R}\homl}}

\newcommand{\iso}{\cong}
\newcommand{\eqv}{\simeq}
\newcommand{\dn}{\downarrow}
\newcommand{\op}{\ensuremath{^{\mathit{op}}}}
\newcommand{\adj}{\dashv}
\newcommand{\Adj}{\;\dashv\;}
\newcommand{\sB}{B_\bullet}

\newcommand{\sX}{\ensuremath{X_\bullet}}
\newcommand{\sY}{\ensuremath{Y_\bullet}}

\newdir{ >}{{}*!/-10pt/@{>}}

\newcommand{\too}[1][]{\ensuremath{\overset{#1}{\longrightarrow}}}
\renewcommand{\to}[1][]{\ensuremath{\overset{#1}{\rightarrow}}}
\newcommand{\we}{\ensuremath{\overset{\sim}{\longrightarrow}}}
\newcommand{\leftwe}{\ensuremath{\overset{\sim}{\longleftarrow}}}
\newcommand{\cof}[1][]{\ensuremath{\overset{#1}{\rightarrowtail}}}

\newcommand{\cohto}{\rightsquigarrow}
\newcommand{\uto}{\ensuremath{\,\text{---}\,}}
\newcommand{\maps}{\ensuremath{\colon}}
\newcommand{\spam}{\,:\!}

\makeatletter

\def\defthm#1#2{\newtheorem{#1}[thm]{#2}}
\ifx\SK@label\undefined\let\SK@label\label\fi
\let\old@label\label
\let\your@thm\@thm
\def\@thm#1#2#3{\gdef\currthmtype{#3}\your@thm{#1}{#2}{#3}}
\def\currthmtype{}
\AtBeginDocument{%
  \let\old@label\label%
  \def\label#1{%
    {\let\your@currentlabel\@currentlabel%
      \edef\@currentlabel{\currthmtype}%
      \old@label{label@name@#1}}%
    \old@label{#1}}
  }
\def\autoref#1{\ref*{label@name@#1}~\ref{#1}}

\newtheorem{thm}{Theorem}[section]

\defthm{cor}{Corollary}
\defthm{prop}{Proposition}
\defthm{lem}{Lemma}
\defthm{metathm}{Meta-theorem}
\theoremstyle{definition}
\defthm{defn}{Definition}
\defthm{metadefn}{Meta-definition}
\defthm{exmp}{Example}
\defthm{ceg}{Counterexample}
\theoremstyle{remark}
\defthm{notn}{Notation}
\defthm{rem}{Remark}

\let\c@equation\c@thm
\numberwithin{equation}{section}

\makeatother

\begin{document}

\title{Homotopy limits and colimits and enriched homotopy theory}
\author{Michael Shulman}
\email{shulman@math.uchicago.edu}
\address{Department of Mathematics,
         University of Chicago,
         5734 S.\ University Ave.,
         Chicago, IL, 60637,
         U.S.A.}

\begin{abstract}
  Homotopy limits and colimits are homotopical replacements for the
  usual limits and colimits of category theory, which can be
  approached either using classical explicit constructions or the
  modern abstract machinery of derived functors.  Our first goal in
  this paper is expository: we explain both approaches and a proof of
  their equivalence.  Our second goal is to generalize this result to
  enriched categories and homotopy weighted limits, showing that the
  classical explicit constructions still give the right answer in the
  abstract sense.  This result partially bridges the gap between
  classical homotopy theory and modern abstract homotopy theory.  To
  do this we introduce a notion of ``enriched homotopical
  categories'', which are more general than enriched model categories,
  but are still a good place to do enriched homotopy theory.  This
  demonstrates that the presence of enrichment often simplifies rather
  than complicates matters, and goes some way toward achieving a
  better understanding of ``the role of homotopy in homotopy theory.''
\end{abstract}

\maketitle

\tableofcontents

\section{Introduction}
\label{sec:introduction}

Limits and colimits are fundamental constructions in category theory
and are used to some degree nearly everywhere in mathematics.
However, when the category in question has a ``homotopy theory,''
limits and colimits are generally not homotopically well-behaved.  For
example, they are not invariant under homotopy equivalence.  One of
the central needs of homotopy theory is thus for a replacement with
better properties, usually called a ``homotopy limit.''  This is part
of a more general need for a ``homotopy theory of diagrams'' which
generalizes the standard category theory of diagrams.

There are two natural candidates for a definition of homotopy limit.
The first is defined for any category enriched over topological
spaces, or over simplicial sets.  In this case there are well-known
explicit constructions of homotopy limits and colimits dating back to
the classical work~\cite{bk}; for a more modern exposition,
see~\cite{hirschhorn}.  These homotopy limits are objects satisfying a
``homotopical'' version of the usual universal property: instead of
representing \emph{commuting} cones over a diagram, they represent
``homotopy coherent'' cones.  This universal property is \emph{local}
in that it characterizes only a single object, although usually the
homotopy limit can be extended to a functor.

More recently, treatments of axiomatic homotopy theory have studied
categories which are equipped with some notion of ``weak
equivalence''.  One can ``invert'' the weak equivalences homotopically
to get a simplicially enriched localization, but it is usually more
convenient to deal with the original category and its weak
equivalences directly.  In this context, we have the notion of
``derived functor,'' which is a universal homotopical approximation to
some given functor.  From this point of view, it is natural to define
a homotopy limit to be a derived functor of the usual limit functor.
This sort of homotopy limit has a \emph{global} universal property
that refers to all possible homotopical replacements for the limit
functor.

Both approaches have advantages.  For instance, the classical
``local'' construction, when expressed as a bar construction, has a
natural filtration which gives rise to spectral sequences; this makes
it very tractable computationally.  On the other hand, the universal
property of the ``global'' constructions makes it easier to obtain
coherence and preservation results.  Ideally, we would like the two
definitions to agree, at least up to homotopy, so that we can use
whichever is most convenient for a given purpose.  Various comparison
proofs have been given, dating back to the original work~\cite{bk}.

The dichotomy between these two definitions of homotopy limit is, in
fact, part of a larger disconnect between the techniques of classical
homotopy theory, on the one hand, and those of modern axiomatic
homotopy theory, on the other.  In order to bridge this gap, we begin
by considering the question of what is meant by ``homotopy theory'' in
the first place.  We may say loosely that homotopy theory is the study
of categories equipped with ``weak equivalences'' using a number of
technical tools, among which are the following.

\begin{itemize}
\item Notions of ``homotopy'' and ``homotopy equivalence''.
\item Subcategories of ``good'' objects on which the weak equivalences
  are homotopy equivalences.
\item Special classes of maps usually called ``fibrations'' and
  ``cofibrations''.
\end{itemize}

The roles of the first two of these, at least, can already be seen in
two basic problems of homotopy theory: the construction of
\emph{derived categories} or \emph{homotopy categories} by inverting
the weak equivalences, and the construction of \emph{derived functors}
that carry the ``homotopical information'' from point-set level
functors.  Of course, a homotopy category can always be constructed
formally, but its hom-sets will in general be large; the problem is to
ensure that its hom-sets are small.  This is most often done by
considering a subcategory of ``good'' objects on which the weak
equivalences are homotopy equivalences, so that the derived hom-sets
are quotients of the original hom-sets by the homotopy relation.

Similarly, derived functors are generally constructed by applying the
original functor to a subcategory of ``good'' objects on which the
original functor preserves weak equivalences; this frequently happens
because the functor preserves homotopies and homotopy equivalences,
and the weak equivalences between good objects are homotopy
equivalences.  The role of the cofibrations and fibrations is less
basic, but it has much to do with the construction of homotopy limits
and colimits and therefore the subject of the present paper; we will
encounter them especially in \S\ref{sec:comp-dhks}--\ref{sec:comp-bar}
and \S\ref{sec:cofibrancy}.

Probably the most powerful and best-known axiomatization of homotopy
theory is Quillen's theory of \emph{model categories}, first
introduced in~\cite{quillen}.  This theory combines all of the above
structures in a neat package that fits together very precisely.  Any
two of the classes of fibrations, cofibrations, and weak equivalences
determine the third via lifting and factorization properties.  The
``good'' objects are the \emph{cofibrant} and/or \emph{fibrant} ones,
which means that the map from the initial object is a cofibration or
that the map to the terminal object is a fibration.  The factorization
properties are used to define \emph{cylinder} and \emph{path} objects,
which are then used to define notions of homotopy, and it turns out
that weak equivalences between fibrant and cofibrant objects are
precisely the homotopy equivalences; in this way the homotopy category
is shown to have small hom-sets.

One also has the notion of a \emph{Quillen adjoint pair} of functors,
which is defined by preservation of cofibrations and fibrations.  It
then follows that the left adjoint preserves weak equivalences between
cofibrant objects, and the right adjoint preserves weak equivalences
between fibrant objects, so that both have canonical derived functors.
Moreover, these derived functors form an adjoint pair between the
homotopy categories.

The theory of Quillen model categories and Quillen adjoints is very
powerful when it applies.  Many interesting categories have useful
model structures and many interesting functors are part of Quillen
adjunctions.  However, not infrequently it happens that a category of
interest does not admit a useful model structure, or that a functor of
interest is not Quillen with respect to any known model structure.
This is frequently true for diagram categories, on which limit and
colimit functors are defined.  Moreover, even when a category has a
model structure, sometimes the cofibrations or fibrations it supplies
are unnecessarily strong for the desired applications.  For example,
in the usual model structure on topological spaces, the cofibrations
are the relative cell complexes, but for many purposes the weaker
classical notion of Hurewicz cofibration suffices.

For these reasons, we are driven to re-analyze classical homotopy
theory and seek new axioms which are not as restrictive as those of
model categories, but which nevertheless capture enough structure to
enable us to ``do homotopy theory.''  Several authors have studied
such more general situations, in which either the fibrations or
cofibrations are missing, or are not characterized as precisely by
lifting properties.

More radically, the authors of the recent book~\cite{dhks} have begun
a study of categories with only weak equivalences, which they call
\emph{homotopical categories}, and of subcategories of ``good
objects,'' which they call \emph{deformation retracts}.  They do not
consider notions of homotopy at all.  On the other side of the coin,
the authors of~\cite{hcct} study abstract notions of homotopy, modeled
by simplicially enriched categories, but without reference to weak
equivalences.

We believe that a synthesis of such approaches is needed, especially
to deal with the important subject of \emph{enriched} homotopy theory,
which has rarely been considered in recent treatments of abstract
homotopy theory, and then usually only for simplicial enrichments.
(For example,~\cite{bourn-cordier} gives a general formulation of
homotopy limits for simplicially enriched categories.)  We make no
claim to have fully achieved such a synthesis.  In the present paper
we only develop these ideas as far as is necessary to compare various
definitions of homotopy limits and colimits, in the enriched as well
as the unenriched setting.  But even this much, we believe, is quite
illuminating, and sheds some light on what one might call ``the role
of homotopy in homotopy theory.''

The relevance of these ideas to homotopy limits may be described
(ahistorically) as follows.  Suppose that we are given a category \M\
equipped with weak equivalences, and also a notion of homotopy and
good objects which realize the weak equivalences as homotopy
equivalences.  We would like categories of diagrams in \M\ to inherit
a similarly good homotopy theory.  If we replace a diagram by one
which is objectwise good, then a natural transformation which is an
objectwise weak equivalence will be replaced by one which is an
objectwise homotopy equivalence, but the homotopy inverses of its
components may not fit together into a natural transformation.
However, they do in general fit together into a transformation which
is ``natural up to coherent homotopy''.  This notion can be made
precise using bar constructions, which in fact allow us to realize
homotopy-coherent transformations as true natural transformations
between ``fattened up'' diagrams.

Thus, taking these fattened-up diagrams to be the ``good'' objects
gives a homotopy theory on the diagram category, and in particular, a
way to construct homotopy limits and colimits in the abstract or
``global'' sense mentioned earlier.  But since coherent
transformations are defined using bar constructions, the homotopy
limits we obtain in this way are essentially the same as the ``local''
constructions from classical homotopy theory.  Then, since the
classical constructions work just as well in an enriched context, they
can be used to define a global notion of ``weighted homotopy limit''
there as well.

The main goal of this paper is to make the above sketch precise,
first showing that the local and global notions agree, and then
extending the comparison to the enriched situation.  If this were all
we had to do, this paper would be much shorter, but due to the
existing disconnect between classical and axiomatic homotopy theory
mentioned above, we must first develop many of the tools we need,
especially in the enriched context.  Moreover, we also want to make it
clear that our approach is essentially equivalent to both existing
approaches.

Thus, there are two threads in this paper interweaving back and forth:
the concrete construction of homotopy limits and colimits, on the one
hand, and the general philosophy and results about abstract homotopy
theory, on the other.  We encourage the reader to keep the above
sketch in mind, to avoid becoming lost in the numerous technicalities.

This paper can be divided into six main parts.  In the first two
parts, we restrict ourselves to simplicially enriched categories and
unenriched diagrams, both for ease of exposition and because this is
the only context in which the existing global approaches are defined.
The remaining four parts are concerned with generalizing these results
to the enriched context.

The first part, comprising
\S\S\ref{sec:global-derived}--\ref{sec:middle-derived}, consists of
general theory about derived functors.  Except for
\S\ref{sec:middle-derived}, this is largely a review of material
from~\cite{dhks}.  The second part, comprising
\S\S\ref{sec:global-hoco}--\ref{sec:coherent-transf}, is devoted to an
exposition of the existing ``global'' and ``local'' approaches to
homotopy limits (given in \S\ref{sec:global-hoco} and
\S\S\ref{sec:local-hoco}--\ref{sec:bar}, respectively), and to two
proofs of comparison.

One general comparison proof was given in an early online draft of
\cite{dhks}, but no longer appears in the published book.  In
\S\ref{sec:comp-dhks} we describe this proof (with details postponed
to appendix~\ref{sec:ugly-proof}).  Then in \S\ref{sec:comp-bar} we
present an alternate proof, making use of the fact that the local
definition can also be described using a two-sided bar construction.
In \S\ref{sec:coherent-transf} we explain how this proof is in line
with the intuition of ``homotopy coherent transformations'' described
above.

In the last four parts of the paper, we consider the question of
\emph{enriched} homotopy theory.  The limits which arise naturally in
enriched category theory are \emph{weighted} limits, but no general
theory of weighted homotopy limits and colimits exists in the
literature.  For several reasons, weighted homotopy limits are
considerably more subtle than ordinary ones.

In particular, the local definition using a two-sided bar construction
enriches very naturally, but all global approaches we are aware of
founder on various difficulties in the enriched case.  Our approach
provides a solution to this difficulty: in favorable cases, for
categories that are both enriched and homotopical, we can show
directly that the enriched version of the classical bar construction
defines a derived functor in the global sense.  We do this in the
third part of the paper, comprising
\S\S\ref{sec:enriched}--\ref{sec:htpy-ten}.  This justifies the
classical enriched bar construction from the global point of view of
abstract homotopy theory.

The results in the third part are likely to be adequate for many of
the applications, but they are not fully satisfactory since they do
not seriously address the relationship between the enrichment and the
homotopical structure.  In particular, while weighted limits and
colimits are themselves enriched functors, the derived functors we
produce are only unenriched functors.  The fourth and fifth parts of
the paper remedy this defect.

In the fourth part, comprising
\S\S\ref{sec:enriched-tvas}--\ref{sec:derived-via-enriched}, we
develop a theory of ``enriched homotopical categories'', applying
ideas developed in~\cite[ch.~4]{hovey} for enriched model categories
to the more general homotopical categories of~\cite{dhks}.  The main
points are an enriched analogue of the theory of ``two-variable
adjunctions'' given in~\cite{hovey}, and a general theory of
enrichments of homotopy categories derived from enriched homotopical
categories.  This abstract framework may be seen as a proposed
axiomatization of enriched homotopy theory, which incorporates weak
equivalences, enrichment, enriched homotopy and homotopy equivalences,
and good objects.

Then in the fifth part, comprising
\S\S\ref{sec:gener-tens-prod}--\ref{sec:homot-theory-enrich}, we apply
this theory to prove that the total derived homotopy limits and
colimits can be enriched over the homotopy category of the enriching
category.  We also study the homotopy theory of diagrams more
generally, including the behavior of enriched homotopy colimits with
respect to functors between domain categories.

The sixth and final part of the paper deals with some technical
details which arise in consideration of enriched homotopy limits.  In
\S\ref{sec:cofibrancy}, we introduce cofibrations and fibrations into
our abstract framework in order to describe the necessary cofibrancy
conditions on the shape category.  Model theoretic cofibrations, when
present, generally suffice, but are often too restrictive for the
applications.  Moreover, there is an added problem that cofibrant
approximation functors generally do not preserve the enrichment; in
\S\ref{sec:objectwise-good}, we consider various ways to get around
this problem.  For clarity, such technical considerations are swept
under the rug with an axiomatization in parts three and five.

It is important to note that everything in this paper has dual
versions for limits and colimits, right and left Kan extensions, and
cotensor and tensor products.  We work mostly with the latter, since
we believe they are easier to conceptualize for various reasons,
leaving dualizations to the reader in many cases.  However, we
emphasize that our results are often more necessary for limits than
colimits, since simplifying model structures on diagram categories are
more often present in the context of colimits.






The reader is assumed to have some knowledge of category theory,
including the formalisms of ends and Kan extensions, as described
in~\cite{maclane}.  Some familiarity with model categories is also
expected, especially in the second part of the paper.  This includes
an acquaintance with the various model structures that exist on
diagram categories.  For good introductions to model categories,
see~\cite{dwyer_spal}, \cite{hovey}, and~\cite{hirschhorn}.  We assume
that all model categories are complete and cocomplete and have
functorial factorizations.

We also use a number of ideas and results from~\cite{dhks}, but no
prior familiarity with that work is required if the reader is willing
to take a few of its results on faith.

In various places, we use the terminology and techniques of enriched
category theory, as described in~\cite{kelly} and~\cite{dubuc}.  Their
use is concentrated in the second half of the paper, however, and we
attempt to explain these concepts as they arise.  Enriched category
theory, being the sort of category theory which nearly always arises
in practical applications, is at least as important as the more
elementary unenriched theory.

To reiterate, then, the upshot of this paper is that the classical
two-sided bar construction is almost always the ``best'' way to
construct homotopy limits and colimits, for the following reasons.
\begin{itemize}
\item It satisfies a local universal property, representing homotopy
  coherent cones;
\item It also satisfies a global universal property, making coherence
  and preservation results easy to obtain;
\item It comes with a natural filtration by simplicial degree, giving
  rise to spectral sequences which make it computable; and
\item It works just as well in enriched situations, which are
  frequently the cases of most interest.
\end{itemize}

\textbf{Acknowledgements.}
I would like to thank my advisor, Peter May, for many useful
conversations about the role of enrichment and innumerable helpful
comments on incomprehensible drafts of this paper; Gaunce Lewis, for
bringing to my attention the subtleties regarding the tensor product
of \V-categories; and Phil Hirschhorn, for supplying me with part of
an early draft of~\cite{dhks} and helping me understand the proof of
\autoref{dhks-comp}.


\section{Global definitions of derived functors}
\label{sec:global-derived}

We begin by considering categories with a suitably well-behaved notion
of \emph{weak equivalence}.  The most common requirement is a
2-out-of-3 property, but in~\cite{dhks} it was found more technically
convenient to assume a slightly stronger property, as follows.

\begin{defn}[{\cite[\S33]{dhks}}]
  A \emph{homotopical category} is a category \M\ equipped with a
  class of morphisms called \emph{weak equivalences} that contains all
  the identities and satisfies the \emph{2-out-of-6 property}: if $hg$
  and $gf$ are weak equivalences, then so are $f$, $g$, $h$, and
  $hgf$.
\end{defn}

It follows easily that the weak equivalences in a homotopical category
include all the isomorphisms and satisfy the usual 2-out-of-3
property.  The stronger 2-out-of-6 property is satisfied by all
examples that arise, including all model categories and categories of
diagrams in model categories.

Any category can be made a homotopical category in a trivial way by
taking the weak equivalences to be precisely the isomorphisms;
in~\cite{dhks} these are called \emph{minimal} homotopical categories.
If $\M$ is a homotopical category and $\D$ is any small category, the
diagram category $\M^\D$ has a homotopical structure in which a map of
diagrams $\alpha\maps  F\to F'$ is a weak equivalence if $\alpha_d\maps Fd\to
F'd$ is a weak equivalence in \M\ for every object $d\in\D$; we call
these \emph{objectwise} weak equivalences.

Any homotopical category \M\ has a \emph{homotopy category} $\Ho\M$,
obtained by formally inverting the weak equivalences, with a
localization functor $\gamma\maps \M\to\Ho\M$ which is universal among
functors inverting the weak equivalences.  In general, $\Ho\M$ need
not have small hom-sets, which is a major impediment in many
applications.  Usually, however, other methods are available to ensure
local smallness, as discussed briefly in \S\ref{sec:introduction}.
For example, as already mentioned, if \M\ admits a Quillen model
structure, then $\Ho\M$ will have small hom-sets.  In
\autoref{vhc->loc-small}, we will see another way to ensure that homotopy
categories have small hom-sets, using an enriched notion of homotopy.

Homotopical categories are a reasonable axiomatization of the first
ingredient for homotopy theory: a subcategory of weak equivalences.
In this section and the next, we will follow~\cite{dhks} in developing
the theory of derived functors for homotopical categories.

Thus, suppose that \M\ and \N\ are homotopical categories.  Let
$\Ho\M$ and $\Ho\N$ be their homotopy categories, with localization
functors $\gamma\maps \M\to\Ho\M$ and $\delta\maps \N\to\Ho\N$.
Following~\cite{dhks}, we say that a functor $F\maps \M\to\N$ which
preserves weak equivalences is \emph{homotopical}.  We will also use
this term for functors $\M\to\Ho\N$ which take weak equivalences to
isomorphisms; this is equivalent to considering $\Ho\N$ as a minimal
homotopical category.  If $F$ is homotopical, then by the universal
property of localization $\delta F$ induces a unique functor
$\widetilde{F}\maps \Ho\M\to\Ho\N$ making the following diagram commute.
\begin{equation}\label{eq:derived}
  \xymatrix{
    \M \ar[r]^F \ar[d]_\gamma & \N \ar[d]^\delta \\
    \Ho\M \ar[r]^{\widetilde{F}} & \Ho\N
  }
\end{equation}

Frequently, however, we encounter functors $\M\to\N$ which are not
homotopical, but for which we would still like a derived functor
$\Ho\M\to\Ho\N$.  One natural approach is to relax the requirement
that~(\ref{eq:derived}) commute and instead ask that it be filled with
a \emph{universal} transformation.  This amounts to asking for a right
or left Kan extension of $\delta F$ along $\gamma$.  The reversal of
handedness in the following definition is unfortunate, but both
terminologies are too well-established to be changed.

\begin{defn}
  A right Kan extension of $\delta F$ along $\gamma$ is called a
  \emph{total left derived functor} of $F$ and denoted $\mathbf{L}F$.
  Dually, a left Kan extension of $\delta F$ along $\gamma$ is called
  a \emph{total right derived functor} of $F$ and denoted
  $\mathbf{R}F$.
\end{defn}

All the types of derived functor we will define in this section come
in two types, ``left'' and ``right.''  For clarity, from now on we
will define only the left versions, leaving the obvious dualizations
to the reader.

Often we are interested in derived functors that, like
$\tilde{F}\circ\gamma$ in~(\ref{eq:derived}), are defined on \M\ and
not just on $\Ho\M$.  By the universal property of localization, the
existence of a total left derived functor is equivalent to the
existence of a ``left derived functor,'' defined as follows.  We also
write $\mathbf{L}F$ for this notion, despite the difference in the
domain.

\begin{defn}
  A \emph{left derived functor} of $F$ is a functor $\mathbf{L}F\maps 
  \M\to\Ho\N$ equipped with a comparison map $\mathbf{L}F\to \delta F$
  such that $\mathbf{L}F$ is homotopical and terminal among
  homotopical functors equipped with maps to $\delta F$.
\end{defn}

Clearly, both left derived functors and total left derived functors
are unique up to unique isomorphism when they exist.

We may also want to ``lift'' the target of the derived functor from
$\Ho\N$ to \N.  This is not always possible, but frequently it is, and
so to clarify the notions we introduce the following oxymoronic
definition.

\begin{defn}\label{def:psder}
  A \emph{point-set left derived functor} is a functor
  $\mathbb{L}F\maps \M\to\N$ equipped with a comparison map $\mathbb{L}F\to
  F$ such that the induced map $\delta\mathbb{L}F\to \delta F$ makes
  $\delta \mathbb{L}F$ into a left derived functor of $F$.
\end{defn}

Such a functor, when it exists, is only unique ``up to homotopy''.  By
no means do all derived functors ``lift'' to the point-set level.  For
example, Quillen functors out of model categories whose factorizations
cannot be made functorial, such as those of~\cite{pross}, have derived
functors, but not point-set derived functors.  However, we will need
our model categories to have functorial factorization in order to have
an ``objectwise cofibrant replacement'' on diagrams, so all our
model-theoretic derived functors will have point-set versions.

Finally, we might replace the universal property on the level of
homotopy categories in the preceding three definitions by a
``homotopical universal property'' on the point-set level.

\begin{defn}[{\cite[41.1]{dhks}}]\label{def:appx}
  A \emph{left approximation} of $F$ is a homotopical functor
  $F'\maps \M\to\N$ equipped with a comparison map $F'\to F$ that is
  ``homotopically terminal'' among such functors, in a sense made
  precise in~\cite[\S38]{dhks}.
\end{defn}

It is proven in~\cite[\S38]{dhks} that objects with such a homotopical
universal property are ``homotopically unique'' in a suitable sense.

Both point-set derived functors and approximations are functors
$\M\to\N$, but in general the existence of one need not imply the
existence of the other.  If $F'$ is an approximation, it has a
universal property referring only to other functors landing in $\N$,
while the universal property of a derived functor refers to all
functors landing in $\Ho\N$, even those that do not factor through \N\
and are thus only ``functorial up to homotopy''; thus $\delta F'$ is
not necessarily a derived functor.  On the other hand, if
$\mathbb{L}F$ is a point-set derived functor, then for any homotopical
functor $F'$ mapping to $F$ we have a unique natural transformation
$\delta F'\to\delta\mathbb{L}F$, but it need not lift to a natural
zig-zag from $F'$ to $\mathbb{L}F$ in $\N$; thus $\mathbb{L}F$ is not
necessarily a left approximation.

However, it is useful to note that if \N\ is a minimal homotopical
category (for example, \N\ could be $\Ho\p$ for some other homotopical
category \p), then $\delta$ is an equivalence and all the above
notions coincide.  Thus the results of~\cite{dhks} about
approximations can usually be applied to derived functors as well.  We
will state our results for point-set derived functors, except when it
becomes necessary to pass all the way to total derived functors.  It
should also be noted that the way in which most derived functors are
constructed, which we will explain in the next section, produces both
derived functors and approximations with equal ease.

\section{Derived functors via deformations}
\label{sec:deformations}

Having defined left and right derived functors in the previous
section, we of course would like to know how to construct such things.
The easiest and most common way to produce a left derived functor of
$F$ applies when the categories involved have model category
structures for which $F$ is left Quillen.  In this case, for any
cofibrant replacement functor $Q$ on \M, the composite $FQ$ can be
shown to be both a point-set left derived functor and a left
approximation.  In this section we follow~\cite{dhks} in axiomatizing
the properties of the cofibrant replacement functor $Q$ which make
this work, thus disentangling the notion of ``good object'' from the
rest of the model category machinery.

\begin{defn}[{\cite[\S40]{dhks}}]\label{def:deform}
  Let $F\maps \M\to\N$ be a functor between homotopical categories.
  \begin{itemize}
  \item A \emph{left deformation} of \M\ is a functor $Q\maps \M\to\M$
    equipped with a natural weak equivalence $q\maps Q\we \Id_\M$; it
    follows from the 2-out-of-3 property that $Q$ is homotopical.
  \item A \emph{left deformation retract} is a full subcategory $\M_Q$
    containing the image of some left deformation $(Q,q)$.  We will
    always use such parallel notation for deformations and deformation
    retracts.
  \item A \emph{left $F$-deformation retract} is a left deformation
    retract $\M_Q$ such that the restriction of $F$ to $\M_Q$ is
    homotopical.
  \item A deformation into an $F$-deformation retract is called an
    \emph{$F$-deformation}.
  \item If there exists a left $F$-deformation retract, we say $F$ is
    \emph{left deformable}.
  \end{itemize}
\end{defn}

\begin{rem}\label{easy-deform}
  By the 2-out-of-3 property, a deformation $(Q,q)$ is an
  $F$-deformation if and only if $FQ$ is homotopical and $Fq Q\maps FQQ\to
  FQ$ is a natural weak equivalence.
\end{rem}

If $F$ is left Quillen, a cofibrant replacement functor serves as an
$F$-deformation.  Our choice of notation for a left deformation is
intended to suggest that it is a ``generalized cofibrant
replacement.''

As an example of a left deformation that is not a cofibrant
replacement, it is well-known that when computing derived tensor
products in homological algebra (i.e.\ the modern version of
$\operatorname{Tor}$), it suffices to use \emph{flat} resolutions
rather than projective ones, although the latter are the cofibrant
objects in the usual model structure.  Similarly, in topological
situations, it often suffices to consider various weaker forms of
cofibrancy, such as the ``$h$-cofibrations'' of~\cite[\S4.1]{pht} or
the ``tame spectra'' of~\cite[\S{}I.2]{ekmm}.  We will also need to
use deformations that do not arise directly from a model structure to
compute homotopy limits and colimits.

The following two results justify the above definitions.  There are,
of course, dual versions for right deformations and right derived
functors.

\begin{prop}\label{deform-hoequiv}
  If $\M_Q$ is a left deformation retract of \M, then the inclusion
  $I\maps \M_Q\hookrightarrow \M$ induces an equivalence of categories
  $\Ho(\M_Q)\eqv \Ho\M$.
\end{prop}
\begin{proof}
  Let $(Q,q)$ be a left deformation of \M\ into $\M_Q$.  Since the
  inclusion $I\maps \M_Q\hookrightarrow \M$ and $Q\maps \M\to\M_Q$ are both
  homotopical, they induce functors $\Ho\M_Q\rightleftarrows \Ho\M$.
  The natural weak equivalence $q\maps \Id_\M\to Q$ then descends to
  homotopy categories to give natural isomorphisms $\Id_{\Ho\M}\iso
  IQ$ and $\Id_{\Ho\M_Q}\iso QI$, so $Q$ and $I$ form the desired
  equivalence.
\end{proof}

\begin{prop}[{\cite[\S41]{dhks}}]\label{deform->derived}
  If $(Q,q)$ is a left $F$-deformation, then $\mathbb{L}F = FQ$ is
  both a point-set left derived functor of $F$ and a left
  approximation of $F$.
\end{prop}
\begin{proof}
  This is an obvious generalization of the well-known proof for
  Quillen functors, e.g.~\cite[8.3.6]{hirschhorn}.
\end{proof}

Intuition from classical homotopy theory leads us to expect that we
could construct the homotopy category $\Ho\M$ using a suitable
deformation $(Q,R)$ of the hom-functor $\M(-,-)$.  However, unless \M\
is \emph{enriched}, this has little chance of working.  This is
because if \M\ is not enriched, $\M(-,-)$ takes values only in
\textbf{Set}, in which the only sensible weak equivalences are the
isomorphisms.  We will see in \S\ref{sec:derived-enrichment}, however,
that if \M\ is enriched over a category with its own suitable notion
of weak equivalence, this intuition is often correct.

In the rest of this section, we will mention a few auxiliary notions
and results about deformations which will be useful in the second half
of the paper.  All of this material is from~\cite{dhks}, and can be
skipped on a first reading.

We first consider the question of when an adjunction between deformable
functors descends to homotopy categories.  We frequently write $F\adj
G$ for an adjunction $F\maps \M_1\rightleftarrows\M_2 \spam G$.

\begin{defn}[{\cite[43.1]{dhks}}]
  An adjunction $F\adj G$ between homotopical categories is
  \emph{deformable} if $F$ is left deformable and $G$ is right
  deformable.
\end{defn}

\begin{prop}[{\cite[43.2 and 44.2]{dhks}}]\label{derived-adjn}
  If $F\maps \M_1\rightleftarrows\M_2 \spam G$ is a deformable adjunction, then
  there is a unique adjunction
  \[\mathbf{L}F\maps  \Ho\M_1 \rightleftarrows \Ho\M_2 \spam \mathbf{R}G\]
  which is compatible with the given one and the localization functors
  of $\M_1$ and $\M_2$.
\end{prop}

\begin{defn}[{\cite[41.6]{dhks}}]
  If $\alpha\maps  F\to F'$ is a natural transformation between functors on
  homotopical categories, a \emph{derived natural transformation} of
  $\alpha$ is a natural transformation $\mathbf{L}\alpha\maps 
  \mathbf{L}F\to\mathbf{L}F'$ such that
  \begin{enumerate}
  \item $\mathbf{L}F\to \delta F$ and $\mathbf{L}F'\to \delta F'$ are
    derived functors of $F$ and $F'$, respectively;
  \item the square
    \begin{equation*}
      \xymatrix{
        \mathbf{L}F \ar[r]\ar[d]_{\mathbf{L}\alpha} &
        \delta F\ar[d]^{\delta(\alpha)} \\
        \mathbf{L}F' \ar[r] &
        \delta F'
      }
    \end{equation*}
    commutes; and
  \item this square is terminal among commutative squares having
    $\delta(\alpha)$ on the right and such that the two functors on
    the left are homotopical.
  \end{enumerate}
\end{defn}

One can similarly define point-set derived natural transformations,
approximations of natural transformations, and so on.

\begin{prop}\label{deform-trans}
  If $\alpha\maps  F\to G$ is a natural transformation between left
  deformable functors $\M\to\N$:
  \[\UseTwocells \xymatrix{\M \rtwocell^F_{G}{\alpha} & \N}\]
  and there exists a deformation $(Q,q)$ of \M\ which is both an
  $F$-deformation and a $G$-deformation, then $\alpha$ has a derived
  natural transformation
  \[\UseTwocells \xymatrix{\Ho\M \rtwocell^{\mathbf{L}F}_{\mathbf{L}G}{\;\mathbf{L}\alpha} & \Ho\N}\]
  given by $\alpha Q$.

  Moreover, if we have another natural transformation
  \[\UseTwocells \xymatrix{\M \rtwocell^{G}_{H}{\beta} & \N}\]
  and $(Q,q)$ is also an $H$-deformation, then
  $\mathbf{L}(\beta\alpha) \iso (\mathbf{L}\beta)(\mathbf{L}\alpha)$.
\end{prop}
\begin{proof}
  See~\cite[41.6]{dhks} for the first statement.  The second statement
  is straightforward.
\end{proof}

We would like to know when the operation of ``left deriving'' is
``functorial'', i.e.\ when it preserves composition of functors.  This
is not the case in general, but there are special circumstances under
which it is true.

\begin{defn}[{\cite[\S42]{dhks}}]\label{def:deform-pair}
  A composable pair $(F_1, F_2)$ of left deformable functors
  $\M\to[F_1]\N\to[F_2]\p$ is called \emph{locally left deformable} if
  there are deformation retracts $\M_Q$ and $\N_Q$ such that
  \begin{itemize}
  \item $F_1$ is homotopical on $\M_Q$;
  \item $F_2$ is homotopical on $\N_Q$; and
  \item $F_2 F_1$ is homotopical on $\M_Q$.
  \end{itemize}
  If in addition, $F_1$ maps $\M_Q$ into $\N_Q$, the pair is called
  \emph{left deformable}.
\end{defn}

Note that a pair of left Quillen functors is always left deformable,
so in the well-behaved world of model categories there is no need for
this notion.

\begin{prop}[{\cite[42.4]{dhks}}]\label{deform-pair}
  If $(F_1,F_2)$ is a left deformable pair, then for any left derived
  functors $\mathbb{L}F_1$ and $\mathbb{L}F_2$ of $F_1$ and $F_2$,
  respectively, the composite $\mathbb{L}F_2\circ \mathbb{L}F_1$ is a
  left derived functor of $F_2F_1$.
\end{prop}

When the functors involved have adjoints, part of the work to check
deformability of pairs can be shifted across the adjunction, where it
often becomes easier.  First we need one more definition.

\begin{defn}[{\cite[33.9]{dhks}}]
  A homotopical category is \emph{saturated} if any map which becomes
  an isomorphism in its homotopy category is a weak equivalence.
\end{defn}

Any model category is saturated, e.g. by~\cite[1.2.10]{hovey}.
Moreover, if \M\ is saturated and \D\ is small, then
by~\cite[33.9(v)]{dhks} the functor category $\M^\D$ (with objectwise
weak equivalences) is also saturated.

The following is the reason for defining ``locally left deformable''
in addition to ``left deformable.''

\begin{prop}[{\cite[42.5]{dhks}}]\label{deform-adjn}
  Let $F_1\maps \M_1\rightleftarrows\M_2 \spam  G_1$ and
  $F_2\maps \M_2\rightleftarrows\M_3 \spam G_2$ be deformable adjunctions in
  which the pairs $(F_1,F_2)$ and $(G_2,G_1)$ are locally left and
  right deformable, respectively, while the $\M_i$ are saturated.
  Then the pair $(F_1,F_2)$ is left deformable if and only if the pair
  $(G_2,G_1)$ is right deformable.
\end{prop}

The composition of left derived and right derived functors is a much
more subtle question, and even under the best of circumstances (such
as left and right Quillen functors) composition need not be preserved.
See, for example,~\cite[Counterexample 0.0.1]{pht}.  In~\cite{pht}
this question is dealt with using the tools of ``middle derived
functors'', which we introduce in the next section.

\section{Middle derived functors}
\label{sec:middle-derived}

We have defined derived functors by a universal property and then
introduced deformations as a means of constructing derived functors,
but a more historically faithful introduction would have been to
remark that although most functors are not homotopical, they do
frequently preserve weak equivalences between ``good objects.''
Therefore, a natural ``homotopical replacement'' for $F$ is to first
replace the domain object by an equivalent good object and then apply
$F$.

From this point of view, we are left wondering what the real
conceptual difference between left derived functors and right derived
functors is, and whether if a functor has both left and right derived
functors, the two must agree or be related in some way.  In search of
a unified notion of derived functor, we might be motivated to drop the
restriction that $q$ be a \emph{single} weak equivalence in favor of a
zigzag of such, leading to the following definitions.

\begin{defn}\label{def:middle-deform}
  Let $F\maps \M\to\N$ be a functor between homotopical categories.
  \begin{itemize}
  \item A \emph{middle deformation} of \M\ is a functor $D\maps \M\to\M$
    equipped with a natural zigzag of weak equivalences
    \[\xymatrix{D\ar@{-}[r]^\sim & D_k\ar@{-}[r]^\sim &
      \dots\ar@{-}[r]^\sim & D_1\ar@{-}[r]^\sim & \Id_\M}.\]
    As before, it follows from the 2-out-of-3 property that $D$ is
    homotopical, as are all the $D_i$.
  \item A \emph{middle deformation retract} is a full subcategory
    $\M_D$ for which there is a middle deformation, as above, such
    that $\M_D$ contains the image of $D$, and each $D_i$ maps $\M_D$
    into itself.
  \item A \emph{middle $F$-deformation retract} is a middle
    deformation retract on which $F$ is homotopical.  A \emph{middle
      $F$-deformation} is a middle deformation into a middle
    $F$-deformation retract.
  \item If $F$ has a middle $F$-deformation $D$, then the composite
    $FD$, which descends to homotopy categories, is a \emph{middle
      derived functor} of $F$.
  \end{itemize}
\end{defn}

The assumption that each $D_i$ preserves $\M_D$ allows us to prove a
version of \autoref{deform-hoequiv} for middle deformations.

\begin{prop}\label{middef-hoequiv}
  If $\M_D$ is a middle deformation retract, then the inclusion
  \[I\maps \M_D\hookrightarrow \M\] induces an equivalence
  $\Ho(\M_D)\eqv \Ho\M$.
\end{prop}
\begin{proof}
  The functors $I\maps \M_D\to\M$ and $D\maps \M\to\M_D$ are both homotopical
  and so induce functors $\Ho\M_D\rightleftarrows \Ho\M$.  The assumed
  natural zigzag of weak equivalences provides an isomorphism $ID\iso
  \Id_{\Ho\M}$, and also $DI\iso \Id_{\Ho\M_D}$ since each
  intermediate functor $D_i$ maps $\M_D$ to itself.
\end{proof}

Clearly, left and right deformations are both special cases of middle
deformations, so if we could prove that middle derived functors are
unique, it would follow that if a functor has both left and right
derived functors, they agree.  Unfortunately, this is not the case, as
the following example shows.

\begin{ceg}\label{eg:middle-nonunique}
  Let \M\ and \N\ be homotopical categories, and let \I\ be the
  category with two objects $0$ and $1$ and one nonidentity arrow
  $0\to 1$ which is a weak equivalence.  Let $F^0, F^1\maps \M\to\N$ be two
  functors such that $F^0$ is left deformable and $F^1$ is right
  deformable and let $\alpha\maps F^0\to F^1$ be a natural transformation.

  Define $F\maps \M\times\I\to\N$ using $\alpha$ in the obvious way, and
  give $\M\times\I$ the product homotopical structure.  Then if $Q$ is
  a left deformation for $F^0$, a left deformation for $F$ is given by
  projecting to $\M\times\{0\}$ and then applying $Q$; thus $F$ is
  left deformable and $\mathbf{L}F \iso \mathbf{L}F^0$.  Dually, $F$
  is also right deformable and $\mathbf{R}F\iso \mathbf{R}F^1$.

  Therefore, $F$ is both left and right deformable, but the only
  relationship between its left and right derived functors is that
  they are connected by a natural transformation.
\end{ceg}

Since left and right derived functors are both middle derived
functors, it follows that middle derived functors are not unique, and
therefore they cannot be expected to satisfy any universal property
similar to that for left and right derived functors.

Nevertheless, there are important cases in which left and right
derived functors \emph{do} agree, and which include many situations
encountered in applications.  One frequently encountered situation is
described by the following definition.

\begin{defn}
  Let $\M_{D'}\subset \M$ and $\M_{D} \subset \M$ be middle
  deformation retracts.  We say $\M_{D'}$ is a \emph{middle
    sub-deformation retract} of $\M_{D}$ if $\M_{D'}\subset \M_{D}$
  and moreover $\M_{D'}$ is a middle deformation retract of $\M_{D}$.
\end{defn}

\begin{rem}
  This condition is equivalent to the existence, for a suitable choice
  of middle deformations $D$ and $D'$, of a zigzag of weak
  equivalences:
  \begin{equation}\label{eq:submiddle-zigzag}
    \xymatrix{D'\ar@{-}[r]^\sim & \dots\ar@{-}[r]^\sim &
      D\ar@{-}[r]^\sim & \dots\ar@{-}[r]^\sim & \Id_\M}
  \end{equation}
  in which all the weak equivalences between $D'$ and $D$ lie in
  $\M_D$.
\end{rem}

\begin{prop}\label{submiddle-unique}
  Let $F\maps \M\to\N$ and let $\M_{D'}$ be a middle sub-deformation
  retract of $\M_D$ such that $\M_D$ is a middle $F$-deformation
  retract (and hence so is $\M_{D'}$).  Then the corresponding middle
  derived functors of $F$ are equivalent.
\end{prop}
\begin{proof}
  Apply $F$ to~(\ref{eq:submiddle-zigzag}).  Since all the weak
  equivalences between $D'$ and $D$ lie in $\M_D$, on which $F$ is
  homotopical, the image of that part of the zigzag gives a zigzag of
  weak equivalences between $FD'$ and $FD$, which are the two middle
  derived functors in question.
\end{proof}

This obvious-seeming result has some less obvious-seeming
consequences, such as the following.

\begin{prop}\label{submiddle-left=right}
  Let $F\maps \M\to\N$ and let $(Q,q)$ and $(R,r)$ be left and right
  $F$-deformations, respectively.  If $Q$ maps $\M_R$ into itself, or
  $R$ maps $\M_Q$ into itself, then $\mathbf{L}F\iso \mathbf{R}F$.
\end{prop}
\begin{proof}
  Assume $Q$ maps $\M_R$ into itself.  Then we have the following
  square of weak equivalences
  \begin{equation*}
    \xymatrix{
      Q \ar[r]^\sim \ar[d]_\sim & QR \ar[d]^\sim \\
      \Id_\M \ar[r]^\sim & R
    }
  \end{equation*}
  in which the image of $QR$ lands in $\M_Q\cap\M_R$.  This shows that
  $\M_Q\cap\M_R$ is a middle sub-deformation retract of both $\M_Q$
  and $\M_R$.  Therefore, both $\mathbf{L}F=FQ$ and $\mathbf{R}F=FR$
  are equivalent to the middle derived functor $FQR$, hence equivalent
  to each other.  The case when $R$ maps $\M_Q$ into itself is dual.
\end{proof}

\begin{cor}\label{model-left=right}
  If \M\ and \N\ are model categories and $F\maps \M\to\N$ is both left and
  right Quillen, then $\mathbf{L}F\iso \mathbf{R}F$.
\end{cor}
\begin{proof}
  If we take $Q$ and $R$ to be the cofibrant and fibrant approximation
  functors coming from the functorial factorization in \M, then
  because $q\maps Q\to\Id_\M$ is a trivial fibration, $Q$ maps $\M_R$ into
  itself, and dually $R$ maps $\M_Q$ into itself.  Therefore we can
  apply \autoref{submiddle-left=right}.
\end{proof}

Thus, if our categories are equipped with a ``preferred'' model
structure, any given functor has at most one ``canonical'' derived
functor.  Later on, in \S\ref{sec:derived-tvas}, when we consider
enriched homotopical categories and corresponding notions of homotopy,
we will see another common situation in which a preferred extra
structure on a homotopical category ensures that there is a
``canonical'' derived functor.

We conclude this section with the following additional example of
\autoref{submiddle-unique}.  This is taken from~\cite{pht}, which has
many examples of functors that have both left and right derived
functors that agree.

\begin{exmp}
  It is shown in~\cite[9.1.2]{pht} that the category $\K_B$ of
  ex-spaces over a space $B$ is equipped with the following collection
  of middle deformation retracts (although the authors of~\cite{pht}
  do not use that terminology).
  
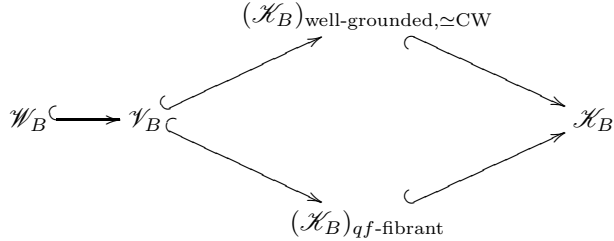
\begin{equation*}
    \xymatrix{
      && (\K_B)_{\text{well-grounded,}\eqv\text{CW}} \ar@{_{(}->}[dr] \\
      \W_B \ar@{^{(}->}[r] &
      \V_B \ar@{^{(}->}[ur] \ar@{^{(}->}[dr] &&
      \K_B\\
      && (\K_B)_{qf\text{-fibrant}} \ar@{^{(}->}[ur].
    }
  \end{equation*}
  The detailed definitions of all these categories, most of which
  refer to a certain ``$qf$-model structure'' on $\K_B$, need not
  concern us.  The important points are that in $\W_B$, weak
  equivalences are homotopy equivalences, so that we have an
  equivalence $h\W_B\eqv\Ho\K_B$, where $h\W_B$ is the quotient of
  $\W_B$ by homotopy, and that the definition of $\W_B$ is entirely
  classical, making no reference to the $qf$-model structure.  The
  result~\cite[9.2.2]{pht}, which says that suitably nice functors
  $\K_A\to\K_B$ have middle derived functors induced by functors
  $h\W_A\to h\W_B$, can then be seen as an instance of
  \autoref{submiddle-unique}.
\end{exmp}

From the perspective of middle derived functors, the ``leftness'' or
``rightness'' of any particular derived functor reveals itself as
merely an accident---albeit a very common and very useful one!  For
example, in \S\ref{sec:wgt-htpy-lim} we will see how useful the
universal properties of left and right derived functors are in proving
compatibility relations.  The question of when a given middle derived
functor is left or right is a subtle one.  Experience, especially from
the case of model categories, leads us to expect that in general, left
adjoints will have left derived functors, while right adjoints will
have right derived functors, but we have no more general theory
explaining why this should be true.  However, see
also~\cite{shulman:dblderived}.


\section{Homotopy colimits as derived functors}
\label{sec:global-hoco}

We now begin the second part of the paper and descend to the specific
case of homotopy limits and colimits.  From the ``global'' point of
view, a homotopy limit should be a derived functor of the limit
functor.  Since the limit is a right adjoint and the colimit is a left
adjoint, we expect the one to have a right derived functor and the
other a left derived functor.

Now, if for all model categories \M\ and all small categories \D\, the
diagram category $\M^\D$ always had model structures for which the
colimit and limit functors were left and right Quillen, respectively,
then no further discussion of this approach would be necessary, at
least for model categories.  There have been steps in the direction of
a notion of model category with these properties; see, for
example,~\cite{thomason}.  But the more common notion of Quillen model
category does not have them.

However, there are important special cases in which the diagram
category does have a model structure and limit or colimit functors are
Quillen.
\begin{enumerate}
\item If \M\ is cofibrantly generated, then all categories $\M^\D$
  have a \emph{projective model structure} in which the weak
  equivalences and fibrations are objectwise
  (see~\cite[ch.~11]{hirschhorn} for details).  When this model
  structure exists, the \emph{colimit} functor is left Quillen on it.
\item If \M\ satisfies the stronger hypothesis of being
  \emph{sheafifiable} (see~\cite{beke}), then all $\M^\D$ have an
  \emph{injective model structure} in which the weak equivalences and
  \emph{co}fibrations are objectwise.  When this model structure
  exists, the \emph{limit} functor is right Quillen on it.
\item If \D\ is a \emph{Reedy category}, then for any model category
  \M, the category $\M^\D$ has a \emph{Reedy model structure} in which
  the weak equivalences are objectwise, but the cofibrations and
  fibrations are generally not (see~\cite[ch.~5]{hovey}
  or~\cite[ch.~15]{hirschhorn} for details).  If in addition \D\ has
  \emph{fibrant constants} as defined in~\cite[15.10.1]{hirschhorn},
  then the colimit functor is left Quillen for this model structure.
  Dually, if \D\ has cofibrant constants, then the limit functor is
  right Quillen.
\end{enumerate}

When a suitable model structure exists on $\M^\D$, the global
definition of homotopy limits or colimits is easy.  We simply apply a
fibrant or cofibrant replacement in the appropriate model structure
and take the usual limit or colimit.  However, for projective and
injective model structures, the cofibrant and fibrant replacements are
constructed by small object arguments and so are difficult to get a
handle on, making a more explicit construction desirable.

Reedy model structures have the advantage that fibrant and cofibrant
replacements are relatively easy to define and compute.  Moreover,
many common diagram shapes have Reedy structures---for example, the
category
\[(\cdot\leftarrow\cdot\rightarrow\cdot)\]
which indexes pushout
diagrams---so many global homotopy limits and colimits can be computed
in this way by simply replacing a few maps by fibrations or
cofibrations.

However, being Reedy is quite a special property of the diagram
category, so in the general case it is far too much to expect.
Moreover, although most model categories arising in nature are
cofibrantly generated, so that projective model structures exist, many
(such as those arising from topological spaces) are not sheafifiable,
so there is no known model structure for which the limit functor is
Quillen.  Thus more technical methods are needed to construct global
homotopy limit functors at all in this context.

The approach taken in~\cite{dhks} is to use a suitable ``homotopical
replacement'' for the shape category \D.  Recall that for any functor
$H\maps \D\to\E$ and any object $e$ of \E, the \emph{comma-category} of $H$
over $e$ has for its objects the arrows $Hd\to e$ in \E, and for its
arrows the arrows in \D\ whose images under $H$ form commutative
triangles.  We write $(H\dn e)$ for the comma-category.  When $H$ is
the identity functor of $\D$, we write $(\D\dn d)$ for the
comma-category; in this case it is also called the ``over-category'' of
$d$.  There is a dual comma-category $(e\dn H)$ which specializes to
the ``under-category'' $(d\dn\D)$.

We define the \emph{category of simplices} of a simplicial set $K$ to
be the comma category $\Delta K = (\Delta\dn K)$, where
$\Delta\maps \DD\to\sS$ is the standard cosimplicial simplicial set.
The objects of $\Delta K$ are all the simplices of $K$, of all
dimensions.  We write $\Delta\op K = (\Delta K)\op \iso (\Delta\op\dn
K)$.  The following is straightforward.
\begin{lem}[{\cite[22.10]{dhks}}]\label{simplices-reedy}
  For any $K$, the categories $\Delta K$ and $\Delta\op K$ are Reedy
  categories with fibrant and cofibrant constants, respectively.
\end{lem}
We have forgetful functors $\Pi\maps \Delta K\to\DD$ and $\Sigma =
\Pi\op\maps \Delta\op K\to\DD\op$, which send an $n$-simplex $\alpha$ to
the object $[n]$ of $\DD$ or $\DD\op$.

We are primarily interested in $\Delta K$ when $K$ is the nerve of a
category.  Recall that the \emph{nerve} of a small category \D\ is a
simplicial set $N\D$ whose $n$-simplices are the strings of $n$
composable arrows in \D.  We have $N\D_n = \mathbf{Cat}([n],\D)$ where
$[n]$ is a string of $n$ composable arrows.  Note that some authors
call $N\D$ the ``classifying space'' of \D\ and denote it $B\D$.
However, in addition to the potential for confusion with the use of
$B$ for the bar construction, we prefer to reserve that notation for
``delooping''-type operations.

If \D\ is a small category, write $\Delta\D = \Delta N\D$ and
$\Delta\op\D=\Delta\op N\D$.  Thus, the objects of $\Delta\D$ are the
strings of zero or more composable arrows in \D.  We have ``target''
and ``source'' functors $T\maps \Delta\D\to\D$ and $S\maps
\Delta\op\D\to\D$ which map an $n$-simplex of $N\D$
\[\alpha_0\too \alpha_1\too \dots\too \alpha_n\]
to $\alpha_n$, respectively~$\alpha_0$.

\begin{rem}\label{op-nerves}
  Not all authors define the nerve of a category the same way; some
  authors' $N\D$ is other authors' $N\D\op$.  Our conventions are
  those of~\cite{hirschhorn} and~\cite{dhks}, but opposite to those
  of~\cite{bk}.  A discussion of the precise effect of the choice of
  convention on the local definition of homotopy colimits can be found
  in~\cite[18.1.11]{hirschhorn}.
  
  Up to homotopy, however, the choice of convention doesn't matter,
  since the nerve of a category is weakly equivalent to the nerve of
  its opposite; see, for example,~\cite[14.1.6]{hirschhorn}.
  Moreover, upon passing to categories of simplices, we actually have
  an isomorphism $\Delta\D\iso\Delta(\D\op)$.  Thus homotopy limits
  and colimits can be defined using either convention, and the
  comparison proofs go through.  Note, however, a potential source of
  confusion: $\Delta\op\D$ and $\Delta(\D\op)$ are quite different
  categories.
\end{rem}

We summarize the functors relating $\D$ and its categories of
simplices in the following diagram.
\begin{equation}\label{eq:simpl-diag}
  \xymatrix{
    \DD & & \D & & \DD\op\\
    & \Delta\D\ar[ul]_\Pi\ar[ur]^T && \Delta\op\D\ar[ul]_S\ar[ur]^\Sigma
  }
\end{equation}
For any category \M, we have a functor $T^*\maps
\M^\D\to\M^{\Delta\D}$ given by precomposing with $T$, and similarly
for all other functors in the diagram.  If \M\ is cocomplete, we also
have \emph{left Kan extensions} such as $\lan_T\maps
\M^{\Delta\D}\to\M^\D$, defined to be a left adjoint to $T^*$.
See~\cite[ch.~X]{maclane} for a development of Kan extensions in
general.  Left Kan extensions are also sometimes called
\emph{prolongation} functors.

\begin{thm}[{\cite{dhks}}]\label{thm:dhks}
  For any model category \M\ and small category \D, the functor
  $\colim_\D$ preserves all weak equivalences between diagrams of the
  form $\lan_T F'$, where $F'$ is a Reedy cofibrant diagram in
  $\M^{\Delta\D}$.  Moreover, if $Q$ is a cofibrant replacement
  functor in the Reedy model structure on $\M^{\Delta\D}$, the
  composite
  \begin{equation}\label{eq:dhks-reedy-deform}
    \xymatrix{
      \M^\D \ar[r]^{T^*} &
      \M^{\Delta\D} \ar[r]^{Q} &
      \M^{\Delta\D} \ar[r]^{\lan_T} &
      \M^\D
    }
  \end{equation}
  is a left deformation for the functor $\colim_\D\maps \M^\D\to\M$,
  and therefore gives rise to a left approximation and left derived
  functor of it.  This is called the \emph{homotopy colimit}.
  Homotopy limits are constructed dually.
\end{thm}
\begin{proof}[Sketch of Proof]
  See~\cite[ch.~IV]{dhks} for the full proof.  The fact that the
  composite
  \[\colim_\D \lan_T Q T^*,\]
  is homotopical follows directly from \autoref{simplices-reedy}.  One
  can also use \autoref{simplices-reedy} to show that $\lan_T$ is
  homotopical on Reedy cofibrant diagrams, and hence the
  composite~(\ref{eq:dhks-reedy-deform}) is also homotopical.  The
  trickier parts are constructing a natural weak equivalence from the
  composite~(\ref{eq:dhks-reedy-deform}) to the identity functor of
  $\M^\D$ and showing that $\colim_\D$ in fact preserves \emph{all}
  weak equivalences between diagrams of the form $\lan_T F'$, where
  $F'$ is a Reedy cofibrant $\Delta\D$-diagram.  These proofs can be
  found in~\cite[\S23]{dhks}.
\end{proof}


A very similar procedure is followed in~\cite{chacholski_scherer},
using essentially the Reedy model structure on $\M^{\Delta\D}$, but
restricted to the subcategory of ``bounded'' diagrams, which are
easier to analyze explicitly.

These technical results give homotopy limit and colimit functors for
any model category (and, moreover, all homotopy Kan extensions, which
fit together into a ``homotopy limit system'' in the terminology
of~\cite{dhks}).  So we have a complete solution for the global
approach, at least in the unenriched context.  But for computations
and examples, as well as for dealing with enrichment, the local
approach is much more tractable; hence the need for a comparison.

\section{Local definitions and tensor products of functors}
\label{sec:local-hoco}


Recall that rather than ask for homotopical properties of the homotopy
limit \emph{functor}, the local approach asks the homotopy limit
\emph{object} to represent something homotopically coherent.  In the
case of colimits, instead of gluing the objects together directly as
in the usual colimit, we glue them together with homotopies in
between.

To have a precise notion of ``homotopy,'' we focus on
\emph{simplicially enriched categories} $\M$.  These are categories
with simplicial sets of morphisms $\M(M,M')$ and simplicial
composition maps $\M(M',M'')\times\M(M,M')\to\M(M,M'')$ which are
unital and associative; we think of the $0$-simplices of $\M(M,M')$ as
the maps $M\to M'$ and the higher simplices as homotopies and higher
homotopies between them.

\begin{rem}\label{simplicial-category}
  Many authors call a simplicially enriched category a ``simplicial
  category.''  However, there is an unfortunate collision with the use
  of the term ``simplicial object in $\mathscr{X}$'' to denote a
  functor $\DD\op\to\mathscr{X}$.  A small simplicially enriched
  category can be identified with a simplicial object in \textbf{Cat}
  which is ``object-discrete.''
\end{rem}

In order to be able to glue these homotopies together in \M, we need
to be able to ``multiply'' objects of \M\ by simplicial sets, as in
the following definition.

\begin{defn}
  A simplicially enriched category \M\ is \emph{tensored} if we have
  for every simplicial set $K$ and object $M$ of \M, an object $K\odot
  M$ of \M, together with an isomorphism of simplicial sets
  \begin{equation}\label{eq:tensor}
    \M(K\odot M, M') \iso \map(K, \M(M,M'))
  \end{equation}
  where $\map$ is the usual simplicial mapping space.  Tensors are
  unique up to unique isomorphism, when they exist, and can be made
  functorial in a unique way such that~(\ref{eq:tensor}) becomes a
  natural isomorphism.
\end{defn}

For example, \sS\ is simplicially enriched by $\map$, and tensored by
the ordinary cartesian product.  The category $\sS_*$ of based
simplicial sets is enriched by the pointed mapping space and tensored
with $K\odot X = K_+ \smsh X$.  The category of (compactly generated)
topological spaces is simplicially enriched by taking the singular
complex of the mapping space, and tensored with $K\odot X = |K|\times
X$.  Similarly, most categories of spectra are simplicially enriched
and tensored.

If \M\ is a tensored simplicially enriched category, there is a
well-known construction called the \emph{tensor product of functors},
which assigns an object
\begin{equation}\label{eq:unenriched-tensor-product}
  G\odot_\D F = 
  \coeq\left( \coprod_{d\to[u] d'} Gd'\odot Fd
  \rightrightarrows \coprod_{d} Gd\odot Fd \right)
\end{equation}
of \M\ to a given pair of functors $F\maps \D\to\M$ and $G\maps
\D\op\to\sS$.  In favorable cases, this is an example of a
\emph{weighted colimit}.  We will study homotopy weighted colimits in
the second part of this paper, but for now we use tensor products of
functors as a tool to construct ordinary homotopy colimits.  Note that
if we take $G$ to be the functor constant at the terminal simplicial
set $*$, which has a unique simplex in each dimension, then we obtain
\begin{equation}\label{eq:colim-as-tensor-product}
  *\odot_\D F \iso \colim F;
\end{equation}
thus the tensor product $G\odot_\D F$ can be thought of as the colimit
of $F$ ``fattened up'' by $G$.

Many of the the usual constructions of homotopy theory, such as
homotopy pushouts, cones, mapping telescopes, and so on, can be
expressed using such tensor products.  Inspecting these examples, all
of which have some claim to be a ``homotopy'' version of some colimit,
we notice that in most cases the shape of the homotopies we glue in
can be described by the nerves of the under-categories $(d\dn\D)$, for
objects $d$ of \D.  We have a functor $N(-\dn\D)\maps  \D\op\to\sS$, and so
we make the following definition.

\begin{defn}\label{def:loc-uhc}
  The \emph{uncorrected homotopy colimit} of a diagram $F\maps
  \D\to\M$ is the tensor product of functors:
  \begin{equation}\label{eq:loc-uhoco}
    \uhc F = N(-\dn\D)\odot_\D F
  \end{equation}
  if it exists.
\end{defn}
The ``corrected'' version will be given in \autoref{def:loc-cor}.  If the
tensor products~(\ref{eq:loc-uhoco}) exist for all $F\maps \D\to\M$,
they define an \emph{uncorrected homotopy colimit functor} $\uhc\maps
\M^\D\to\M$.

\begin{rem}\label{rmk:hococones}
  Just as an ordinary colimit object is a representing object for
  cones under a diagram, the uncorrected homotopy colimit is a
  representing object for \emph{homotopy coherent cones}.  This
  ``homotopical universal property'' will be made precise in
  \S\ref{sec:coherent-transf}; see~\cite{hcct} for a similar approach.
\end{rem}

Before moving on, we briefly mention the dual notions.  A simplicially
enriched category \M\ is \emph{cotensored} if we have for every
simplicial set $K$ and object $M$ of \M, an object $\cten{K,M}$ of \M\ 
and an isomorphism of simplicial sets
\begin{equation*}
  \M(M, \cten{K,M'}) \iso \map(K, \M(M,M')).
\end{equation*}
For example, \sS\ is cotensored over itself by $\map$.  Many
simplicially enriched categories, including based simplicial sets,
(compactly generated) topological spaces, and spectra, are also
cotensored.  Like tensors, cotensors are unique and functorial when
they exist.

If \M\ is a cotensored simplicially enriched category and we have
functors $F\maps \D\to\M$ and $G\maps \D\to\sS$, we define the
\emph{cotensor product of functors} $\cten{G,F}^\D$ using equalizers
and products in a completely dual way.  Finally, the \emph{uncorrected
  homotopy limit} is defined to be the cotensor product of functors
\begin{equation}\label{eq:loc-uholim}
  \uhl F = \cten{N(\D\dn-),F}^\D.
\end{equation}
When all these cotensor products exist, they define a functor
$\uhl\maps \M^\D\to\M$.

\section{The bar construction}
\label{sec:bar}

The ``generalized bar construction'' is quite unreasonably useful in
mathematics; here we apply it to rephrase the definition of homotopy
colimits in a computationally useful way, which we will show in
\S\ref{sec:coherent-transf} can also be conceptually helpful.  There
are many ways to approach the bar construction, and different readers
will find different motivations more natural.  The introduction given
here is chosen to clarify the connection with ordinary colimits and
explain why the bar construction is a natural thing to consider in the
context of homotopy colimits.  Another motivation can be found in
\S\ref{sec:coherent-transf}.

Refer to the diagram~(\ref{eq:simpl-diag}).  The functor
$S^*\maps \M^\D\to\M^{\Delta\op\D}$ is full and faithful, so that the
counit $\lan_S S^*\to \Id$ of the adjunction $\lan_S\adj S^*$ is an
isomorphism.  Hence by the universal properties of Kan extensions
(see~\cite[X.4~exercise~3]{maclane}), for any $F\in\M^\D$ we have:
\begin{equation}\label{eq:bar-colims}
  \begin{split}
    \colim_\D F &\iso \colim_\D \lan_S S^* F\\
    &\iso \colim_{\Delta\op\D} S^* F\\
    &\iso \colim_{\Delta\op} \lan_\Sigma S^* F.
  \end{split}
\end{equation}
Thus the colimit of a diagram $F$ of any shape is isomorphic to the
colimit of a diagram $\lan_\Sigma S^* F$ of shape $\DD\op$.  This
simplicial object has been studied for many years under another name.
\begin{defn}
  The \emph{simplicial bar construction} of $F$ is a simplicial object
  of \M, denoted $\sB(*,\D,F)$, whose object of $n$-simplices is
  \begin{equation}
    B_n(*,\D,F) = \coprod_{\alpha\maps [n]\to\D} F(\alpha_0)\label{eq:bar}
  \end{equation}
  and whose faces and degeneracies are defined using composition in
  \D, the evaluation maps $\D(d,d')\times F(d)\to F(d')$, and
  insertion of identity arrows.
\end{defn}
We refer the reader to~\cite[\S9-10]{may:goils}
and~\cite[\S7~and~\S12]{may:csf} for May's original definitions, and
to~\cite{meyer_i} and~\cite{meyer_ii} for a more abstract point of
view.  These references deal with the two-sided version, which will be
introduced in \S\ref{sec:comp-bar}.  The term ``bar construction''
comes from the use of bars for tensor products in the pre-\TeX\ days
of homological algebra, when the symbol $\ten$ was difficult to
typeset.  The symbol $*$ on the left is just a placeholder for now; we
will see what it represents in \S\ref{sec:comp-bar}.

\begin{rem}\label{bar-nerve}
  Note that when $\M=\mathbf{Set}$, the bar construction is just the
  nerve of the category of elements of $F$.  In particular, we have
  $B(*,\D,*) \iso N\D$, where here $*$ denotes the functor constant at
  a one-element set.  In other ``set-like'' cases (such as topological
  spaces), $\sB(*,\D,F)$ is the nerve of an internal ``enriched
  category of elements,'' which is sometimes called the
  \emph{transport category}.  However, for general \M\ the bar
  construction need not be the nerve of anything.
\end{rem}

Simplicial bar constructions appear ubiquitously in mathematics.  One
clue as to the reason is the following characterization in terms of
Kan extensions: they can be regarded as a way to replace a diagram of
any shape by a ``homotopically equivalent'' simplicial diagram.
Another clue is \autoref{unenriched-sb-coft}, which will tell us why this
is often such a useful replacement to make.  Yet another will be
explained in \S\ref{sec:coherent-transf}: they provide a construction
of ``homotopy coherent natural transformations.''  Finally, we will
see in \S\ref{sec:enrich-htpy-tens} that the bar construction is the
``universal deformation'' for diagram categories, which is perhaps the
most satisfying explanation of all.

\begin{prop}\label{bar=may}
  $\lan_\Sigma S^* F \iso \sB(*,\D,F)$.
\end{prop}
\begin{proof}
  By~\cite[X.3 theorem 1]{maclane}, we have the following formula for
  computing the left Kan extension.
  \begin{equation}\label{eq:lan-fmla}
    (\lan_{\Sigma} S^* F)(n) \iso \colim_{(\Sigma\dn [n])} S^* F.
  \end{equation}
  However, the comma-category $(\Sigma\dn [n])$ is a disjoint union of
  categories with a terminal object.  These terminal objects are the
  simplices $\alpha\maps [n]\to\D$ in $\Delta\op\D$, and hence the above
  Kan extension is a coproduct of $S^*F(\alpha) = F(\alpha_0)$ over
  all such $\alpha$.  This provides exactly the desired natural
  isomorphism
  \begin{equation*}
    \big(\lan_\Sigma S^* F\big)(n) \iso B_n(*,\D,F).
    \qed
  \end{equation*}
  \renewcommand{\qed}{}
\end{proof}

Recall from~(\ref{eq:bar-colims}) that the colimit of $\lan_\Sigma S^*
F$ is isomorphic to the colimit of $F$.  But since the colimit of a
simplicial diagram is just the coequalizer of the two face maps
$d_0,d_1\maps  X_1\rightrightarrows X_0$, we
recover the classical formula for computing colimits using coproducts
and coequalizers:
\begin{equation}\label{eq:compute-colim}
  \colim_\D F \iso \coeq\left(\coprod_{d\to d'} Fd \rightrightarrows
  \coprod_d Fd\right).
\end{equation}
Dually, of course, the limit of any functor is isomorphic to the limit
of the \emph{cosimplicial cobar construction}, which has a
corresponding characterization using the left half
of~(\ref{eq:simpl-diag}), and this reduces to the analogue
of~(\ref{eq:compute-colim}) for computing limits using products and
equalizers.

The above shows that for usual limits and colimits, the ``higher
information'' contained in the bar construction above level one is
redundant.  However, in the homotopical situation, this information
plays an important role and we have to incorporate it.


In the classical context of topological spaces, the ``correct homotopy
type'' associated to a simplicial space (or spectrum) is usually its
geometric realization.  We can define the \emph{geometric realization}
of a simplicial object $\sX$ in any tensored simplicially enriched
category as the tensor product of functors
\begin{equation*}
  |\sX| = \Delta \odot_{\DD\op} \sX
\end{equation*}
where $\Delta\maps \DD\to\sS$ is the canonical cosimplicial simplicial
set.  When $\M$ is the category of topological spaces, this reduces to
the usual geometric realization.  When $\M=\sS$, the realization of a
bisimplicial set is isomorphic to its diagonal.


Then we can make the following definition.
\begin{defn}
  The \emph{bar construction} is the geometric realization of the
  simplicial bar construction:
  \begin{align*}
    B(*,\D,F) &= |\sB(*,\D,F)|\\
    &\equiv \Delta\odot_{\DD\op} \sB(*,\D,F)
  \end{align*}
\end{defn}

Before stating the next lemma, we need to clear up some notation.  For
any object $d$ of \D, we have the hom-functor $\D(d,-)\maps
\D\to\mathbf{Set}$, which we can regard as a functor $\D\to\sS$
landing in discrete simplicial sets.  Since \sS\ is tensored over
itself, we can form the bar construction $\sB(*,\D,\D(d,-))$, which by
\autoref{bar-nerve} is the nerve of the category of elements of
$\D(d,-)$, viewed as a horizontally discrete bisimplicial set.  Hence
its diagonal $B(*,\D,\D(d,-))$ is precisely that nerve.  Varying $d$,
we obtain a functor $\D\op\to\sS$ which we denote $B(*,\D,\D)$.
\begin{lem}\label{one-sided-bar-pullout}
  \begin{equation*}
    B(*,\D,F) \iso B(*,\D,\D)\odot_\D F.
  \end{equation*}
\end{lem}
\begin{proof}
  Because tensor products of functors, including geometric
  realization, are all colimits, they all commute with each other.
  Therefore the right-hand side can be written as $B(*,\D,\D\odot_\D
  F)$, where $\D\odot_\D F\maps \D\to\M$ is defined by $(\D\odot_\D
  F)(d) = \D(-,d)\odot_\D F$.  But writing out the definitions, it is
  easy to see that $\D\odot_\D F \iso F$.
\end{proof}

The punchline is that the bar construction is precisely the
uncorrected homotopy colimit from~\S\ref{sec:local-hoco}.
\begin{prop}\label{prop:loc=bar}
  For $F\maps \D\to\M$ where \M\ is simplicially enriched, tensored,
  and cocomplete,
  \begin{align*}
    \uhc F
    &= N(-\dn\D)\odot_\D F\\
    &\iso B(*,\D,F).
  \end{align*}
\end{prop}
\begin{proof}
  Recall that $B(*,\D,\D)$ is the nerve of the category of elements of
  $\D(d,-)$.  But that category of elements is precisely the
  comma-category $(d\dn\D)$; thus we have a natural isomorphism
  \[\sB(*,\D,\D) \iso N(-\dn\D).\]
  Now applying \autoref{one-sided-bar-pullout}, we obtain
  \begin{equation*}
    \begin{split}
      B(*,\D,F) &\iso B(*,\D,\D) \odot_\D F\\
      &\iso N(-\dn\D) \odot_\D F = \uhc F.
      \qed
    \end{split}
  \end{equation*}
  \renewcommand{\qed}{}
\end{proof}
So far, this is just a rephrasing of \autoref{def:loc-uhc}, and was known
as far back as~\cite{bk}, who called the simplicial bar construction
the \emph{simplicial replacement}.  However, as we will see
in~\S\ref{sec:comp-bar}, the bar construction provides an alternate
approach to the comparison of local and global homotopy colimits.

\section{A comparison based on Reedy model structures}
\label{sec:comp-dhks}

In order to compare the local and global definitions, we need to put
ourselves in a context in which both are defined, and we also require
some compatibility conditions.  So suppose now that \M\ is a
simplicial model category.  A simplicial model category is a model
category which is also simplicially enriched, tensored, and
cotensored, and such that the lifting properties in the definition of
model category are extended to ``homotopy lifting properties'' with
respect to the simplicial mapping spaces.
See~\cite[ch.~9]{hirschhorn} for a definition.

One of the consequences of the definition is the following result,
which will be crucial later on for analyzing bar constructions.
\begin{lem}\label{realiz-quillen}
  If \M\ is a simplicial model category, then the geometric
  realization functor $|\cdot|\maps \M^{\DD\op}\to\M$ is left Quillen
  when $\M^{\DD\op}$ is given the Reedy model structure.  Therefore,
  it preserves weak equivalences between Reedy cofibrant simplicial
  objects.
\end{lem}
\begin{proof}
  This is a special case of~\cite[18.4.11]{hirschhorn} (which we state
  below as \autoref{simplicial-reedy-2va}).  The consequence is stated
  as~\cite[18.6.6]{hirschhorn}.
\end{proof}

Now, since \autoref{def:loc-uhc} uses no properties of the model
structure, it is unsurprising that we have to modify it a bit.  Recall
that we assume all model categories have functorial factorizations, so
that \M\ has fibrant and cofibrant replacement functors, which we call
$R$ and $Q$ respectively.
\begin{defn}\label{def:loc-cor}
  The \emph{corrected homotopy limit} and \emph{colimit} of $F$ are:
  \begin{align}
    \hocolim F &= \uhc QF\\
    \holim F &= \uhl RF
  \end{align}
  where $Q$ and $R$ are applied objectwise to $F$.
\end{defn}

The approach to analyzing these functors in~\cite[\S18.4 and
18.5]{hirschhorn} goes as follows.  Because $N(-\dn\D)$ is cofibrant
in the projective model structure on $\sS^{\D\op}$, the functor $\uhc$
preserves weak equivalences between objectwise cofibrant diagrams (and
dually for $\uhl$).  It follows that the functors $\hocolim$ and $\holim$ are
homotopical.


We will not go into this argument in detail, because while this
approach shows that $\hocolim$ is homotopical, it is unclear from the
argument whether it is defined by applying a deformation, and hence
whether it is in fact a derived functor in the sense of
\S\ref{sec:global-derived}.  The results in this section and the next
give alternate proofs of the fact that $\hocolim$ is homotopical which also
show this stronger result.

In this section we outline a proof which appears in an early draft
of~\cite{dhks} (but not in the published book).  The full proof will
be given in Appendix~\ref{sec:ugly-proof}.  I would like to thank Phil
Hirschhorn for providing me with a copy of the early draft and helping
me understand the proof.

Recall from \autoref{thm:dhks} that any cofibrant replacement functor in
the Reedy model structure on $\M^{\Delta\D}$ gives rise to a derived
functor of $\colim$.  The goal of this approach is to show that $\hocolim$
is in fact the colimit of a specific Reedy cofibrant replacement for
$F$ in $\M^{\Delta\D}$.  Define a functor $\Q\maps
\M^\D\to\M^{\Delta\D}$ as follows.
\begin{align*}
  \Q F(\alpha) &= \hocolim \alpha^* F\\
  &\iso N(-\dn[n])  \odot_{[n]} \alpha^* Q F.
\end{align*}
Here we consider an $n$-simplex $\alpha$ as a functor $\alpha\maps
[n]\to\D$.  The stated isomorphism follows from the fact that
$\alpha^* QF = Q\alpha^* F$, since $Q$ is applied objectwise.

Thus the values of $\Q F$ are ``generalized mapping cylinders'' of $F$
over all the simplices of \D.  The mapping cylinder of a map $f$ is
naturally homotopy equivalent to the codomain of $f$, and the same is
true for these generalized mapping cylinders.  Specifically, the maps
$N(-\dn [n])\to *$ induce maps
\begin{equation*}
  \xymatrix{
    N(-\dn [n])  \odot_{[n]} \alpha^* Q F \ar[r]\ar[drr]_{\delta} &
    {*} \odot_{[n]} \alpha^* Q F \ar[r]^(0.55){\iso} &
    QF(\alpha_n) \ar[d]^{q}\\
    & & F(\alpha_n) 
  }
\end{equation*}
which fit together into a natural transformation
\begin{equation*}
  \delta\maps  \Q F \too T^* F.
\end{equation*}

We then prove the following theorem.

\newsavebox{\dhkstheorem}
\begin{lrbox}{\dhkstheorem}
\begin{minipage}{\linewidth}
\begin{thm}\label{dhks-comp}
  The diagram $\Q F$ is Reedy cofibrant for all $F$, the functors $\Q$
  and $\lan_T \Q$ are homotopical, and $\delta$ is a natural weak
  equivalence.  Moreover,
  \begin{equation}\label{eq:dhks-repl}
    \begin{split}
      \hocolim F &\iso \colim_{\Delta\D}\Q F\\
      &\iso \colim_\D \lan_T \Q F.
    \end{split}
  \end{equation}
  Therefore, $(\lan_T\Q,\lan_T\delta)$ is a left deformation of
  $\colim$, and hence $\hocolim$ is a derived functor of the usual colimit.
\end{thm}
\end{minipage}
\end{lrbox}

\medskip\noindent\usebox{\dhkstheorem}

\begin{proof}
  It is straightforward to show that $\delta$ is a weak equivalence.
  To wit, the diagram $N(-\dn[n])$ is always projective-cofibrant, and
  because $[n]$ has a terminal object $n$, the constant \sS-diagram at
  a point is also projective-cofibrant.  Thus
  by~\cite[18.4.5]{hirschhorn}, this implies that the functor $-
  \odot_\D QF$ preserves the weak equivalence $N(-\dn [n])\to *$, so
  upon composing with the weak equivalence $q$, we see that $\delta$
  is a weak equivalence.
  
  It follows from the 2-out-of-3 property that \Q\ and $\lan_T\Q$ are
  homotopical.  Thus the bulk of the proof consists of
  proving~(\ref{eq:dhks-repl}) and showing that $\Q F$ is Reedy
  cofibrant.  This involves a lengthy manipulation of Kan extensions
  and tensor products of functors, which we defer to
  Appendix~\ref{sec:ugly-proof}.
\end{proof}

We do, however, note the following fact, which is crucial to the proof
that $\Q F$ is Reedy cofibrant, and will be needed in the next
section.
\begin{lem}\label{objw-coft->reedy}
  If $F\maps \D\to\M$ is objectwise cofibrant, then $S^* F$ is Reedy
  cofibrant in $\M^{\Delta\op\D}$.
\end{lem}
\begin{proof}
  Recall from~\cite{hovey} or~\cite{hirschhorn} that a diagram $X$ on
  a Reedy category is Reedy cofibrant when all the ``latching maps''
  $L_\alpha X\to X_\alpha$ are cofibrations.  Here $L_\alpha X$ is the
  ``latching object'' for $X$ at $\alpha$.

  In our situation, let $\alpha$ be an $n$-simplex in \D.  We have
  $S^*F(\alpha) = F(\alpha_0)$.  If $\alpha$ contains no identity
  maps, then the latching object $L_\alpha S^* F$ is the initial
  object.  Otherwise, it is simply $F(\alpha_0)$.  Since $F$ is
  objectwise cofibrant, in either case the map $L_\alpha S^* F \to
  S^*F(\alpha)$ is a cofibration, as required.
\end{proof}

Recall that $S^*F$ is a character who has already figured in our
discussion of the bar construction in \S\ref{sec:bar}.  Thus, this
result suggests we may be able to use the bar construction and deal
directly with objectwise cofibrant diagrams, without needing to invoke
the replacement category $\Delta\D$.  In the next section, we will use
this idea to give an alternate proof that $\hocolim$ is a derived functor
of the usual colimit.

\section{A comparison using the two-sided bar construction}
\label{sec:comp-bar}

Let \M\ be a simplicial model category.  We shall give a direct proof
of the following theorem (although some final technical steps will be
postponed to \S\ref{sec:htpy-ten}).

\begin{thm}\label{2bar-comp}
  If $Q$ is any cofibrant replacement functor on \M, then
  \[\hocolim \iso B(*,\D,Q-)\]
  is a left derived functor of $\colim$.
\end{thm}

We stress that $Q$ is only assumed to be a cofibrant replacement
functor for the model structure on $\M$, applied objectwise; \emph{no
  model structure on $\M^\D$ is required}.

We begin by showing that $\hocolim$ is homotopical.
\begin{lem}\label{unenriched-sb-coft}
  If $F\maps \D\to\M$ is objectwise cofibrant, then the simplicial bar
  construction $\sB(*,\D,F)$ is Reedy cofibrant.  Hence the functor
  \[B(*,\D,-)\maps  \M^\D \too \M\]
  preserves weak equivalences between objectwise cofibrant diagrams.
\end{lem}
\begin{proof}
  Note that when $\M^{\DD\op}$ and $\M^{\Delta\op\D}$ are given the
  Reedy model structures, the adjunction $\lan_\Sigma\adj\Sigma^*$ is
  Quillen.  In fact, $\Sigma^*$ preserves matching objects---that is
  to say, $M_\alpha \Sigma^* G \iso M_{\Sigma\alpha} G$---from which
  it follows trivially that $\Sigma^*$ preserves Reedy fibrations and
  Reedy trivial fibrations.  But by \autoref{objw-coft->reedy}, $S^*F$ is
  Reedy cofibrant whenever $F$ is objectwise cofibrant, so
  $\lan_\Sigma S^* F$ is also Reedy cofibrant.  The second statement
  follows from the first together with \autoref{realiz-quillen}.
\end{proof}

A more general version of this result will be proven in
\S\ref{sec:cofibrancy}.

Results of this sort are as old as the bar construction; it has always
been thought of as a sort of cofibrant replacement.
May~(\cite{may:goils} and~\cite{may:csf}) shows that under suitable
conditions, simplicial bar constructions (on topological spaces) are
``proper simplicial spaces,'' which is equivalent to being ``Reedy
$h$-cofibrant.''  We will have more to say about this in
\S\ref{sec:cofibrancy}.  For $\M=\sS$, this result is well-known in
simplicial homotopy theory, especially in its dual formulation:
see~\cite[Lemma~VIII.2.1]{goerss_jardine}, the classical~\cite{bk}, or
more recent~\cite{hcct}.

To understand its relevance here, recall from \autoref{prop:loc=bar} that
\[\uhc F \iso B(*,\D,F).\]
Thus, as promised after \autoref{def:loc-cor},
we have an alternate proof of the fact that the functor that takes $F$
to $\hocolim F \iso B(*,\D,QF)$ is homotopical:
\autoref{unenriched-sb-coft} tells us that $B(*,\D,-)$ preserves weak
equivalences between objectwise cofibrant diagrams, but $Q$ (applied
objectwise) is a homotopical functor replacing every diagram by an
objectwise cofibrant one.

Now we want to use the bar construction to show that $\hocolim$ is in
fact a derived functor, but first we need to introduce a two-sided
generalization of the bar construction.  Just as we replaced $*$ by a
functor $G$ to generalize ordinary colimits to tensor products of
functors in~(\ref{eq:unenriched-tensor-product}), here we do the same
for the bar construction.

\begin{defn}[{\cite[\S12]{may:csf}}]\label{def:two-sided-bar}
  Given functors $F\maps \D\to\M$ and $G\maps \D\op\to\sS$, the
  \emph{two-sided simplicial bar construction} $\sB(G,\D,F)$ is a
  simplicial object of $\M$ whose object of $n$-simplices is
  \begin{equation*}
    B_n(G,\D,F) = \coprod_{\alpha\maps [n]\to\D} G(\alpha_n)\odot F(\alpha_0).
  \end{equation*}
  and whose faces and degeneracies are defined using composition in
  \D, the evaluation maps $\D(d,d')\times F(d)\to F(d')$ and
  $G(d')\times \D(d,d')\to G(d)$, and insertion of identity arrows.
\end{defn}

As with the one-sided version, the two-sided bar construction can be
rephrased using the functors of~(\ref{eq:simpl-diag}).  If we write
$G\xodot F\maps \D\op\times\D\to\M$ for the ``external tensor
product'' of functors:
\[(G\xodot F)(d',d) = Gd'\odot Fd\]
then we have the following analogue of \autoref{bar=may}:
\begin{equation}\label{eq:2bar=may}
  \sB(G,\D,F) \iso \lan_\Sigma (T,S)^* (G\xodot F)
\end{equation}
and the following
analogues of~(\ref{eq:bar-colims}) and~(\ref{eq:compute-colim}):
\begin{equation}\label{eq:compute-wgt-colim}
  \begin{split}
    G\odot_\D F
    &= \coeq\left(
      \coprod_{d\to[u] d'} Gd'\odot Fd
      \rightrightarrows \coprod_d Gd\odot Fd
    \right)\\
    &\iso \colim_{\DD\op} \lan_\Sigma (T,S)^* G\xodot F.
  \end{split}
\end{equation}

As in the one-sided case, we define the \emph{two-sided bar
  construction} to be the realization of the two-sided simplicial bar
construction:
\begin{equation}\label{eq:2bar}
  B(G,\D,F) = |\sB(G,\D,F)|.
\end{equation}

The comparison hinges on the following standard result about bar
constructions.
\begin{lem}\label{unenriched-bar-pullout}
  \begin{equation*}
    \begin{split}
      \colim B(\D,\D,F) &\iso
      * \odot_\D B(\D,\D,F)\\
      &\iso B(*,\D,F)
    \end{split}
  \end{equation*}
\end{lem}
\begin{proof}
  The first isomorphism was stated
  as~(\ref{eq:colim-as-tensor-product}), and follows directly from the
  definitions.  Now, as in the proof of \autoref{one-sided-bar-pullout},
  because tensor products of functors, including geometric
  realization, are all colimits, they all commute with each other.
  Thus, we have $*\odot_\D B(\D,\D, F) \iso B(*\odot_\D \D, \D,F)$.
  But it is easy to see that $*\odot_\D \D \iso *$.
\end{proof}

\autoref{unenriched-bar-pullout} suggests we can take $B(\D,\D,Q-)$ as a
deformation for the colimit, but before showing that this works, we
must recall some facts about simplicial homotopy.  Given two maps
$\sX\rightrightarrows Y_\bullet$ of simplicial objects in any category
\M, there is an old, purely combinatorial, notion of simplicial
homotopy between them; see, for example,~\cite[\S5]{may:simp-obj}.
When \M\ is complete and cocomplete, we can reformulate this in modern
terms by saying that $\M^{\DD\op}$ is simplicially enriched, tensored,
and cotensored, and thus has an attendant notion of simplicial
homotopy using simplicial cylinders $\Delta^1\odot \sX$.  The tensor
has a simple description as $(K\odot X)_n = K_n \times X_n$, and the
combinatorial description follows by inspection.

For general \M\ this combinatorial notion of homotopy in $\M^{\DD\op}$
is \emph{a priori} unrelated to any homotopy theory that \M\ might
have.  However, when \M\ is tensored and simplicially enriched, we
have
\begin{equation}\label{eq:realiz-pres-tensor}
  |K\odot \sX| \iso K\odot |\sX|
\end{equation}
for any simplicial set $K$ and simplicial object $\sX$ in \M.  See,
for example,~\cite[5.4]{rss}.  It follows that ``geometric realization
preserves simplicial homotopy'' in the sense that a simplicial
homotopy between maps $\sX\rightrightarrows Y_\bullet$ gives rise to a
simplicial homotopy between the realizations $|\sX|\rightrightarrows
|Y_\bullet|$ of these maps.

In particular, simplicial homotopy \emph{equivalences} in
$\M^{\DD\op}$, in the above purely combinatorial sense, realize to
simplicial homotopy equivalences in \M.  Moreover, if \M\ is a
simplicial \emph{model} category, simplicial homotopy equivalences are
necessarily also weak equivalences for the model structure---see, for
example~\cite[9.5.16]{hirschhorn}.

With this as background, we can now continue on the path to relating
the two-sided bar construction to homotopy colimits, by showing that
$B(\D,\D,Q-)$ does in fact define a deformation for the colimit.  As
with \autoref{unenriched-bar-pullout}, this is a standard result about
bar constructions which we merely rephrase in abstract language.  The
idea is found in~\cite[9.8]{may:goils} and is clarified and
generalized in~\cite[\S6-7]{meyer_i}.

\begin{lem}\label{unenriched-bar-repl}
  There is a natural weak equivalence $\ep\maps B(\D,\D,F) \we F$.
  Thus the functor $B(\D,\D,Q-)\maps  \M^\D\to\M^\D$ is a deformation.
\end{lem}
\begin{proof}
  Since $F$ is the colimit of $\sB(\D,\D,F)$, we have a natural map
  \[\tep\maps  \sB(\D,\D,F)\to \tilde{F},\]
  where $\tilde{F}$ is the constant simplicial diagram on $F$.  Since
  $F$ is the realization of $\tilde{F}$, the realization of $\tep$ is
  our desired map $B(\D,\D,F)\to F$.  It remains only to show that
  $\ep$ is a weak equivalence, and by the above discussion, for this
  it suffices to show that $\tep$ is objectwise a simplicial homotopy
  equivalence.
  
  We have a map $\eta\maps  \tilde{F}\to\sB(\D,\D,F)$ given by inserting
  identities, and clearly $\tep\eta = 1_F$.  Moreover,
  $\sB(\D,\D,F)$ has an \emph{extra degeneracy} induced by inserting
  identities at the beginning of $\alpha$ (inserting identities
  elsewhere induces the usual degeneracies).  An extra degeneracy for
  a simplicial object $\sX$ is a collection of maps
  $s_{-1}\maps X_n\to X_{n+1}$ such that $d_0 s_{-1} = 1_{X_n}$,
  $d_{i+1}s_{-1} = s_{-1}d_i$, and $s_{j+1}s_{-1} = s_{-1}s_j$ for all
  $i$ and $j$.  This extra degeneracy furnishes a simplicial homotopy
  \begin{equation}\label{eq:extra-deg}
    \eta\tep \eqv 1 \maps  \sB(\D(-,d),\D,F) \to \sB(\D(-,d),\D,F)
  \end{equation}
  for each fixed $d\in\D$.  Thus each map $B(\D(-,d),\D,F)\to Fd$ is a
  simplicial homotopy equivalence in \M, hence a weak equivalence, as
  desired.  (Note, though, that the simplicial
  homotopies~(\ref{eq:extra-deg}) are not natural in $d$.)
\end{proof}


Taken together, Lemmas~\ref{unenriched-bar-pullout}
and~\ref{unenriched-bar-repl} say that $\hocolim F$ is actually the
colimit of a weakly equivalent replacement diagram $B(\D,\D,QF)$.  To
complete the proof of \autoref{2bar-comp}, it remains to show that
$B(\D,\D,Q-)$ is a left $\colim$-deformation in the sense of
\autoref{def:deform}, that is, that $\colim$ is homotopical on its full
image.  This is a special case of \autoref{thm:hoten} below, which
applies to the more general enriched case.  We do not repeat the proof
here, as it does not simplify much in this case and in any case is not
particularly enlightening.

This completes the proof of \autoref{2bar-comp}.  The reader may find
this route to the comparison theorem to be simpler than the comparison
proof sketched in \S\ref{sec:comp-dhks} and detailed in
Appendix~\ref{sec:ugly-proof} (combined with the proof
in~\cite[ch.~IV]{dhks} of our \autoref{thm:dhks}).  However, the real
merit of the approach in this section is that it applies nearly
verbatim to the more general enriched context, in which the bar
construction makes perfect sense but nerves and categories of
simplices do not.

\section{Homotopy coherent transformations}
\label{sec:coherent-transf}

There is a different approach we might take for the analysis of the
homotopy theory of diagrams, including homotopy limits and colimits.
Recall that weak equivalences are frequently handled in homotopy
theory by replacing the objects with ``good'' objects such that weak
equivalences become \emph{homotopy} equivalences, which already have
inverses up to homotopy.  Thus, suppose that \M\ has this (vaguely
defined) property, and consider how we might show that $\M^\D$ also
has it.

One approach is to imitate the theory of ``cell complexes'' used in
classical homotopy theory, constructing ``cell diagrams'' using
diagrams freely generated by cells placed at different positions.
This idea leads to the projective model structure, discussed briefly
in \S\ref{sec:global-hoco}, with all the attendant disadvantages.  In
this section we consider a different approach.

Let $F$ and $G$ be diagrams in \M\ of shape \D, and suppose that we
have a natural weak equivalence $\alpha\maps  G\to F$.  We can easily
replace $F$ and $G$ objectwise by ``good'' objects, so that each
component $\alpha_d\maps Gd\to Fd$ has a homotopy inverse, but in general
these homotopy inverses will not fit together into a natural
transformation.  However, we can reasonably expect that they will fit
together into a \emph{homotopy coherent transformation}.  Intuitively,
this means that the squares
\begin{equation}\label{eq:coh-trans}
  \xymatrix{
    Gd_0 \ar[r]^{Gf} & Gd_1 \ar@{}[dl]|{\eqv}\\
    Fd_0 \ar[r]_{Ff} \ar[u]^{\alpha_{d_0}^{-1}} & Fd_1 \ar[u]_{\alpha_{d_1}^{-1}}
  }
\end{equation}
commute only up to homotopy; that for a pair of composable arrows
$d_0\to[f] d_1\to[g] d_2$ in \D, the corresponding homotopies are
related by a higher homotopy; and so on.

The notion of homotopy coherent transformation has been made precise
in many places, such as~\cite{hcct}, where it is defined using what is
essentially the cobar construction.  In our terminology, the
definition is as follows.  We assume that $\M$ is tensored and
cocomplete, so that the bar construction can be defined.

\begin{defn}
  A \emph{homotopy coherent transformation} $F\cohto G$ is a natural
  transformation $B(\D,\D,F)\to G$.
\end{defn}

To justify this definition, we unravel the first few levels.  Note
that by definition of geometric resolution, such a natural
transformation is given by a collection of $n$-simplices
$B_n(\D,\D,F)\to G$ for each $n$, related by face and degeneracy maps
in an appropriate way.  For small values of $n$, what this translates
to is the following.

\begin{itemize}
\item For each arrow $d_0\to[f] d_1$ in \D, we have a morphism
  $\alpha_f\maps Fd_0\to Gd_1$.
\item These morphisms have the property that for any $d_0\to[f] d_1
  \to[g] d_2$ in \D, we have $\alpha_{gf} = G(g)\circ \alpha_f$.  This
  comes from the naturality of the transformation $B(\D,\D,F)\to G$.
  In particular, we have $\alpha_f = G(f)\circ \alpha_{1_{d_0}}$, so
  the 0-simplices are determined by morphisms $\alpha_d\maps  Fd\to Gd$.
\item For any $d_0\to[f] d_1 \to[g] d_2$ in \D, we have a homotopy
  $\alpha_{gf} \eqv \alpha_g\circ F(f)$.  In particular, for any
  $d_0\to[f] d_1$ we have $G(f)\circ \alpha_{d_0} = \alpha_f \eqv
  \alpha_{d_1}\circ F(f)$, so the squares~(\ref{eq:coh-trans}) commute
  up to homotopy.
\item The above homotopies are preserved by composition with the
  action of $G$ on arrows, in a straightforward way.
\item For any $d_0\to[f] d_1\to[g] d_2\to[h] d_3$ in \D, there is a
  higher homotopy (a 2-simplex) comparing the homotopy $\alpha_{hgf}
  \eqv \alpha_{h} \circ F(gf)$ to the composite homotopy $\alpha_{hgf}
  \eqv \alpha_{hg} \circ F(f) \eqv \alpha_h \circ F(g) \circ F(f) =
  \alpha_h \circ F(gf)$.
\item And so on\dots
\end{itemize}

The above remarks lead us to expect that we could produce a ``good''
replacement for a diagram by replacing its values objectwise by
``good'' objects in \M, and then applying the bar construction
$B(\D,\D,-)$.  The results of \S\ref{sec:comp-bar} show that, at least
when considering the colimit functor, this strategy works.  In
\S\ref{sec:enrich-htpy-tens} we will see that it also works to
construct the homotopy category of a diagram category, as we would
expect: objectwise ``homotopy equivalences'' really do have ``homotopy
coherent'' inverses.

We can now also make precise the notion of ``homotopy coherent cones''
referred to in \autoref{rmk:hococones}.  A homotopy coherent cone under a
diagram $F$ should consist of a vertex $C$ and a coherent
transformation $F\cohto C$, where by $C$ we mean the diagram constant
at the vertex $C$.  By definition, this is the same as a natural
transformation $B(\D,\D,F)\to C$.  Therefore, a representing object
for homotopy coherent cones under $F$ is precisely a colimit of
$B(\D,\D,F)$; but as we saw in \autoref{unenriched-bar-pullout}, this is
precisely the homotopy colimit of $F$.  Thus the latter really does
have a ``homotopical universal property.''

\begin{rem}
  We would like to take this opportunity to briefly mention a
  misconception which can arise.  Since ordinary colimits can be
  constructed using coproducts and coequalizers, it may seem natural
  to ask whether homotopy colimits can be constructed using homotopy
  coproducts and homotopy coequalizers.  After having read
  \S\ref{sec:bar}, the reader may be able to guess that this is not
  true: the naive version of such a construction would correspond to
  taking the geometric realization of the bar construction truncated
  at level 1, while we have seen that in general, all the higher
  information must be included to get the correct result.  In the
  language of this section, the naive construction would be a
  representing object for cones which commute up to homotopy, but not
  up to \emph{coherent} homotopy.

  There are special cases, however, in which a naive construction does
  give an equivalent result.  For instance, if the category \D\ is
  freely generated by some graph, then a cone which commutes with the
  \emph{generating} arrows up to specified homotopies can be extended,
  essentially uniquely, to a cone which commutes up to coherent
  homotopy.  In such cases we can use the naive construction, as long
  as we glue in homotopies only for the generating arrows.  A specific
  example of this is a mapping telescope, which is a homotopy colimit
  over the category freely generated by the graph
  $(\cdot\to\cdot\to\cdot\to\dots)$.
\end{rem}

This completes the second part of this paper, and we now turn to
enriched category theory and weighted homotopy colimits.


\section{Enriched categories and weighted colimits}
\label{sec:enriched}

In the majority of categories encountered by the working
mathematician, the hom-sets are not just sets but have some extra
structure.  They may be the vertices of simplicial mapping spaces, as
for the simplicially enriched categories considered in
\S\ref{sec:local-hoco}, but they may also be topological spaces, or
abelian groups, or chain complexes, or small categories.  A general
context for dealing with such situations is given by \emph{enriched
  category theory}.

The importance of enriched category theory for us is twofold.
Firstly, as remarked above, most categories which arise in practice
are enriched; thus it is essential for applications to be able to deal
with enriched situations.  Secondly, however, the presence of
enrichment often simplifies, rather than complicates, the study of
homotopy theory, since enrichment over a suitable category
automatically provides well-behaved notions of homotopy and homotopy
equivalence.  We will see examples of this in
\S\ref{sec:derived-enrichment}--\ref{sec:derived-via-enriched}.

In this section, we will review some notions of enriched category
theory.  The standard references for the subject are~\cite{kelly}
and~\cite{dubuc}.

Let $\V_0$ be a bicomplete closed symmetric monoidal category with
product $\ten$ and unit object $E$.  Often the unit object is written
$I$, but in homotopy theory there is potential for confusion with an
interval, so we prefer a different notation.  Here \emph{closed} means
that we have internal hom-objects $[B,C]$ satisfying an adjunction
\[\V_0(A\ten B,C) \iso \V_0(A, [B,C]).\]
All the above examples give possible choices of $\V_0$.  A category
\M\ \emph{enriched over \V}, or a \emph{\V-category}, is defined by a
collection of objects $M, M', \dots$ together with objects $\M(M,M')$
of \V, unit morphisms $E\to\M(M,M)$, and composition morphisms
\[\M(M',M'')\ten\M(M,M')\to\M(M,M''),\]
satisfying associativity and unit axioms.  Note that if $\V_0$ is the
category \textbf{Set} with its cartesian monoidal structure, then a
$\V$-category is the same as an ordinary category (with small
hom-sets).

Such a \V-category gives rise to an \emph{underlying} ordinary
category $\M_0$ defined to have the same objects and $\M_0(M,M') =
\V_0(E,\M(M,M'))$; experience shows that this always gives the correct
answer.  The internal-hom of $\V_0$ gives rise to a \V-category,
denoted \V, with $\V(A,B)=[A,B]$.  It is easy to check that the
underlying category of the \V-category $\V$ is, in fact, the original
category $\V_0$, justifying the notation.  We often abuse notation by
referring to the original closed symmetric monoidal category by \V\ as
well.

If \M\ and \N\ are \V-categories, a \emph{\V-functor} $F\maps \M\to\N$
consists of an assignation of an object $F(M)$ of \N\ to each object
$M$ of \M\ together with morphisms
\[\M(M,M')\to\N(F(M),F(M'))\]
in $\V_0$ that satisfy functoriality diagrams.  As with \V-categories,
a \V-functor gives rise to an underlying ordinary functor $F_0\maps
\M_0\to\N_0$.  Similarly, we have notions of a \emph{\V-adjoint pair}
and a \emph{\V-natural transformation} which give rise to underlying
ordinary adjoint pairs and natural transformations, but are in general
stronger.

Now let \D\ be a small \V-category (i.e.\ it has a small set of
objects), and $F\maps \D\to\M$ a \V-functor.  We view $F$ as a diagram
of shape \D\ in \M, and would like to take its limit.  However, for
general \V\ it may not be possible to define a ``cone'' over such a
diagram in the usual manner; thus we are forced to work with weighted
limits from the very beginning.  Our approach to weighted colimits is
via tensor products of functors, generalizing those used in
\S\ref{sec:local-hoco}.



We say that a \V-category is \emph{tensored} if we have for every
object $K$ of \V\ and $M$ of \M, an object $K\odot M$ of \M\ and
\V-natural isomorphisms:
\[\V(K,\M(M,M'))\iso \M(K\odot M,M')\]
and \emph{cotensored} if in the same situation we have an object
$\cten{K,M}$ and \V-natural isomorphisms:
\[\V(K,\M(M',M))\iso \M(M',\cten{K,M}).\]
The \V-category \V\ is always tensored and cotensored; its tensor is
the monoidal product $\ten$ and its cotensor is its internal-hom.
When $\V=\mathbf{Set}$, tensors are copowers and cotensors are powers.
As in the simplicial context, tensors are unique up to unique
isomorphism when they exist, and can be made functorial in a unique
way so that the defining isomorphisms become natural.

We can now make the following definition.

\begin{defn}\label{def:wgt-colim}
  If \M\ is a tensored \V-category, \D\ is a small \V-category, and
  $G\maps \D\op\to\V$ and $F\maps \D\to\M$ are \V-functors, their
  \emph{tensor product} is the following coequalizer:
  \begin{equation}\label{eq:wgt-colim}
    G\odot_\D F =
    \coeq\left( \coprod_{d,d'} \big(Gd'\ten \D(d,d')\big) \odot Fd
      \rightrightarrows \coprod_{d} Gd\odot Fd \right).
  \end{equation}
  When \M\ is cotensored, we define the \emph{cotensor product of
    functors} $\cten{G,F}^\D$ of $F\maps \D\to\M$ and $G\maps \D\to\V$
  to be a suitable equalizer.
\end{defn}

\begin{rem}\label{tensor-assoc}
  By the universal properties of tensors, we have the following
  ``associativity'' isomorphism:
  \begin{equation*}
    \big(Gd'\ten \D(d,d')\big) \odot Fd 
    \;\iso\;
    Gd'\odot \big(\D(d,d') \odot Fd \big).
  \end{equation*}
  This isomorphism is needed to define one of the two arrows
  in~(\ref{eq:wgt-colim}).  Compared
  to~(\ref{eq:unenriched-tensor-product}), the hom-object $\D(d,d')$ has
  moved out from under the $\coprod$, since it is no longer simply a
  set we can take a coproduct over.
\end{rem}

Of course, \autoref{def:wgt-colim} uses ordinary colimits: coproducts and
coequalizers.  The indexing categories for these colimits are not
enriched, however, so we can take these to mean ordinary colimits in
the underlying category $\M_0$.  However, generally we need to require
that their universal property is also suitably enriched, as follows.
Recall that the universal property of a colimit in $\M_0$ can be
expressed by a natural bijection of sets
\begin{equation*}
  \M_0\Big(\operatornamewithlimits{colim}_{d\in D} F(d), \, Y\Big)
  \iso \lim_{d\in D} \M_0\big(F(d),\,Y\big).
\end{equation*}
We say that this colimit is a \emph{\V-colimit} if this bijection is
induced by a \V-natural isomorphism in \V\ 
\begin{equation*}
  \M\Big(\operatornamewithlimits{colim}_{d\in D} F(d), \, Y\Big)
  \iso \lim_{d\in D} \M\big(F(d),Y\big),
\end{equation*}
where now the limit on the right is an ordinary limit in the category
\V.

When the coproducts and coequalizer in~(\ref{eq:wgt-colim}) are
\V-colimits, it is proven in~\cite[\S3.10]{kelly} that the tensor
product of functors coincides with the \emph{weighted colimit of $F$
  weighted by $G$}, as defined there.  The general definition of
weighted colimits includes tensors and \V-colimits as special cases,
but we will not need this definition, so we refer the reader
to~\cite[\S3]{kelly} for its development.  A \V-category \M\ is said
to be \emph{cocomplete} if it is tensored and admits all small
\V-colimits; this is equivalent to its admitting all small weighted
colimits.

If \M\ is cotensored as well, the distinction between ordinary
colimits in $\M_0$ and \V-colimits vanishes; see~\cite[\S3.8]{kelly}.
This applies in particular to the \V-category \V, so that all colimits
in \V\ are \V-colimits.  Thus, if \M\ is both tensored and cotensored,
it is cocomplete if and only if $\M_0$ is cocomplete in the ordinary sense.  This
will usually be the case for our \V-categories.

Dually, of course, we define a limit to be a \emph{\V-limit} when its
defining bijection is induced from a \V-natural isomorphism in \V.
When products and equalizers are \V-limits, the cotensor product
$\cten{W,F}^\D$ of functors $F\maps \D\to\M$ and $W\maps \D\to\V$
coincides with the \emph{weighted limit of $F$ weighted by $G$}.  We
say a \V-category \M\ is \emph{complete} if it is cotensored and
admits all small \V-limits; this is equivalent to admitting all small
weighted limits.  If \M\ is both tensored and cotensored, completeness
of \M\ is equivalent to ordinary completeness of $\M_0$.

We can perform a construction similar to~(\ref{eq:wgt-colim}) on two
functors $F_1\maps \D\to\M$ and $F_2\maps \D\to\M$, both landing in
\M, by replacing the tensor or cotensor with the enriched hom-objects
of \M.  This gives an object of \V\
\begin{equation}\label{eq:nat-trans-as-end}
  \M^\D(F_1,F_2) = \eqlzr\left(
    \prod_{d} \M(F_1 d,F_2d) \rightrightarrows
    \prod_{d,d'} \M\Big(\D(d,d') \odot F_1 d, F_2 d'\Big)
  \right),
\end{equation}
which we think of as the \V-object of \V-natural transformations from
$F_1$ to $F_2$.  These are the hom-objects of a \V-category $\M^\D$
whose objects are \V-functors $\D\to\M$; its underlying category has
\V-natural transformations as its morphisms.  When \M\ is cocomplete,
the $G$-weighted colimit for a fixed $G\maps \D\op\to\V$ extends to a
\V-functor $\M^\D\to \M$, and similarly for weighted limits.

Now suppose that $\M_0$ is a homotopical category.  Then $(\M^\D)_0$
is a homotopical category in which the weak equivalences are
objectwise.  As in the unenriched case, if $\M_0$ is saturated, then
so is $(\M^\D)_0$.  We may then ask, for a given weight $G$, does the
underlying \textbf{Set}-functor $(\M^\D)_0 \to \M_0$ of the
$G$-weighted limit or colimit have a (global) derived functor?  In
other words, is there a \emph{weighted homotopy colimit}?

Pure model category theory has little to say about this question.
Even when $\V_0$ and $\M_0$ are model categories, related in the most
favorable ways known, the category $(\M^\D)_0$ may not have any
relevant model structure.  When $\M_0$ is cofibrantly generated, then
$(\M^\D)_0$ often has a projective model structure, although some
smallness and cofibrancy conditions are necessary; we will consider
this model structure in \autoref{proj-model-str} for a different reason.
However, cofibrant replacements in this model structure are as
unwieldy as ever.  Moreover, injective and Reedy model structures have
no known analogues in the enriched context.

The more technical solution of~\cite{dhks} does not seem to generalize
at all to the enriched case.  There seems no way to define a category
of simplices for a \V-category \D, so we can't even get off the
ground.  However, the results of \S\ref{sec:comp-bar} suggest a
different approach: define an enriched two-sided bar construction and
show that it satisfies the universal property to be a derived functor.
We carry out this idea in the next two sections.

\section{The enriched two-sided bar construction}
\label{sec:enriched-bar}

We now proceed to construct enriched homotopy tensor and cotensor
products of functors, which will give us weighted homotopy limits and
colimits as a special case.  As alluded to previously, our approach is
to give explicit definitions using bar and cobar constructions and
then prove that they are, in fact, derived functors.

We begin by introducing the enriched two-sided bar construction.
Henceforth, we assume all our (large) \V-categories are complete and
cocomplete.

\begin{defn}\label{def:ten-enr-simp-bar}
  Let \M\ be a \V-category, let \D\ be a small \V-category, and let
  $G\maps \D\op\to\V$ and $F\maps \D\to\M$ be \V-functors.  The
  \emph{two-sided simplicial bar construction} is a simplicial object
  of \M\ whose object function is
  \begin{equation*}
    B_n(G,\D,F) = \coprod_{\alpha\maps [n]\to\D}
    \big(G(\alpha_n) \ten \D(\alpha_{n-1},\alpha_n)\ten \dots \ten \D(\alpha_{0},\alpha_1)\big)
    \odot F(\alpha_0)
  \end{equation*}
  and whose faces and degeneracies are defined using composition in
  \D, the evaluation maps $\D(d,d')\odot F(d)\to F(d')$ and
  $\D(d,d')\ten G(d')\to G(d)$, and insertion of identities $E\to
  \D(d,d)$.
\end{defn}

\begin{rem}\label{rmk:not-cartesian}
  When \V\ is not cartesian monoidal, we may not be able to define a
  one-sided bar construction $B(*,\D,F)$.  We would want $*$ to be the
  functor constant at the unit object of \V, but when \D\ is enriched
  such a thing need not exist.  Thus the two-sided version is forced
  on us by the enriched point of view.
\end{rem}

We would now like to geometrically realize the enriched simplicial bar
construction, as we did for the unenriched case in
\S\ref{sec:comp-bar}.  So far, we only know how to take the geometric
realization in a simplicially enriched category, but now we are
already enriched over another category \V.  Inspecting the definition
of geometric realization, however, we see that all we need is a
canonical functor $\Delta^\bullet\maps \DD\to\V$.  We can then define
the \emph{geometric realization} of any simplicial object $\sX$ of \p\
to be
\[|\sX| = \Delta^\bullet\odot_{\DD\op} \sX.\]

Since \V\ is cocomplete, giving such a functor $\Delta$ is equivalent
to giving an adjunction $\sS\rightleftarrows\V$
(see~\cite[3.1.5]{hovey}).  The left adjoint maps a simplicial set $K$
to the tensor product $K_\bullet \times_{\DD} \Delta^\bullet$, while
the right adjoint maps an object $V$ of \V\ to the simplicial set
$\V(\Delta^\bullet, V)$.

We will often need to assume more about this adjunction, however.  If
the left adjoint $\sS\to\V$ is strong monoidal (as defined
in~\cite[XI.2]{maclane}), then the adjunction makes \V\ into a
monoidal \sS-category (an \emph{\sS-algebra} in the sense
of~\cite[4.1.8]{hovey}).  In this case, every \V-enriched category \M\
also becomes simplicially enriched, and if \M\ is tensored or
cotensored over \V, it is also so over \sS.  As explained after
\autoref{unenriched-bar-pullout}, this condition ensures that geometric
realization preserves simplicial homotopy.


Furthermore, if $\V_0$ is a monoidal model category and the adjunction
is Quillen (in addition to being strong monoidal), then any \V-model
category becomes a simplicial model category.  This condition is
useful for us primarily because of \autoref{realiz-quillen}, and because
then simplicial homotopy equivalences in \M\ are necessarily weak
equivalences.  In a more general homotopical context, we will need to
assume the latter property explicitly.

Now we can make the following definitions.
\begin{defn}
  Assume the situation of \autoref{def:ten-enr-simp-bar} and that \V\ has
  a canonical cosimplicial object $\Delta^\bullet\maps \DD\to\V$.  The
  \emph{two-sided bar construction} is the geometric realization of
  the two-sided simplicial bar construction:
  \begin{equation*}
    B(G,\D,F) = |\sB(G,\D,F)|.
  \end{equation*}
\end{defn}

As in the unenriched case, we have the following result.
\begin{lem}\label{bar-pullout}
  \begin{equation*}
    \begin{split}
      B(G,\D,F) &\iso G\ten_\D B(\D,\D,F)\\
      &\iso B(G,\D,\D)\odot_\D F.
    \end{split}
  \end{equation*}
\end{lem}
\begin{proof}
  Just like the proofs of \autoref{unenriched-bar-pullout} and
  \autoref{one-sided-bar-pullout}.
\end{proof}

This result means that $B(G,\D,F)$, which we think of as an
``uncorrected weighted homotopy colimit'', is in fact an ordinary
weighted colimit with a ``fattened-up'' weight $B(G,\D,\D)$.  This
should be familiar from the unenriched case, when we saw that the
uncorrected homotopy colimit could be defined as an ordinary weighted
colimit.

\section{Weighted homotopy colimits, I}
\label{sec:htpy-ten}

We think of the two-sided bar construction introduced in the last
section as an \emph{uncorrected homotopy tensor product of functors}.
In this section, we will construct a corrected version and show that
it defines a derived tensor product.

Recall that in the unenriched case, the bar construction (uncorrected
homotopy colimit) could easily be shown to preserve weak equivalences
between objectwise cofibrant diagrams.  The purpose of the
`correction' was to construct from it a functor preserving \emph{all}
weak equivalences, simply by first replacing all diagrams by
objectwise cofibrant ones.  We want to do the same thing in the
enriched case.

Thus, the first thing we need is a notion of weak equivalence for our
\V-categories.  In \S\ref{sec:derived-enrichment} we will consider a
very well-behaved notion of ``enriched homotopical category'', but for
now, to convey the basic ideas and results, we simply suppose that the
underlying categories $\V_0$ and $\M_0$ are homotopical categories.
We also assume in this section that \V\ comes equipped with a strong
monoidal adjunction $\sS\rightleftarrows\V$, as in
\S\ref{sec:enriched-bar}, so that we have a notion of geometric
realization in \V-categories which preserves simplicial homotopy.

There are two technical issues in the enriched context that did not
arise in the unenriched case.  The first is that we need to assume
some cofibrancy conditions on our shape category \D\ in order for even
the \emph{uncorrected} bar construction to be well-behaved.  These
conditions can be somewhat technical, so for now we sidestep the issue
by defining \D\ to be ``good'' if the bar construction over \D\ turns
out be well-behaved on objectwise cofibrant diagrams.

In \S\ref{sec:cofibrancy} we will discuss what sort of ``cofibrancy''
we must impose on \D\ to make this true.  We will prove, for example,
that for a \V-model categories, any suitably cofibrant \V-category is
good.  However, in many situations weaker conditions suffice.  For
example, in topological situations, instead of cofibrancy in the model
structure, a sort of Hurewicz cofibrancy is usually enough.  The basic
idea is that the simplicial bar construction is ``Reedy cofibrant'' in
a suitable sense, just as in \autoref{unenriched-sb-coft}, but by writing
this out out explicitly we can make it apply to types of
``cofibration'' that aren't necessarily part of a model structure.

We now give the formal definition of goodness.  If $\M_Q$ is a full
subcategory of $\M$, such as a deformation retract, we agree to write
$\M_Q^\D$ for the full sub-\V-category of $\M^\D$ consisting of the
functors which are objectwise in $\M_Q$.  These will be the
counterpart of the ``objectwise cofibrant'' diagrams on which we
expect the bar construction to be well-behaved.

\begin{defn}\label{def:pre-good}
  We say that \D\ is \emph{good} for the tensor $\odot$ of \M\ if there
  exist left deformation retracts $\V_Q$ and $\M_Q$ of \V\ and \M,
  respectively, such that the following conditions hold:
  \begin{itemize}
  \item $B(-,\D,-)$ is homotopical on $\V_Q^{\D\op} \times \M_Q^\D$;
  \item If $F\in \M_Q^\D$ and $G\in \V_Q^{\D\op}$, then
    \begin{itemize}
    \item $B(\D,\D,F)\in \M_Q^\D$ and
    \item $B(G,\D,\D)\in \V_Q^{\D\op}$.
    \end{itemize}
  \end{itemize}
\end{defn}

The second technical issue that arises in the enriched situation is
that in general we may not be able to replace a diagram by one which
is objectwise good.  Even though there is a deformation $Q$ into
$\M_Q$, the functor $Q$ is not in general a \V-functor, so we cannot
apply it objectwise to a \V-functor $F\in \M^\D$ to produce a
\V-functor $QF\in\M_Q^{\D}$.  The problem is that most functorial
deformations are constructed using the small-object argument, which
involves taking a coproduct over a \emph{set} of maps, thereby losing
all information about the enrichment.

We do not have a fully satisfactory solution to this problem.
However, there are many cases in which an objectwise cofibrant
replacement can be constructed.  Sometimes the deformation $Q$ can be
made enriched and applied objectwise, while in other cases we can
replace diagrams by objectwise cofibrant ones without doing it
objectwise.  Situations in which objectwise cofibrant replacements
exist include the following.
\begin{itemize}
\item When \D\ is unenriched.
\item When $\M_Q=\M$ (every object is cofibrant).  Note that the dual
  of this situation, $\M_R=\M$, frequently occurs in topology; this
  makes homotopy limits often easier to deal with than homotopy
  colimits.
\item When \M\ is a cofibrantly generated \V-model category and every
  object of \V\ is cofibrant (for example, when \V=\sS).
\item When $\M^\D$ has a good projective model structure.
\end{itemize}

In \S\ref{sec:objectwise-good} we will consider various partial
solutions of this sort.  For now, we again sidestep the issue by
assuming, when necessary, that such a replacement is possible.
Formally, what we assume is that the subcategories $\V_Q^{\D\op}$ and
$\M_Q^\D$ are left deformation retracts of $\V^{\D\op}$ and $\M^\D$,
respectively.  We use the notations $Q_\V^{\D\op}$ and $Q_\M^\D$ for a
corresponding pair of left deformations.

In many of these cases, the deformations $Q_\V^{\D\op}$ and $Q_\M^\D$
are produced by small object arguments, and therefore can be difficult
to work with.  However, if our given functors $F$ and $G$ are already
in $\M_Q^\D$ and $\V_Q^{\D\op}$, then the corrected homotopy tensor
product to be defined below is weakly equivalent to the uncorrected
version $B(G,\D,F)$.  Thus in this case it is not necessary to apply
the replacement first.

Finally, with all these technical details firmly shoved under the rug,
we can construct global homotopy tensor products.

\begin{defn}\label{def:hoten}
  Let \D\ be a small \V-category which is good for the tensor $\odot$
  of \M, and assume that $\V_Q^{\D\op}$ and $\M_Q^\D$ are left
  deformation retracts of $\V^{\D\op}$ and $\M^\D$, respectively, with
  corresponding left deformations $Q_\V^{\D\op}$ and $Q_\M^\D$.
  Define the \emph{corrected homotopy tensor product} of $F\in \M^\D$
  and $G\in \V^{\D\op}$ to be
  \begin{equation}\label{eq:ho-ten-all}
    G \hoodot_\D F = B(Q_\V^{\D\op} G,\D,Q_\M^\D F).
  \end{equation}
\end{defn}

The cocompleteness of \M\ ensures that the corrected homotopy tensor
product defines a functor
\begin{equation*}
  \hoodot_\D\maps  \big(\V^{\D\op}\big)_0 \times \big( \M^\D \big)_0 \too \M_0.
\end{equation*}
This is the functor that we want to show is the left derived functor
of $\odot_\D$.  Note that since $Q_\V^{\D\op}$ and $Q_\M^\D$ will not
in general be \V-functors, this bifunctor is \emph{not} in general
induced by a \V-bifunctor $\V^{\D\op}\ten\M^\D\to\M$.  We will have
more to say about this later.

The main idea of the proof is the same as for the unenriched case in
\S\ref{sec:comp-bar}: we show that the deformations $Q_\V^{\D\op}$ and
$Q_\M^\D$, followed by $B(\D,\D,-)$, give a left deformation for the
tensor product.  We first introduce a notation for the corresponding
deformation retracts.

\begin{notn}
  Recall that we write $\M_Q^{\D}$ for the full subcategory of $\M^\D$
  determined by the diagrams which are objectwise in $\M_Q$.  We write
  $\M_{BQ}^{\D}$ for the full image of $B(\D,\D,-)$ applied to
  $\M_Q^\D$.  Thus, $\M_{BQ}^{\D}$ is the full subcategory of $\M^\D$
  determined by diagrams isomorphic to a diagram of the form
  $B(\D,\D,F)$ for some $F\in \M_Q^\D$.  Similarly, using
  $B(-,\D,\D)$, we have $\V_{BQ}^{\D\op}$.
\end{notn}

The assumptions of goodness ensure that $\V_{BQ}^{\D\op}\subset
\V_Q^{\D\op}$ and $\M_{BQ}^{\D}\subset \M_Q^{\D}$.  When these are
inclusions of left deformation retracts, and when $\V_Q^{\D\op}$ and
$\M_Q^\D$ are also left deformation retracts of $\V^{\D\op}$ and
$\M^\D$, respectively (as discussed above), we can conclude that
$\V_{BQ}^{\D\op}$ and $\M_{BQ}^{\D}$ are also left deformation
retracts of $\V^{\D\op}$ and $\M^\D$.  The following lemma shows when
this holds.

\begin{lem}\label{bar-repl}
  For any \V-functor $F\maps \D\to\M$, there is a natural map
  \begin{equation*}
    \ep\maps  B(\D,\D,F)\too F.
  \end{equation*}
  which is an objectwise simplicial homotopy equivalence.  Therefore,
  if simplicial homotopy equivalences in $\M_0$ are weak equivalences,
  it is an objectwise weak equivalence, and so $B(\D,\D,-)$ is a left
  deformation.
\end{lem}
\begin{proof}
  Exactly the same as for \autoref{unenriched-bar-repl}.
\end{proof}

Recall, as remarked after \autoref{unenriched-bar-pullout}, that if \M\
is a simplicial model category, then simplicial homotopy equivalences
are necessarily weak equivalences.  This is the case if \V\ is a
simplicial monoidal model category and \M\ is a \V-model category, for
example.

\begin{rem}\label{two-sides}
  Note from the proof of \autoref{unenriched-bar-repl} that the map \ep\
  is induced by tensoring with $F$ from a natural weak equivalence
  $\ep\maps B(\D,\D,\D)\to \D$.  Moreover, $\sB(\D,\D,\D)$ has
  \emph{two} extra degeneracies---one on each side---and thus this
  latter map $\ep$ is a weak equivalence for two different reasons.
  However, each of those reasons is only natural on one side, so only
  one of them can be used to prove that $\ep\maps B(\D,\D,F)\to F$ is
  a weak equivalence.  Dually, however, the other can prove that
  $\ep\maps B(G,\D,\D)\to G$ is a weak equivalence.  This fact will be
  used in the proof of \autoref{thm:hoten} which follows.
\end{rem}

We can now prove the general theorems.  Notice that we must choose, in
defining the deformation, whether to fatten up $F$ or $G$, but it
doesn't matter which.  This is analogous to the fact that in defining
a derived tensor product of modules, we can replace either module by a
resolution.  

\begin{thm}\label{thm:hoten}
  Let \M\ be a \V-category and \D\ a small \V-category, and make the
  following assumptions:
  \begin{itemize}
  \item $\M_0$ and $\V_0$ are homotopical categories;
  \item \D\ is good for the tensor $\odot$ of \M;
  \item Simplicial homotopy equivalences in \V\ and \M\ are weak
    equivalences; and
  \item $\V_Q^{\D\op}$ and $\M_Q^\D$ are left deformation retracts of
    $\V^{\D\op}$ and $\M^\D$, respectively.
  \end{itemize}
  Then
  \begin{equation}\label{eq:tendef-BQ}
    \big(\V_{Q}^{\D\op}\big)_0 \times \big(\M_{BQ}^\D\big)_0
  \end{equation}
  is a left deformation retract for the functor
  \begin{equation}\label{eq:ten}
    \odot_\D\maps \big(\V^{\D\op}\big)_0 \times \big(\M^\D\big)_0 \too \M_0.
  \end{equation}
  Therefore the corrected homotopy tensor product
  \[\hoodot_\D\maps  \big(\V^{\D\op}\big)_0 \times \big(\M^\D\big)_0 \too \M_0\]
  is a left derived functor of~(\ref{eq:ten}).
\end{thm}

\begin{proof}
  By \autoref{bar-repl}, the subcategory~(\ref{eq:tendef-BQ}) is a
  left deformation retract of
  \begin{equation}\label{eq:tendef-Q}
    \big(\V_Q^{\D\op}\big)_0 \times \big(\M_Q^\D\big)_0.
  \end{equation}
  But by assumption, this is in turn a left deformation retract of the
  whole category
  \begin{equation}\label{eq:tendef-all}
    \big(\V^{\D\op}\big)_0 \times \big(\M^\D\big)_0.
  \end{equation}
  Since a left deformation retract of a left deformation retract is a
  left deformation retract, it follows that~(\ref{eq:tendef-BQ}) is a
  left deformation retract of~(\ref{eq:tendef-all}).  Thus, it
  suffices to show that $\odot_\D$ is homotopical on the
  subcategory~(\ref{eq:tendef-BQ}).

  For brevity, write $B_\D F = B(\D,\D,F)$.  Then by \autoref{bar-repl},
  the functor $\Id\times B_\D$ is a left deformation
  of~(\ref{eq:tendef-Q}) into~(\ref{eq:tendef-BQ}), and and since \D\
  is assumed good, the composite
  \[\odot_\D \circ (\Id\times B_\D)\]
  is homotopical on~(\ref{eq:tendef-Q}). Thus, by \autoref{easy-deform},
  it suffices to show that $\odot_\D$ preserves the weak equivalence
  \[(1,\ep B_\D)\maps (G, B_\D B_\D F) \we (G, B_\D F).\]
  for $F\in \M_Q^\D$ and $G\in\V_Q^{\D\op}$.  When we apply $\odot_\D$
  to $(1,\ep B_\D)$, we obtain the left-hand map in the following
  diagram.
  \begin{equation*}
    \xymatrix{
      G \odot_\D B(\D,\D,\D) \odot_\D B(\D,\D,F)
      \ar[d]^{1 \odot \ep \odot 1} \ar[r]^(.54)\iso&
      B(G,\D,\D) \odot_\D B(\D,\D,F)
      \ar[d]^{\ep\odot 1}\\ 
      G \odot_\D \D(-,-) \odot_\D B(\D,\D,F)
      \ar[r]^(.57){\iso} &
      G \odot_\D B(\D,\D, F)
    }
  \end{equation*}
  But because \D\ is good, the map $\ep\maps B(G,\D,\D)\to G$ is a weak
  equivalence between diagrams in $\V_Q^{\D\op}$.  (Here we use
  \autoref{two-sides} and the ``other'' extra degeneracy of \ep.)
  Therefore, once again because \D\ is good, the right-hand map is a
  weak equivalence, and thus so is the left-hand one, as desired.
\end{proof}


\begin{cor}\label{cor:pre-wgtcolim}
  Under the hypotheses of \autoref{thm:hoten}, if $G\in
  \V_Q^{\D\op}$, then the functor
  \begin{equation*}
    G\hoodot_\D- = B(G,\D,Q-)
  \end{equation*}
  is a left derived functor of the weighted colimit $G\odot_\D -$.
\end{cor}
\begin{proof}
  Straightforward from \autoref{thm:hoten}.
\end{proof}

Applying \autoref{cor:pre-wgtcolim} to the case $\V=\sS$ and $G=*$, we have
completed the proof of \autoref{2bar-comp}, since in the unenriched case,
we can apply cofibrant and fibrant approximation functors objectwise,
so that $\M_Q^\D$ and $\M_R^\D$ are always deformation retracts of
$\M^\D$.

The same methods applied to the cobar construction will produce
homotopy weighted limits.  The proofs are all dual, so we merely state
the main definition and results.

\begin{defn}\label{def:pre-co-good}
  Say that \D\ is \emph{good} for the cotensor $\cten{-,-}$ of \M\ if
  there exist a left deformation retract $\V_Q$ of \V\ and a right
  deformation retract $\M_R$ of \M\ such that the following conditions
  hold:
  \begin{itemize}
  \item $C(-,\D,-)$ is homotopical on $\V_Q^{\D} \times \M_R^\D$;
  \item If $F\in \M_R^\D$ and $G\in \V_Q^{\D}$, then
    \begin{itemize}
    \item $C(\D,\D,F)\in \M_R^\D$ and
    \item $B(G,\D,\D)\in \V_Q^{\D}$.
    \end{itemize}
  \end{itemize}
\end{defn}

We write $\M_{CR}^{\D}$ for the full image of $C(\D,\D,-)$ applied to
$\M_R^\D$.

\begin{thm}\label{thm:hocoten}
  Let \M\ be a \V-category and \D\ a small \V-category, and make the
  following assumptions:
  \begin{itemize}
  \item $\M_0$ and $\V_0$ are homotopical categories;
  \item \D\ is good for the cotensor $\cten{-,-}$ of \M;
  \item Simplicial homotopy equivalences in \V\ and \M\ are weak
    equivalences; and
  \item $\V_Q^{\D}$ is a left deformation retract of $\V^{\D}$ and
    $\M_R^\D$ is a right deformation retract of $\M^\D$.
  \end{itemize}
  Then
  \begin{equation}\label{eq:cotendef-CR}
    \big(\V_{Q}^{\D}\big)_0\op \times \big(\M_{CR}^\D\big)_0
  \end{equation}
  is a right deformation retract for the functor
  \begin{equation}\label{eq:coten}
    \cten{-,-}^\D\maps \big(\V^{\D}\big)_0\op \times \big(\M^\D\big)_0 \too \M_0.
  \end{equation}
  Therefore the corrected homotopy cotensor product
  \[\hocten{-,-}^\D\maps  \big(\V^{\D}\big)_0\op \times \big(\M^\D\big)_0 \too \M_0\]
  is a right derived functor of~(\ref{eq:coten}).
\end{thm}

\begin{cor}\label{cor:pre-wgtlim}
  Under the hypotheses of \autoref{thm:hocoten}, if $G\in \V_Q^{\D}$,
  then the functor
  \begin{equation*}
    \hocten{G,-}^\D = C(G,\D,R-)
  \end{equation*}
  is a right derived functor of the weighted limit $\cten{G,-}^\D$.
\end{cor}


\section{Enriched two-variable adjunctions}
\label{sec:enriched-tvas}

The results in the last section convey the basic point of this paper:
that the bar and cobar constructions, suitably corrected, define
homotopy colimits and limits in the ``global'' derived sense.  These
results are sufficient for most applications.  However, they are still
not fully satisfactory, since the original functor $\odot_\D$ was
actually a \V-functor, but our derived functor $\hoodot_\D$ is not
enriched.

The importance of enrichment on homotopy categories has only been
realized relatively recently.  One reason for this is that it was not
until the introduction of closed symmetric monoidal categories of
spectra in the 1990s that stable homotopy theory could be done in a
fully enriched way.
See~\cite{ss:stable-modules,ss:stable-uniq,dug:spectral-enr} for
some examples of the advantages of this perspective.

Thus, we would like to show that the derived functor $\hoodot_\D$ is
itself enriched in a suitable sense.  To achieve this, in the next few
sections we develop a theory of ``enriched homotopical categories'',
which we then leverage to prove that at least the total derived
functor $\thoodot_\D$ can be enriched over the homotopy category
$\Ho\V$.  We will then use this to prove a number of useful results
about homotopy limits and colimits and the homotopy theory of enriched
diagrams.

To get a feel for what enriched homotopical categories should look
like, we consider first the case of model categories.  Even in this
case there is no known truly satisfactory notion of an ``enriched
model category,'' only of an enriched category with a model structure
on its underlying ordinary category that interacts well with the
enrichment.  The most convenient way to ensure this is with Quillen
conditions on the enriched-hom, tensor, and cotensor.  The condition
on the enriched-hom was first written down by Quillen himself for the
special case $\V=\sS$.

The case of more general \V\ is considered in~\cite[ch.~4]{hovey}, in
which are defined the notions of \emph{monoidal model category},
\emph{\V-model category}, \emph{\V-Quillen pair}, and so on.
Moreover, it is proven that the homotopy category of a \V-model
category is enriched over $\Ho\V$, as is the derived adjunction of a
\V-Quillen pair, and so on.

\begin{rem}
  In~\cite{lewis_mandell:mmmc} Lewis and Mandell have described a more
  general notion of an enriched model category which does not require
  that it be tensored and cotensored.  However, since model categories
  are usually assumed to have enough limits and colimits, it seems
  reasonable in the enriched case to assume the existence of weighted
  limits and colimits, including tensors and cotensors.  This is what
  we will do.
\end{rem}

Our notion of ``enriched homotopical category'' is related to this
notion of enriched model category in the same way that deformable
functors are related to Quillen functors.  Namely, instead of lifting
and extension conditions, we assume directly the existence of certain
subcategories on which weak equivalences are preserved.  In
\S\ref{sec:derived-enrichment} we will give this definition and prove
corresponding results about homotopy categories and derived functors,
but first we need to build up some technical tools for constructing
enrichments.  There is no homotopy theory in this section, only
(enriched) category theory.  We omit standard categorical proofs, but
the reader is encouraged to work at least some of them out.

The following notion is the basic tool we will be using.

\begin{defn}[{\cite[ch.~4]{hovey}}]\label{def:tva}
  Let $\M_0$, $\N_0$, and $\p_0$ be unenriched categories.  An
  \emph{adjunction of two variables} $(\oast,\homl,\homr):
  \M_0\times\N_0\to\p_0$ consists of bifunctors
  \begin{align*}
    \oast &\maps  \M_0\times\N_0\to\p_0 \\
    \homl &\maps  \M_0\op\times\p_0\to\N_0\\
    \homr &\maps \N_0\op\times\p_0\to\M_0
  \end{align*}
  together with natural isomorphisms:
  \begin{equation*}
    \p_0(M\oast N, P)\iso \M_0(N,\homl(M, P))) \iso \N_0(M,\homr(N,P)).
  \end{equation*}
\end{defn}

If $\V_0$ is a closed symmetric monoidal category and \M\ is a
tensored and cotensored \V-category, as in \S\ref{sec:enriched}, then
the tensor, cotensor, and hom-objects define a \tva\
$\V_0\times\M_0\to\M_0$.  In particular, this applies to $\M=\V$.  We
encourage the reader to keep these examples in mind.  We have chosen
the (perhaps peculiar-looking) notation $\oast$ for a general \tva\ to
avoid conflict with the monoidal product $\ten$ of \V\ and the tensor
$\odot$ of \M, since sometimes all three will occur in the same
formula.

We can use two-variable adjunctions to give alternate
characterizations of enriched structures.  For instance, a closed
symmetric monoidal category is essentially a symmetric monoidal
category whose monoidal product $\ten$ is part of a \tva.  Similarly,
for enriched categories, we have the following notion
from~\cite[ch.~4]{hovey}.

\begin{defn}\label{def:vmodule}
  Let $\V_0$ be a closed symmetric monoidal category.  A \emph{closed
    \V-module} \M\ consists a category $\M_0$ together with the
  following data.
  \begin{itemize}
  \item A \tva\ $(\odot,\cten{},\hom)\maps \V_0\times\M_0\to\M_0$;
  \item A natural isomorphism $K\odot (L\odot M)\iso (K\ten L)\odot
    M$; and
  \item A natural isomorphism $E\odot M\iso M$;
  \end{itemize}
  such that three obvious coherence diagrams commute.
\end{defn}

Given a closed \V-module, there is an essentially unique way to enrich
$\M_0$ over \V\ (using the given $\hom$ to define the hom-objects)
such that the given \tva\ becomes the tensor-hom-cotensor adjunction.
Thus in the sequel, we frequently blur the distinction between closed
\V-modules and tensored-cotensored \V-categories.  Note that our
closed \V-modules will all be ``left'' modules, in contrast to the
``right'' modules of~\cite{hovey}.

This correspondence carries over to functors as well, via the
following categorical observations.

\begin{prop}\label{lax-vmod-functor}
  Let \M\ and \N\ be closed \V-modules and $F_0\maps \M_0\to\N_0$ an
  ordinary functor.  Then the following data are equivalent.
  \begin{enumerate}
  \item A \V-functor $F\maps \M\to\N$ whose underlying ordinary
    functor is $F_0$.\label{item:5}
  \item A natural transformation $m\maps K\odot F_0X\to F_0(K\odot X)$
    satisfying associativity and unit axioms.\label{item:7}
  \item A natural transformation $n\maps
    F_0(\cten{K,X})\to\cten{K,F_0Y}$ satisfying associativity and unit
    axioms.\label{item:8}
  \end{enumerate}
\end{prop}

We call any of the above data an \emph{enrichment} of $F_0$.  If
in~(\ref{item:7}) $m$ is an isomorphism, we say that $F$ is a
\emph{colax \V-module functor}; this is the only one of these notions
considered in~\cite{hovey}.  Colax \V-module functors correspond to
enrichments of $F_0$ which preserve tensors.  Similarly, if
in~(\ref{item:8}) $n$ is an isomorphism, we call $F$ a \emph{lax
  \V-module functor}; these correspond to enrichments which preserve
cotensors.

The following technical observation will be needed in \autoref{lan-ran-chain}.

\begin{prop}\label{mod-opmod-agree}
  Let  $F_0\maps \M_0\to\N_0$ and let $m$ and $n$ be two
  transformations as in \autoref{lax-vmod-functor}(\ref{item:7})
  and~(\ref{item:8}) such that the following diagram commutes.
  \begin{equation*}
    \xymatrix{
      \M_0(K\odot X, Y) \ar[r]^(0.45){F_0} \ar[d]_\iso &
      \N_0(F_0(K\odot X), F_0Y) \ar[r]^m &
      \N_0(K\odot F_0X, F_0Y) \ar[d]^\iso\\
      \M_0(X,\cten{K,Y}) \ar[r]_(0.45){F_0} &
      \N_0(F_0X,F_0(\cten{K,Y})) \ar[r]_n &
      \N_0(F_0X,\cten{K,F_0Y})
    }
  \end{equation*}
  Then the two induced enrichments on $F_0$ are the same.
\end{prop}

Most of the functors we are interested in are part of adjunctions, and
since in an enriched adjunction the left adjoint automatically
preserves tensors and the right adjoint preserves cotensors, these
results take on an especially simple form for adjunctions.

\begin{prop}\label{vmod-adjn-agree}
  Let \M\ and \N\ be closed \V-modules and let
  $F_0\maps \M_0\rightleftarrows \N_0\spam G_0$ be an ordinary adjunction.
  Then the following data are equivalent.
  \begin{enumerate}
  \item A \V-adjunction $F\maps \M\rightleftarrows\N\spam G$ whose underlying
    ordinary adjunction is the given one;\label{item:6}
  \item A colax \V-module structure on $F_0$; and\label{item:10}
  \item A lax \V-module structure on $G_0$.\label{item:9}
  \end{enumerate}
\end{prop}

This result can be viewed as a special case of the theory of
``doctrinal adjunctions'' presented in~\cite{kelly:doc-adjn}.  In view
of~(\ref{item:6}), we call such data an \emph{enrichment} of the
adjunction $F_0\adj G_0$.  An adjunction $F_0\adj G_0$ together
with~(\ref{item:10}) or~(\ref{item:9}) is also known as a
\emph{\V-module adjunction}.  We summarize the above propositions as
follows.

\begin{prop}\label{enr-as-modules}
  There is an equivalence between the 2-categories of closed
  \V-modules, \V-module adjunctions, and (suitably defined) \V-module
  transformations and of tensored-cotensored \V-categories,
  \V-adjunctions, and \V-natural transformations.
\end{prop}

We can obtain a ``fully enriched'' notion of \tva\ by replacing
everything in \autoref{def:tva} by its \V-analogue.  There is a slight
subtlety in that the ``correct'' definition of an enriched bifunctor
does not use the cartesian product of \V-categories, but rather the
\emph{tensor product}.  The tensor product of two \V-categories \M\
and \N, as defined in~\cite[\S1.4]{kelly}, is a \V-category $\M\ten\N$
whose objects are pairs $(M,N)$ where $M$ and $N$ are objects of \M\
and \N\ respectively, and whose hom-objects are
\begin{equation*}
  (\M\ten\N)\big((M,N),(M',N')\big) = \M(M,M')\ten \N(N,N').
\end{equation*}
Note that unless \V\ is \emph{cartesian} monoidal, we do not in
general have $(\M\ten\N)_0 \eqv \M_0\times\N_0$; all we have is a
functor $\M_0\times\N_0 \to (\M\ten\N)_0$ which is bijective on
objects.

\begin{defn}
  A \emph{\V-adjunction of two variables} $(\oast,\homl,\homr)\maps
  \M\ten\N\to\p$ between \V-categories consists of \V-bifunctors
  \begin{align*}
    \oast &\maps  \M\ten\N\to\p \\
    \homl &\maps  \M\op\ten\p\to\N\\
    \homr &\maps \N\op\ten\p\to\M
  \end{align*}
  together with \V-natural isomorphisms between hom-objects in \V:
  \begin{equation*}
    \p(M\oast N, P)\iso \M(N,\homl(M, P))) \iso \N(M,\homr(N,P)).
  \end{equation*}
\end{defn}
The \tva{}s which arise in \V-category theory, such as the
hom-tensor-cotensor adjunction, are generally \tvva{}s.  Moreover,
just as in the one-variable case, a \tvva\ $\M\ten\N\to\p$ gives rise
to an underlying ordinary \tva\ $\M_0\times\N_0\to\p_0$, and we can
recover the enriched structure from the underlying \tva\ together with
some ``module'' data, as follows.

\begin{defn}
  Let \V\ be a closed symmetric monoidal category and \M, \N, and \p\
  be closed \V-modules.  A \emph{\V-bilinear \tva} consists of a \tva\
  $(\oast,\homl,\homr)\maps \M_0\times\N_0\to\p_0$ together with natural
  isomorphisms
  \begin{align*}
    m_1\maps  K\odot (M\oast N) &\too[\iso] (K\odot M)\oast N\\
    m_2\maps K\odot (M\oast N) &\too[\iso] M\oast(K\odot N)
  \end{align*}
  which make the adjunctions
  \begin{align*}
    (-\oast N) &\Adj \homr(N,-)\\
    (M\oast -) &\Adj \homl(M, -)
  \end{align*}
  into \V-module adjunctions for all $M$ and $N$ and which are such
  that the following diagram commutes:
  \begin{equation}\label{eq:bilin}
    \xymatrix{
      (K\ten L)\odot(M\oast N) \ar[r]^{a} \ar[dd]_s &
      K\odot(L\odot(M\oast N)) \ar[r]^{m_2} &
      K\odot(M\oast(L\odot N)) \ar[d]^{m_1}\\
      && (K\odot M)\oast (L\odot N)\\
      (L\ten K)\odot(M\oast N) \ar[r]_{a} &
      L\odot(K\odot(M\oast N)) \ar[r]_{m_1} &
      L\odot((K\odot M)\oast N) \ar[u]_{m_2}
    }
  \end{equation}
  (Here $a$ denotes use of the associativity isomorphisms of the
  \V-module structures, and $s$ denotes use of the symmetry
  isomorphism of \V).
\end{defn}

\begin{rem}\label{bilin-duals}
  By \autoref{vmod-adjn-agree}, either of the morphisms $m_1$
  and $m_2$ can be replaced by their conjugates $n$ without
  changing the substance of the definition.  For example, we could
  replace $m_1$ by an isomorphism
  \[n_1\maps \homr(N,\cten{K,P}) \too[\iso] \cten{K,\homr(N,P)}.\]
  The Axiom~(\ref{eq:bilin}) would then be replaced by a
  correspondingly dual version.
\end{rem}

\begin{prop}\label{bilin}
  A \tvva\ $(\oast,\homl,\homr)\maps \M\ten\N\to\p$ gives rise to, and can be
  reconstructed from, a \V-bilinear \tva\ $\M_0\times\N_0\to\p_0$.
\end{prop}
\begin{proof}
  This is a sequence of diagram chases left to the reader.  It is
  convenient to use the single-variable description of maps out of a
  tensor product of \V-categories given in~\cite[\S1.4]{kelly}.
\end{proof}

\section{Derived two-variable adjunctions}
\label{sec:derived-tvas}

We now want to define a notion of enriched homotopical category in
such a way that homotopy categories and derived functors can be
enriched.  The way such results are proven for enriched model
categories in~\cite{hovey} is by phrasing enrichment in terms of
\tva{}s, as we have done in \S\ref{sec:enriched-tvas}, and then
showing that every ``Quillen \tva'' has a derived \tva.  Thus, the
central notion for us will be, instead of a Quillen \tva, a
``deformable \tva''.  In this section, we study derived \tva{}s and
derived monoidal structures, remaining in an unenriched context.  Then
in \S\ref{sec:derived-enrichment}, we will apply these results to
construct derived enrichments.

These two sections are somewhat long and technical, and the reader who
is not interested in the details may safely skim them.  The following
two sections \S\ref{sec:enriched-homotopy} and
\S\ref{sec:derived-via-enriched} then use these ideas to provide a
perspective on the general role of enrichment in homotopy theory,
before we return to the homotopy theory of diagrams in
\S\ref{sec:gener-tens-prod}.


Recall from \autoref{def:deform} the notions of \emph{deformation} and
\emph{deformation retract} which are used to produce derived functors.
\autoref{derived-adjn} gives us a general way to produce derived
single-variable adjunctions, but we need derived \tva{}s; thus we
begin by defining a notion of deformation for these.

\begin{defn}\label{def:deform-tva}
  Let \M, \N, and \p\ be homotopical categories.  A \emph{deformation
    retract} for a \tva\ $(\oast,\homl,\homr)\maps \M\times\N\to\p$
  consists of left deformation retracts $\M_Q$ and $\N_Q$ for \M\ and
  \N, respectively, and a right deformation retract $\p_R$ of \p, such
  that
  \begin{itemize}
  \item $\oast$ is homotopical on $\M_Q \times \N_Q$;
  \item $\homl$ is homotopical on $\M_Q\op \times \p_R$; and
  \item $\homr$ is homotopical on $\N_Q\op \times \p_R$.
  \end{itemize}
  A \emph{deformation} for a \tva\ consists of a deformation retract,
  as above, together with corresponding deformations
  \begin{alignat*}{2}
    Q_\M &\maps  \M\too \M_Q &\qquad q_\M&\maps Q_\M\we \Id_\M\\
    Q_\N &\maps  \N\too \N_Q &\qquad q_\N&\maps Q_\N\we \Id_\N\\
    R_\p &\maps  \p\too \p_R &\qquad r_\p&\maps \Id_\p\we R_\p
  \end{alignat*}
  of \M, \N, and \p.
\end{defn}
Of course, for a Quillen \tva, a deformation is given by cofibrant
replacements on \M\ and \N\ and a fibrant replacement on \p.

It is proven in~\cite[40.4]{dhks} than for any single deformable
functor $F$, there is a unique maximal $F$-deformation retract.
Unfortunately, a corresponding result is not true for \tva{}s, as we
will see in \autoref{eg:enriched-mid-nonunique}.  Thus, in the sequel, we
do not consider ``deformability'' to be a \emph{property} of a \tva;
instead, we regard a deformation retract to be \emph{structure} it can
be equipped with.  However, it does follow from~\cite[40.5]{dhks} that
once a deformation retract for a \tva\ is chosen, it does not matter
what deformation we choose; thus we often tacitly choose such a
deformation without mentioning it.

\begin{prop}\label{deform-tva}
   Given a deformation for a \tva\ as above, the functors
   \begin{align*}
     M\hooast N &= Q_\M M \oast Q_\N N\\
     \hohoml(M,P) &= \homl(Q_\M M , R_\p P)\\
     \hohomr(N,P) &= \homr(Q_\N N, R_\p P)
   \end{align*}
   are derived functors of the given ones, and they descend to
   homotopy categories to define a derived \tva\
   $(\thooast,\thohoml,\thohomr)$.
\end{prop}
\begin{proof}
  That they are derived functors follows trivially from
  \autoref{def:deform-tva}.  To construct the derived adjunction
  isomorphisms, we appeal to \autoref{derived-adjn}, which tells us that
  deformable adjunctions descend to derived adjunctions on the
  homotopy categories.  Therefore, for fixed $M$, $N$, and $P$, we have
  derived adjunctions:
  \begin{align*}
    (M\thooast -) &\Adj \thohoml(M,-)\\
    (-\thooast N) &\Adj \thohomr(N,-)\\
    \thohoml(-,P) &\Adj \thohomr(-,P).
  \end{align*}
  Moreover, \autoref{derived-adjn} also tells us that the derived
  adjunction isomorphisms are the \emph{unique} such natural
  isomorphisms which are compatible with the given adjunctions.  It
  follows that the composition of two of the above adjunction
  isomorphisms is the third; hence they give rise to a unique derived
  \tva.
\end{proof}

We can now apply the above result to the defining \tva{}s of a closed
monoidal category to conclude that its homotopy category is also
closed monoidal.  The notation $\V_0$ for our closed monoidal category
is chosen in anticipation of the enriched context of the next section,
but everything is still unenriched for the time being.

\begin{defn}\label{def:mhc}
  Let $\V_0$ be a closed (symmetric) monoidal category which is also a
  homotopical category.  A \emph{closed monoidal deformation retract}
  for $\V_0$ is a deformation retract $(\V_Q,\V_Q,\V_R)$ for the \tva\
  \[(\ten,[-,-],[-,-])\maps  \V_0\times\V_0\to\V_0\]
  such that the following conditions are satisfied:
  \begin{itemize}
  \item $\ten$ maps $\V_Q\times\V_Q$ into $\V_Q$;
  \item $[-,-]$ maps $\V_Q\op\times\V_R$ into $\V_R$;
  \item For all $X\in\V_Q$, the natural maps
    \begin{alignat*}{2}
      QE\ten X  &\too E\ten X\iso X  & \qquad\text{and}\\
      X \ten QE &\too X\ten E\iso X
    \end{alignat*}
    are weak equivalences; and
  \item For all $Y\in\V_R$, the natural map
    \begin{equation*}
      Y\iso [E,Y]  \too{} [QE, Y]
    \end{equation*}
    is a weak equivalence.
  \end{itemize}
  If $\V_0$ is equipped with a closed monoidal deformation retract, we
  say that it is a \emph{closed (symmetric) monoidal homotopical
    category}.
\end{defn}

As for monoidal model categories, the unit conditions are necessary to
ensure that the unit isomorphism descends to the homotopy category.
Any monoidal model category (which, recall, is defined to be closed as
well) is a closed monoidal homotopical category, using cofibrant and
fibrant replacements.


\begin{prop}\label{mhc->ho}
  If $\V_0$\ is a closed (symmetric) monoidal homotopical category,
  then $\Ho(\V_0)$ is a closed (symmetric) monoidal category.
\end{prop}
\begin{proof}
  The proof is exactly like the proof for monoidal model categories
  given in~\cite[4.3.2]{hovey}.  By \autoref{deform-tva}, the monoidal
  product and internal-hom of $\V_0$ descend to the homotopy category
  and define a derived \tva, so it remains to construct the
  associativity, unit, and (if $\V_0$ is symmetric) symmetry
  isomorphisms.
  
  Consider the associativity constraint, which is a natural
  isomorphism
  \begin{equation*}
    a\maps  \ten\circ(\ten\times\Id) \too[\iso] \ten\circ(\Id\times\ten).
  \end{equation*}
  Because $\ten$ maps $\V_Q\times\V_Q$ into $\V_Q$, these two functors
  have the common deformation retract $\V_Q\times\V_Q\times\V_Q$;
  therefore by \autoref{deform-trans}, $a$ gives rise to a derived
  natural isomorphism
  \begin{equation*}
    \mathbf{L}a\maps  
    \mathbf{L}(\ten\circ(\ten\times\Id))
    \too[\iso] \mathbf{L}(\ten\circ(\Id\times\ten)).
  \end{equation*}
  Again because $\ten$ maps $\V_Q\times\V_Q$ into $\V_Q$, both pairs
  $((\Id\times\ten),\ten)$ and $((\ten\times\Id),\ten)$ are
  deformable.  Thus by \autoref{deform-pair} there are canonical natural
  isomorphisms
  \begin{alignat*}{2}
    \mathbf{L}(\ten\circ(\ten\times\Id)) &\too[\iso]
    \thoten\circ(\thoten\times\Id)
    &\qquad\text{and}\\
    \mathbf{L}(\ten\circ(\Id\times\ten)) &\too[\iso]
    \thoten\circ(\Id\times\thoten).
  \end{alignat*}
  Composing these, we get the desired natural associativity
  isomorphism on $\Ho(\V_0)$.  Similarly,
  $\V_Q\times\V_Q\times\V_Q\times\V_Q$ is a common deformation retract
  for all the functors in the pentagon axiom, hence by
  \autoref{deform-trans} its derived version still commutes on the
  homotopy category level.
  
  For the unit isomorphisms, observe that by the 2-out-of-3 property,
  the unit axiom for a closed monoidal deformation retract ensures
  that the pair
  \begin{equation*}
    \xymatrix{
      \V_0 \ar[r]^(0.4){E\times\Id}
      & \V_0\times\V_0 \ar[r]^(0.6)\ten
      & \V_0}
  \end{equation*}
  is deformable.  Thus we can apply similar reasoning to construct the
  unit isomorphisms and verify their axioms.  The symmetry
  isomorphism, in the case when $\V_0$ is symmetric, follows in the
  same way.
\end{proof}

\section{Derived enrichment}
\label{sec:derived-enrichment}

For all of this section, we assume that $\V_0$ is a \csmhc.  Since, by
\autoref{mhc->ho}, $\Ho(\V_0)$ is also a closed symmetric monoidal
category, we can consider categories enriched over $\Ho(\V_0)$.  In
particular, as for any closed symmetric monoidal category, $\Ho(\V_0)$
has an internal hom which produces a $\Ho(\V_0)$-enriched category,
which we denote $\Ho\V$.  The underlying category of $\Ho\V$ is
precisely $\Ho(\V_0)$; in other words, we have
\[(\Ho\V)_0 \iso \Ho(\V_0).\]
Since there is no real ambiguity, from now on we will write $\Ho\V_0$
for this category.  In this section, we will show that the homotopy
theory of \V-categories gives rise to derived $\Ho\V$-categories.

In the next few sections, we will need to be rather pedantic about the
distinction between enriched categories and their underlying ordinary
categories.  It is common in enriched category theory, and usually
harmless, to abuse notation and ignore the distinction between \V\ and
$\V_0$ and between \M\ and $\M_0$, but to avoid confusion amid the
proliferation of notation, here we emphasize that distinction.

\begin{defn}\label{def:vhc}
  Let \M\ be a tensored and cotensored \V-category such that $\M_0$ is
  a homotopical category.  A \emph{deformation retract for the
    enrichment} of \M\ consists of left and right deformation retracts
  $\M_Q$ and $\M_R$ of $\M_0$ such that $(\V_Q,\M_Q,\M_R)$ is a
  deformation retract for the \tva\
  $(\odot,\cten{},\hom)\maps \V_0\times\M_0\to\M_0$, and also the following
  conditions are satisfied.
  \begin{itemize}
  \item $\odot$ maps $\V_Q\times\M_Q$ into $\M_Q$;
  \item $\cten$ maps $\V_Q\op\times\M_R$ into $\M_R$;
  \item For all $X\in\M_Q$, the natural map
    \begin{equation*}
      QE\odot X  \too E\odot X\iso X
    \end{equation*}
    is a weak equivalence; and
  \item For all $Y\in\M_R$, the natural map
    \begin{equation*}
      Y\iso \cten{E,Y}  \too \cten{QE, Y}
    \end{equation*}
    is a weak equivalence.
  \end{itemize}
  If \M\ is equipped with a deformation retract for its enrichment, we
  say it is a \emph{\V-homotopical category}.
\end{defn}

There are several possible choices we could make in this definition,
especially with regards to the dual pairs of axioms relating to the
tensor and cotensor.  We have chosen to include both for simplicity,
but most results depend only on one or the other.

Note that \V\ itself (meaning the \V-enriched category \V) is a
\V-homotopical category.  Its closed monoidal deformation retract
serves as a deformation retract for its self-enrichment.


\begin{prop}\label{vhc->ho}
  If \M\ is a \V-homotopical category, then $\Ho(\M_0)$ is the
  underlying category of a tensored and cotensored $\Ho\V$-category,
  which we denote $\Ho\M$.
\end{prop}
\begin{proof}
  By \autoref{enr-as-modules}, it suffices to construct a $\Ho\V$-module
  structure on $\Ho\M_0$.  We do this in an way entirely analogous to
  \autoref{mhc->ho}.  The tensor-hom-cotensor \tva\ descends to homotopy
  categories because it is assumed deformable, and the other axioms of
  \autoref{def:vhc} ensure that the module isomorphisms and coherence
  descend to derived versions.
\end{proof}

As for $\V$, we have $\Ho(\M_0)\iso (\Ho\M)_0$ so there is no
ambiguity in writing $\Ho\M_0$ for this category.

We now consider \V-enriched adjunctions between \V-homotopical
categories.

\begin{defn}\label{def:vmod-def}
  Let $\Phi\maps \M\to\N$ be a \V-functor between \V-homotopical
  categories.  We say $\Phi$ is \emph{left \V-deformable} if
  \begin{itemize}
  \item $\M_Q$ is a left deformation retract for $\Phi_0$; and
  \item $\Phi$ maps $\M_Q$ into $\N_Q$.
  \end{itemize}
  Similarly, we say $\Phi$ is \emph{right \V-deformable} if
  \begin{itemize}
  \item $\M_R$ is a right deformation retract for $\Phi_0$; and
  \item $\Phi$ maps $\M_R$ into $\N_R$.
  \end{itemize}
\end{defn}

\begin{prop}\label{vmod->ho}
  If $\Phi\maps \M\to\N$ is left \V-deformable, then it has a left derived
  $\Ho\V$-functor
  \[\mathbf{L}\Phi\maps \Ho\M\to\Ho\N.\]
  Similarly, if $\Phi$ is right \V-deformable, then it has a right
  derived $\Ho\V$-functor
  \[\mathbf{R}\Phi\maps \Ho\M\to\Ho\N.\]
\end{prop}
\begin{proof}
  We assume $\Phi$ is left \V-deformable; the other case is dual.
  Since $\Phi_0$ is left deformable, it has a left derived functor.
  We now construct a $\Ho\V$-enrichment on $\mathbf{L}\Phi_0$ using
  \autoref{lax-vmod-functor}.  Let
  \begin{equation*}
    m\maps (\odot\circ(\Id\times \Phi_0))\too (\Phi_0\circ\odot)
  \end{equation*}
  be as in \autoref{lax-vmod-functor}(\ref{item:7}) for $\Phi_0$.
  \autoref{def:vmod-def} guarantees that the composites $\odot \circ
  (\Id\times\Phi_0)$ and $\Phi_0\circ\odot$ have a common deformation
  retract $(\V_Q,\M_Q)$.   Thus, by \autoref{deform-trans}, $m$ descends to
  the homotopy category to give a natural transformation
  \begin{equation}\label{eq:vmod->ho}
    \mathbf{L}(\odot\circ(\Id\times \Phi_0))\too \mathbf{L}(\Phi_0\circ\odot).
  \end{equation}
  \autoref{def:vmod-def} also guarantees that the pairs $((\Id\times
  \Phi_0),\odot)$ and $(\odot,\Phi_0)$ are left deformable, so using
  \autoref{deform-pair}, we can identify the source and target
  of~(\ref{eq:vmod->ho}) with the composite derived functors
  $\thoodot\circ(\Id\times\mathbf{L}\Phi_0)$ and
  $\mathbf{L}\Phi_0\circ\thoodot$.  It is straightforward to verify
  the axioms to show that this transformation is an enrichment of
  $\Phi_0$, as desired.
\end{proof}

\begin{rem}\label{psder-not-enriched}
  Note that the point-set derived functor $\mathbb{L}\Phi_0$ is
  \emph{not} in general a \V-functor, because we have not required
  deformations of \V-homotopical categories to be \V-functors.  In
  fact, cofibrant and fibrant replacements in \V-model categories are
  not in general \V-functors, so this would be too much to assume.
  \autoref{vmod->ho} tells us that the resulting derived functors can
  nevertheless be enriched ``up to homotopy.''
\end{rem}

Since most of the functors we are interested in are part of
\V-adjunctions, from now on we restrict ourselves to such adjunctions
for simplicity.  The following definition is the analogue of a Quillen
adjoint pair for \V-homotopical categories.

\begin{defn}\label{def:vadj}
  Let $\Phi\maps \M\rightleftarrows\N\spam \Psi$ be a \V-adjunction between
  \V-homotopical categories.  Say the adjunction is \emph{left
    \V-deformable} if
  \begin{itemize}
  \item $\M_Q$ is a left deformation retract for $\Phi$;
  \item $\N_R$ is a right deformation retract for $\Psi$; and
  \item $\Phi$ maps $\M_Q$ into $\N_Q$.
  \end{itemize}
  If instead of the last condition, $\Psi$ maps $\N_R$ into $\M_R$, we
  say the adjunction is \emph{right \V-deformable}.
\end{defn}

\begin{prop}\label{vadj->ho}
  If $\Phi\maps \M\rightleftarrows\N\spam \Psi$ is left or right \V-deformable,
  then it has a $\Ho\V$-enriched derived adjunction
  $\mathbf{L}\Phi\maps \Ho\M\rightleftarrows \Ho\N\spam \mathbf{R}\Psi$.
\end{prop}
\begin{proof}
  We assume the adjunction is left \V-deformable; the other case is
  dual.  Since $\Phi\adj \Psi$ is deformable, by
  \autoref{derived-adjn} it has a derived adjunction.  The left
  \V-deformability of the adjunction ensures that $\Phi$ is left
  \V-deformable.  Moreover, since $\Phi$ is a left \V-adjoint, it
  preserves tensors, and hence is actually a colax \V-module functor.
  This implies that $\mathbf{L}\Phi_0$ is also a colax $\Ho\V$-module
  functor, so that by \autoref{vmod-adjn-agree} we have a canonical
  enriched adjunction $\mathbf{L}\Phi\adj \mathbf{R}\Psi$.
\end{proof}

\begin{prop}\label{compose-vadj}
  Let \M\ and \N\ be saturated \V-homotopical categories and let
  $\Phi_1\maps \M\rightleftarrows\N\spam \Psi_1$ and
  $\Phi_2\maps \N\rightleftarrows\p\spam \Psi_2$ be either both left
  \V-deformable or both right \V-deformable.  Then we have
  $\Ho\V$-natural isomorphisms
  \begin{align*}
    \mathbf{L}\Phi_2 \circ \mathbf{L}\Phi_1 &\iso
    \mathbf{L}(\Phi_2\circ \Phi_1)\\
    \mathbf{R}\Psi_1 \circ \mathbf{R}\Psi_2 &\iso
    \mathbf{R}(\Psi_1\circ \Psi_2)
  \end{align*}
\end{prop}
\begin{proof}
  Suppose both are left \V-deformable; the other case is dual.  Then
  the pair $(\Phi_1,\Phi_2)$ is left deformable in the sense of
  \autoref{def:deform-pair}, so the first isomorphism follows from
  \autoref{deform-pair}.  By \autoref{deform-adjn}, since \M\ and \N\ are
  saturated, the pair $(\Psi_2,\Psi_1)$ is right deformable, so the
  second isomorphism follows also from \autoref{deform-pair}.  The fact
  that these isomorphisms are $\Ho\V$-natural can be deduced by
  tracing through the construction of the enrichments and the proof of
  \autoref{deform-pair} in~\cite{dhks}.
\end{proof}

We now consider the analogue of Quillen equivalences.

\begin{prop}\label{vadj-equiv}
  Let $\Phi\maps \M\rightleftarrows\N\spam \Psi$ be a left or right
  \V-deformable \V-adjunction between saturated \V-homotopical
  categories.  Then the following conditions are equivalent.
  \begin{enumerate}
  \item For all $X\in\M_Q$ and $Y\in \N_R$, a map $\Phi X\to Y$ is a
    weak equivalence if and only if its adjunct $X\to \Psi Y$ is.\label{item:1}
  \item For all $X\in\M_Q$ and all $Y\in \N_R$, the composites
    \begin{gather*}
      \xymatrix{X\ar[r]^(0.4)\eta & \Psi\Phi X \ar[r] & \Psi R \Phi X}\\
      \xymatrix{\Phi Q \Psi Y \ar[r] & \Phi\Psi Y\ar[r]^(0.55)\ep & Y}
    \end{gather*}
    are weak equivalences.\label{item:2}
  \item The derived adjunction
    $\mathbf{L}\Phi_0\maps \Ho\M_0\rightleftarrows
    \Ho\N_0\spam \mathbf{R}\Psi_0$ is an adjoint equivalence of
    categories.\label{item:3}
  \item The derived $\Ho\V$-adjunction $\mathbf{L}\Phi\maps
    \Ho\M\rightleftarrows \Ho\N\spam \mathbf{R}\Psi$ is an adjoint
    equivalence of $\Ho\V$-categories.\label{item:4}
  \end{enumerate}
  When these conditions hold, we say that the given adjunction is a
  \emph{\V-homotopical equivalence of \V-homotopical categories}.
\end{prop}
\begin{proof}
  The equivalence of~(\ref{item:1}),~(\ref{item:2}),
  and~(\ref{item:3}) is the same as the corresponding result for
  Quillen equivalences (see, for example,~\cite[1.3.13]{hovey}).  This
  result is stated for general deformable adjunctions between
  homotopical categories in~\cite[\S45]{dhks}.  Note that saturation
  of \M\ and \N\ is only needed for the
  implication~(\ref{item:3})$\Rightarrow$~(\ref{item:1})
  and~(\ref{item:2}).  Finally, it follows from the definition of an
  equivalence of enriched categories (see, for
  example,~\cite[\S1.11]{kelly}) that a $\Ho\V$-adjunction is an
  adjoint equivalence of $\Ho\V$-categories precisely when its
  underlying adjunction is an adjoint equivalence; this proves
  that~(\ref{item:3}) is equivalent to~(\ref{item:4}).
\end{proof}

As in \S\ref{sec:enriched-tvas}, the case of enriched \tva{}s is
slightly more subtle.  In general, the category $(\M\ten\N)_0$ does
not even have an obvious notion of weak equivalence.  However,
\autoref{bilin} tells us that we should really be considering the induced
\V-bilinear \tva\ $\M_0\times\N_0\to\p_0$, and $\M_0\times\N_0$ does
have a notion of weak equivalence induced from $\M_0$ and $\N_0$.
Thus we make the following definition.

\begin{defn}\label{def:tvva-vdr}
  Let $(\oast,\homl,\homr)\maps \M\ten\N\to\p$ be a \tvva\ between
  \V-homotopical categories and suppose that $(\M_Q,\N_Q,\p_R)$ is a
  deformation retract for the underlying \tva\
  $\M_0\times\N_0\to\p_0$.  Then:
  \begin{itemize}
  \item If $\oast$ maps $\M_Q\times\N_Q$ into $\p_Q$, we say the
    adjunction is \emph{$\oast$-\V-deformable};
  \item If $\homr$ maps $\N_Q\op\times\p_R$ into $\M_R$, we say it is
    \emph{$\homr$-\V-deformable}; and
  \item If $\homl$ maps $\M_Q\op\times\p_R$ into $\N_R$, we say it is
    \emph{$\homl$-\V-deformable}.
  \end{itemize}
  If any one of the above conditions obtains we say the adjunction is
  \emph{\V-deformable}.
\end{defn}

Note that the tensor-hom-cotensor \tvva\ of a \V-homotopical category
\M\ is automatically $\odot$-\V-deformable and
$\cten{}$-\V-deformable.

\begin{rem}
  All the notions we have defined so far are generalizations to
  homotopical categories of notions defined for model categories
  in~\cite[ch.~4]{hovey}.  However, enriched \tva{}s are not
  considered there.  We define a \emph{Quillen \tvva} to be a \tvva\
  between \V-model categories whose underlying ordinary \tva\ is
  Quillen.  It is straightforward to show that such a Quillen \tvva\
  is all three kinds of \V-deformable.
\end{rem}

\begin{prop}\label{tvva->ho}
  If a \tvva\ $\M\ten\N\to\p$ is \V-deformable, then it gives rise to
  a derived \tvhva\ $\Ho\M\thoten\Ho\N\to\Ho\p$.
\end{prop}
\begin{proof}
  Entirely analogous to \autoref{vadj->ho}, using \autoref{bilin} in place
  of \autoref{enr-as-modules}.  The $\oast$-\V-deformable case is most
  straightforward; for the other two we must invoke
  \autoref{bilin-duals}.
\end{proof}

Note that we write $\Ho\M\thoten\Ho\N$ for the tensor product of
$\Ho\V$-categories, since the tensor product of $\Ho\V$ is $\thoten$.
\autoref{psder-not-enriched} applies in the two-variable context as well.


\begin{rem}
  While we are dealing with enriched diagrams and enriched functors,
  and producing $\Ho\V$-enriched derived functors and adjunctions, the
  basic notion of ``derived functor'' in use is still an unenriched
  one; that is, its universal property is still \textbf{Set}-based.
  The correct enriched universal property is not entirely clear, but
  in~\cite{lewis_mandell:mmmc}, Lewis and Mandell introduce some
  candidate universal properties for enriched derived functors and
  bifunctors.  It would be interesting to interpret tensor products of
  functors in that framework, especially since monoids in a monoidal
  category (which they are primarily interested in) can be considered
  as one-object enriched categories, with left and right modules
  corresponding to functors of appropriate variance.
\end{rem}


We finish this section by mentioning the following remarkable
consequence of \autoref{vhc->ho}.
\begin{cor}\label{vhc->loc-small}
  If $\Ho\V_0$ has small hom-sets and \M\ is a \V-homotopical
  category, then $\Ho\M_0$ also has small hom-sets.
\end{cor}
\begin{proof}
  For any two objects $M$ and $M'$ of \M, we have
  \begin{equation*}
    \begin{split}
      \Ho\M_0(M,M') &\iso \Ho\M_0\left(E\thoodot M, M'\right)\\
      &\iso \Ho\V_0\big(E,\mathbf{R}\M(M,M')\big).
    \end{split}
  \end{equation*}
  But the last set is a hom-set in $\Ho\V_0$ which was assumed to be
  small; thus $\Ho\M_0(M,M')$ is also small.
\end{proof}

From one point of view, this result is unexpected; we seem to have
gotten something for nothing.  However, if we take the point of view
that homotopy categories are constructed by replacing the objects by
``good'' ones so that weak equivalences become homotopy equivalences,
and then quotienting by homotopy, we can see that in a certain sense
this is precisely what we have done.  Here $Q$ and $R$ give the
``good'' replacement objects, and the notion of homotopy is that
induced from the enrichment over \V.

This is a good example of how the presence of enrichment simplifies,
rather than complicates, the study of homotopy theory.  We remarked in
\S\ref{sec:deformations} that for an unenriched category, we cannot
hope to construct the homotopy category by simply deforming the
hom-functor, but \autoref{vhc->loc-small} shows that in the enriched
situation we can often do precisely that.  In the next two sections,
we will see more examples of this phenomenon.

\section{Enriched homotopy equivalences}
\label{sec:enriched-homotopy}

Let \V\ be a \csmhc\ and let \M\ be a \V-category.  We do \emph{not}
assume yet that \M\ is equipped with any notion of weak equivalence.
Rather, in this section we want to explain how such an \M\ is
automatically equipped with a notion of ``homotopy equivalence''
induced from its enrichment over \V.  In the next section, we will
return to considering categories \M\ that have their own notion of
weak equivalence compatible with that of \V.

Observe that by definition of the derived tensor product in $\Ho\V_0$,
the localization functor $\gamma\maps \V_0\to\Ho\V_0$ is lax symmetric
monoidal.  Therefore, it can be applied to the hom-objects of any
\V-enriched category \M\ to give a $\Ho\V$-enriched category, which we
call denote $h\M$.  We denote its underlying category $(h\M)_0$ by
simply $h\M_0$, since there is no sense to the notation $h(\M_0)$.
Note that there is a canonical functor $\overline{\gamma}\maps \M_0
\to h\M_0$ which is just $\gamma$ applied to the hom-sets $\M_0(M,M')
= \V_0(E,\M(M,M'))$.

It may perhaps seem strange that $h\M$ is enriched over $\Ho\V$,
especially when we take $\M=\V$ and obtain $h\V$ enriched over
$\Ho\V$.  However, the $\Ho\V$-category $h\M$ should be viewed as
merely a convenient way to produce the ordinary category $h\M_0$,
which is really a classical object, as shown by the following
examples.

\begin{exmp}\label{eg:hM}
  If $\V$ is \textbf{Top} with the cartesian product, so that \M\ is
  topologically enriched, then $h\M_0$ is the quotient of $\M_0$ by
  enriched homotopy, since we have
  \begin{align*}
    h\M_0(X,Y)
    &= \Ho\V_0(*,\M(X,Y))\\
    &= \pi_0(\M(X,Y)).
  \end{align*}

  Similarly, if $\V=\sS$, then $h\M_0$ is the quotient of $\M_0$ by
  simplicial homotopy, and if $\V$ is chain complexes with tensor
  product, so that $\M$ is a dg-category, then $h\M_0$ is the quotient
  of $\M_0$ by chain homotopy.
\end{exmp}

These examples lead us to expect more generally that $h\M_0$ should be
the localization of $\M_0$ at a class of weak equivalences induced
from the enrichment.  A natural definition to make is the following.

\begin{defn}\label{def:vequiv}
  A morphism $f\maps X\to Y$ in $\M_0$ is a \emph{\V-equivalence} if
  $\overline{\gamma}(f)$ is an isomorphism in $h\M_0$.
\end{defn}

It is clear that \V-equivalences satisfy the 2-out-of-6 property,
because isomorphisms do, so that $\M_0$ becomes a homotopical category
when equipped with the \V-equivalences.  We do not yet know that it is
\V-homotopical, however.

\begin{exmp}
  In the example $\V=\mathbf{Top}$, we claim that the \V-equivalences
  in any topologically enriched category \M\ are precisely the
  homotopy equivalences.  Both directions follow from the fact that,
  by \autoref{eg:hM}, $h\M_0$ is the quotient of \M\ by enriched
  homotopy.  Suppose first that $f$ is a homotopy equivalence with
  homotopy inverse $g$.  Then $\overline{\gamma}(fg) = 1_Y$ and
  $\overline{\gamma}(gf) = 1_X$ in $h\M_0$, so $\overline{\gamma}(f)$
  is an isomorphism, as desired.

  Now suppose that $f\maps X\to Y$ is a \V-equivalence in $\M_0$.
  Then $\overline{\gamma}(f)$ is an isomorphism in $h\M_0$, so it has
  an inverse $\overline{g}\in h\M_0(Y,X) \iso \pi_0(\M(Y,X))$.  If we
  choose a $g\in\M_0(Y,X)$ such that $\overline{\gamma}(g) =
  \overline{g}$, then $\overline{\gamma}(fg) = 1_Y$ and
  $\overline{\gamma}(gf) = 1_X$.  It follows by \autoref{eg:hM} that $fg$
  and $gf$ are homotopic to $1_Y$ and $1_X$, respectively, in \M, and
  hence that $f$ is a homotopy equivalence, as desired.

  Similarly, one can prove that when \V\ is simplicial sets or chain
  complexes, the \V-equivalences are the simplicial homotopy
  equivalences or chain homotopy equivalences, respectively.  Although
  it seems quite strange, we have managed to access the classical
  notions of homotopy and homotopy equivalence without referring to
  intervals or cylinders!
\end{exmp}

One of the reasons homotopy equivalences are nicer than weak homotopy
equivalences in classical homotopy theory is that any topological
functor will preserve homotopy equivalences, simply by virtue of its
enrichment.  A similar fact is true for all \V-equivalences.

\begin{prop}\label{pres-vequiv}
  If $F\maps \M\to\N$ is any \V-functor between \V-categories, then
  $F_0\maps \M_0\to\N_0$ preserves \V-equivalences.
\end{prop}
\begin{proof}
  The following square of functors commutes by definition.
  \begin{equation*}
    \xymatrix{
      \M_0 \ar[r]^{F_0} \ar[d] & \N_0 \ar[d] \\
      h\M_0 \ar[r]_{hF_0} & h\N_0
      }
  \end{equation*}
  If $f$ in $\M_0$ becomes an isomorphism in $h\M_0$, it stays an
  isomorphism in $h\N_0$.  Therefore, $F_0(f)$ becomes an isomorphism
  in $h\N_0$, so $F_0(f)$ is a \V-equivalence.
\end{proof}

We now consider under what conditions \M\ is a \V-homotopical category
when equipped with the \V-equivalences.  By the definition of
\V-equivalence, the functor
\[h\M(X,Y) = \gamma(\M(X,Y))\]
takes all \V-equivalences to isomorphisms in $\Ho\V_0$.  Thus, if \V\
is saturated, the functor $\M(-,-)$ takes all \V-equivalences to weak
equivalences in $\V_0$.  Moreover, \autoref{pres-vequiv} implies that the
functors $K\odot -$ and $\cten{K,-}$ preserve all \V-equivalences in
\M\ for any $K\in\V$.  Thus, to make \M\ a \V-homotopical category
with the \V-equivalences as weak equivalences and $\M_Q=\M_R=\M_0$, it
would suffice to ensure that the functors $-\odot X$ and $\cten{-,Y}$
preserve weak equivalences in $\V_Q$.  Unfortunately this is not in
general true, except in very special cases.

\begin{prop}\label{vhc-with-vequiv}
  Suppose that all the hom-objects $\M(X,Y)$ are in $\V_R$.  Then \M\
  is a saturated \V-homotopical category with the \V-equivalences as weak
  equivalences.  Moreover, we then have an isomorphism $h\M\iso\Ho\M$
  of $\Ho\V$-categories, and therefore $h\M_0$ is the localization of
  $\M_0$ at the \V-equivalences.
\end{prop}
\begin{proof}
  By the remarks above, to show that \M\ is a \V-homotopical category
  with the \V-equivalences, it remains only to check that $-\odot X$
  and $\cten{-,Y}$ preserve weak equivalences in $\V_Q$.  We prove the
  first; the other is dual.  Let $f\maps K\we L$ be a weak equivalence
  in $\V_Q$.  Let $X,Y\in\M$ and consider the following sequence of
  naturality squares in \V.
  \begin{small}
    \begin{equation*}
      \xymatrix{
        \M(L\odot X,Y) \ar[r]^(0.45)\iso \ar[d]_{f\odot X} &
        \V(L,\M(X,Y)) \ar[r] \ar[d]_{f} &
        \V(QL,R\M(X,Y)) \ar@{=}[r] \ar[d]_{Qf}&
        \Ho\V(L, \M(X,Y)) \ar[d]^{f}\\
        \M(K\odot X,Y) \ar[r]_(0.45)\iso &
        \V(K,\M(X,Y)) \ar[r] &
        \V(QK,R\M(X,Y)) \ar@{=}[r]&
        \Ho\V(K, \M(X,Y))
      }
    \end{equation*}
  \end{small}

  The first pair of natural isomorphisms come from the definition of
  the tensor of \M, and the third pair of equalities come from the
  definition of the closed structure on $\Ho\V_0$.

  The second pair of natural transformations
  \begin{align*}
    \V(L,\M(X,Y)) &\too \V(QL,R\M(X,Y))\\
    \V(K,\M(X,Y)) &\too \V(QK,R\M(X,Y))
  \end{align*}
  is induced by the natural transformations $q\maps Q\we \Id$ and
  $r\maps \Id\we R$ for the closed monoidal deformation of \V.  Since
  $\M(X,Y)\in\V_R$ and $K,L\in \V_Q$, these natural transformations
  are weak equivalences.  Moreover, since the functor
  \[\Ho\V(-,-)\maps \V_0\op\times\V_0\to \V_0\]
  is homotopical, the right vertical map is a weak equivalence as
  well.  By the 2-out-of-3 property, it follows that the left vertical
  map is a weak equivalence.  Applying $\gamma$ to this map, we obtain
  an isomorphism
  \[h\M(K\odot X,Y)\iso h\M(L\odot X,Y)\]
  in $\Ho\V_0$.  Since the left vertical map is \V-natural and
  $\gamma$ is lax symmetric monoidal, this induced map is
  $\Ho\V$-natural in $X$ and $Y$.  Therefore, it induces a natural
  isomorphism
  \[h\M_0(K\odot X,Y)\iso h\M_0(L\odot X,Y).\]
  Since this is true naturally in $Y$, by the Yoneda Lemma the map
  $\overline{\gamma}(f\odot X)$ must be an isomorphism in $h\M_0$, and
  therefore $f\odot X$ is a \V-equivalence in $\M_0$, as desired.

  For the second statement, note that since we are taking
  $Q=R=\Id_{\M_0}$, we have $h\M(X,Y)=\gamma(\M(X,Y))=\Ho\M(X,Y)$.
  Thus $h\M\iso\Ho\M$, and hence $h\M_0\iso\Ho\M_0$.  Since $\Ho\M_0$
  is the localization of $\M_0$ at the \V-equivalences, so is $h\M_0$.
  This implies directly that \M\ is saturated, by definition of the
  \V-equivalences.
\end{proof}

The condition that all hom-objects are in $\V_R$ corresponds to the
frequent requirement in the literature of simplicially enriched
categories (such as in~\cite{hcct}) that they be \emph{locally Kan}.
In the topological literature this requirement is rarely found, since
all topological spaces are fibrant.  In general, we expect to have to
deform the objects of \M\ to make it a \V-homotopical category, and
this may be no easier for the \V-equivalences than for some more
general type of weak equivalence.  Thus, in the next section, we
consider the role of the \V-equivalences in a category that already
has a \V-homotopical structure.

\section{Derived functors via enriched homotopy}
\label{sec:derived-via-enriched}

In this section, let \V\ be a \csmhc\ and let \M\ be a \V-homotopical
category.  In this case, $h\M$ is generally different from $\Ho\M$,
but just as in classical homotopy theory, it is often a useful
stepping-stone to it.  The following result is simply a restatement in
our more general context of a standard result from model category
theory (which, itself, is a generalization of a standard result from
classical homotopy theory).  We agree from now on to write $\M_{QR} =
\M_Q\cap \M_R$.

\begin{prop}\label{we-vequiv}
  Let \M\ be a \V-homotopical category.  Then there is a
  $\Ho\V$-functor $\varphi\maps h\M\to \Ho\M$ such that the following
  diagram commutes:
  \begin{equation*}
    \xymatrix{
      \M_0 \ar[d]_{\overline{\gm}}\ar[dr]\\
      h\M_0 \ar[r]_(0.45){\varphi_0} & \Ho\M_0.
    }
  \end{equation*}
  Moreover,
  \begin{enumerate}
  \item The restriction of $\varphi$ to $h\M_{QR}$ is fully faithful.
    If either $Q$ maps $\M_R$ into itself or $R$ maps $\M_Q$ into
    itself, then this restriction is an
    equivalence.\label{item:we-vequiv-1}
  \item If $X,Y\in \M_{QR}$ and $f\maps X\to Y$ is a weak equivalence,
    then it is a \V-equivalence.  Conversely, if \M\ is saturated,
    then \emph{all} \V-equivalences in \M\ are weak
    equivalences.\label{item:we-vequiv-2}
  \end{enumerate}
\end{prop}
\begin{proof}
  First we construct the functor $\varphi$.  By definition of the
  enrichment in \autoref{vhc->ho}, we have $\Ho\M(X,Y) = \M(QX,RY)$.
  Thus the transformations $QX\to X$ and $Y\to RY$ induce maps
  $\M(X,Y)\to\M(QX,RY)$.  Tracing through the definitions of the
  composition in $h\M$ and $\Ho\M$, it is straightforward to check
  that these maps are functorial, and that the given diagram commutes.

  We now prove~(\ref{item:we-vequiv-1}).  Since $X\in \M_Q$ and
  $Y\in\M_R$, the hom-functor $\M(-,-)$ preserves the weak
  equivalences $QX\we X$ and $Y\we RY$, so $\varphi$ restricted to
  $h\M_{QR}$ is clearly fully faithful.  It remains to show it is
  essentially surjective.  If $Q$ preserves $\M_R$, then for any
  $M\in\M$ we have $QRM\in \M_{QR}$ and a zigzag of weak equivalences
  $M\we RM \leftwe QRM$, so $M$ is isomorphic in $\Ho\M$ to $QRM\in
  h\M_{QR}$.  The other case is similar.

  Finally, to prove~(\ref{item:we-vequiv-2}), let $f\maps X\to Y$ be a
  weak equivalence in $\M_{QR}$.  Then it becomes an isomorphism in
  $\Ho\M_0$, and hence also an isomorphism in $h\M_0$, since
  $\varphi_0$ is fully faithful.  Thus $f$ is a \V-equivalence by
  definition.  Conversely, if $f\maps X\to Y$ is \emph{any}
  \V-equivalence in \M, then it becomes an isomorphism in $h\M_0$,
  hence also in $\Ho\M_0$.  Thus if \M\ is saturated, it must have
  been a weak equivalence.
\end{proof}


Recall that we have seen condition~(\ref{item:we-vequiv-1}) before, in
\autoref{submiddle-left=right}, and that we saw in
\autoref{model-left=right} that this condition is satisfied in all model
categories if $Q$ and $R$ are well chosen.  Similarly, in practice
most homotopical categories are saturated, so this result is actually
quite general.  It implies the following fact, which is familiar from
classical homotopy theory.

\begin{prop}\label{enriched->middlederived}
  Let \M\ and \N\ be \V-homotopical categories such that
  \begin{enumerate}
  \item \N\ is saturated, and\label{enriched-middlederived-1}
  \item either $Q$ preserves $\M_R$ or $R$ preserves $\M_Q$.\label{enriched-middlederived-2}
  \end{enumerate}
  Then \emph{any} \V-functor $F\maps \M\to\N$ has a middle derived
  functor $\Ho\M_0\to\Ho\N_0$.
\end{prop}
\begin{proof}
  Suppose that $Q$ preserves $\M_R$; the other case is dual.  Then the
  natural zigzag $\Id\leftwe Q \we QR$ is a middle deformation and so
  $\M_{QR}$ is a middle deformation retract.  We aim to show it is a
  middle $F$-deformation retract.

  By \autoref{pres-vequiv}, $F$ preserves all \V-equivalences.  But by
  \autoref{we-vequiv}, all weak equivalences in $\M_{QR}$ are
  \V-equivalences, so $F$ maps them to \V-equivalences.  Since \N\ is
  saturated, by \autoref{we-vequiv}, all \V-equivalences in \N\ are weak
  equivalences; thus $F$ is homotopical on $\M_{QR}$ as desired.
\end{proof}

\begin{cor}
  Let \M\ and \N\ be \V-homotopical categories such that
  \begin{itemize}
  \item \N\ is saturated; and
  \item $\M_R=\M$ (``every object is fibrant'').
  \end{itemize}
  Then every \V-functor $F\maps \M\to\N$ has a left derived functor.
  Dually, if the second condition is replaced by $\M_Q=\M$ (``every
  object is cofibrant'') then every \V-functor $F\maps \M\to\N$ has a
  right derived functor.
\end{cor}
\begin{proof}
  Straightforward.
\end{proof}

\begin{rem}
  Even though they were produced using the enrichment, the above
  middle derived functors are \emph{not} in general $\Ho\V$-enriched.
  The functors $Q$ and $R$ are generally not \V-functors, so $FQR$ is
  not a \V-functor.  In adjoint cases like \autoref{vadj->ho}, we can use
  the notion of \V-module functor to get around this problem and
  construct a $\Ho\V$-enrichment, but no such tricks are available in
  the general situation.
\end{rem}

\autoref{enriched->middlederived} shows that the deformations $Q$ and $R$
of a \V-homotopical category are, in a sense, ``universal:'' they work
for \emph{every} (enriched) functor.  And while formally speaking,
middle derived functors are not unique, as we saw in
\autoref{eg:middle-nonunique}, certainly these middle derived functors
produced from the enrichment have a good claim to be the ``correct''
ones.  However, we must beware of too much enthusiasm; we still have
the following enriched version of \autoref{eg:middle-nonunique}.

\begin{ceg}\label{eg:enriched-mid-nonunique}
  Let \V\ be a closed \emph{cartesian} monoidal homotopical category;
  thus we write its product as $\times$ and its unit (which is the
  terminal object) as $*$.  Let \M\ and \N\ be \V-homotopical
  categories satisfying the conditions of
  \autoref{enriched->middlederived}, and let \I\ be the \V-category with
  two objects $0$ and $1$ and with $\I(0,0)=\I(1,1)=\I(0,1)=*$ and
  $\I(1,0)=\emptyset$ (the initial object of \V).  The underlying
  category $\I_0$ is the ordinary category $\I$ considered in
  \autoref{eg:middle-nonunique}.  Let $F^0, F^1\maps \M\to\N$ be two
  \V-functors and let $\alpha\maps F^0\to F^1$ be a \V-natural
  transformation; then $\alpha$ gives rise to a \V-functor $F\maps
  \M\times\I\to\N$ in an obvious way.

  Now, because \V\ is cartesian monoidal, we have $(\M\times\I)_0 \iso
  \M_0\times\I_0$, which we give the product homotopical structure.
  The \V-category $\M\times\I$, however, has two \V-homotopical
  structures, as follows.

  In the ``left'' \V-homotopical structure, we define deformations
  \begin{align*}
    Q^L(M,i) &= (Q_\M M, 0)\\
    R^L(M,i) &= (R_\M M, i)
  \end{align*}
  The \V-category $\M\times\I$ is tensored and cotensored with tensors
  and cotensors defined in the obvious way:
  \begin{align*}
    K\odot (M,i) &= (K\odot M, i)\\
    \cten{K,(M,i)} &= (\cten{K,M}, i).
  \end{align*}
  It is straightforward to check that this satisfies \autoref{def:vhc}
  and condition
  \ref{enriched->middlederived}(\ref{enriched-middlederived-2})
  because \M\ does.  Similarly, in the ``right'' \V-homotopical
  structure, we define the deformations to be
  \begin{align*}
    Q^L(M,i) &= (Q_\M M, i)\\
    R^L(M,i) &= (R_\M M, 1);
  \end{align*}
  this also satisfies \autoref{def:vhc} and condition
  \ref{enriched->middlederived}(\ref{enriched-middlederived-2}).
  Applying \autoref{enriched->middlederived} to these two \V-homotopical
  structures, we produce two different middle derived functors of $F$;
  one is the ``canonical'' middle derived functor of $F^0$ and the
  other is that of $F^1$.
\end{ceg}

We are forced to conclude that even in the enriched situation, the
deformations $Q$ and $R$ must be regarded as part of the structure, in
the sense that different choices can produce genuinely different
behavior.

This completes our study of general enriched homotopy theory in
general.  We now return to the case of most interest to us and
consider diagram categories.


\section{Generalized tensor products and bar constructions}
\label{sec:gener-tens-prod}

In the next section, we will construct homotopy tensor and cotensor
products of functors in a ``fully enriched'' way, using the techniques
of \S\ref{sec:derived-enrichment}.  But now that we have the
perspective of \S\ref{sec:enriched-tvas}, it is clear that there is no
reason to restrict our attention to the tensor-hom-cotensor of a
\V-category; we can define homotopy tensor products of functors in the
more general setting of an arbitrary \tvva.  This also has the
advantage of making the symmetry of the construction more evident.  In
this section we record the obvious definitions of the tensor product
of functors and the bar and cobar constructions in this more general
context.

\begin{defn}\label{def:tensor-product}
  Let $(\oast,\homl,\homr)\maps \M\ten\N\to\p$ be a \tvva, where \p\
  is cocomplete, and let \D\ be a small \V-category.  Given
  \V-functors $G\maps \D\op\to\M$ and $F\maps \D\to\N$, their
  \emph{tensor product} is an object of \p\ defined as the following
  \V-coequalizer:
  \begin{equation}\label{eq:tensor-product}
    G\oast_\D F =
    \coeq\left( \coprod_{d,d'} \D(d,d') \odot \big(Gd'\oast Fd\big)
      \rightrightarrows \coprod_{d} Gd\oast Fd \right).
  \end{equation}
\end{defn}
\begin{prop}
  In the above situation, the tensor product defines a functor
  $\M^{\D\op}\ten\N^\D\to\p$.  Moreover, this functor is part of a
  \tvva.
\end{prop}
The right adjoints to $\oast_\D$ are the functors $\homr$ and $\homl$
applied objectwise; that is, $\homr(F,Z)(d) = \homr(F(d),P)$ and
similarly.

Clearly, in the case when the given \tvva\ is the tensor-hom-cotensor
of a \V-category \M, the tensor product of functors $G\odot_\D F$
reduces to that of \autoref{def:wgt-colim}.  Just as in that case, there
is an obvious dual notion of \emph{cotensor product of functors}.

\begin{defn}\label{def:enr-simp-bar}
  In the situation of \autoref{def:tensor-product}, the \emph{two-sided
    simplicial bar construction} is a simplicial object of \p\ whose
  object function is
  \begin{equation*}
    B_n(G,\D,F) = \coprod_{\alpha\maps [n]\to\D_0}
    \big(\D(\alpha_{n-1},\alpha_n)\ten \dots \ten \D(\alpha_{0},\alpha_1)\big)
    \odot \big(G(\alpha_n)
    \oast F(\alpha_0)\big)
  \end{equation*}
  and whose faces and degeneracies are defined using composition in
  \D, the evaluation maps $\D(d,d')\odot F(d)\to F(d')$ and
  $\D(d,d')\odot G(d')\to G(d)$, and insertion of identities $E\to
  \D(d,d)$.
\end{defn}
Note the three different bifunctors appearing: the monoidal product
$\ten$ in \V, the tensor $\odot$ of $\p$, and the given bifunctor
$\oast$.  Because $\oast$ is part of a \tvva, it ``doesn't matter
which way we parenthesize this expression,'' so we can think of it as:
\begin{equation*}
  \text{``}\quad
  G(\alpha_n) \ten \D(\alpha_{n-1},\alpha_n)\ten\dots
  \ten\D(\alpha_0,\alpha_1) \ten F(\alpha_0).
  \quad\text{''}
\end{equation*}

\begin{defn}\label{def:enr-cosimp-cobar}
  Let $(\oast,\homl,\homr)\maps  \M\ten\N\to\p$ be a \tvva, let \D\ be a
  small \V-category, and let $G\maps \D\to\N$ and $F\maps \D\to\p$ be
  \V-functors.  The \emph{two-sided cosimplicial cobar construction} is a
  cosimplicial object of \M\ whose object function is
  \begin{equation*}
    C^n(G,\D,F) = \prod_{\alpha\maps [n]\to\D_0}
    \homr\Big(
    \big(\D(\alpha_{n-1},\alpha_n)\ten \dots \ten \D(\alpha_{0},\alpha_1)\big)
    \odot G(\alpha_0),
    F(\alpha_n)\Big)
  \end{equation*}
  and whose faces and degeneracies are defined in a way similar to the
  simplicial bar construction.
\end{defn}
We refer the reader to~\cite{meyer_i} for more about the cosimplicial
cobar construction; in almost all respects it is dual to the
simplicial bar construction.

\begin{defn}
  Assume the situation of \autoref{def:enr-simp-bar} and that \V\ has a
  canonical cosimplicial object $\Delta^\bullet\maps \DD\to\V$.  The
  \emph{two-sided bar construction} is the geometric realization of
  the two-sided simplicial bar construction:
  \begin{equation*}
    B(G,\D,F) = |\sB(G,\D,F)|.
  \end{equation*}
  Similarly, in the situation of \autoref{def:enr-cosimp-cobar}, the
  \emph{two-sided cobar construction} is the \emph{totalization} of
  the cosimplicial cobar construction:
  \begin{equation*}
    \begin{split}
      C(G,\D,F) &= \operatorname{Tot}\big(C^\bullet(G,\D,F)\big)\\
      &\equiv \big\{\Delta^\bullet, C^\bullet(G,\D,F)\big\}^{\DD}.
    \end{split}
  \end{equation*}
\end{defn}

\begin{lem}\label{gen-bar-pullout}
  In the above situations, we have
  \begin{equation*}
    \begin{split}
      B(G,\D,F) &\iso G\oast_\D B(\D,\D,F)\\
      &\iso B(G,\D,\D)\oast_\D F.
    \end{split}
  \end{equation*}
  and
  \begin{align*}
    C(G,\D,F) &\iso \homr^\D(G,C(\D,\D,F))\\
    &\iso \homr^\D(B(\D,\D,G), F)
  \end{align*}
\end{lem}

\section{Enriched homotopy tensor and cotensor products}
\label{sec:enrich-htpy-tens}

As promised, we now put the results of \S\ref{sec:htpy-ten} in the
context of \S\ref{sec:derived-enrichment}, so that we can apply the
results of the latter section to produce $\Ho\V$-enrichments of the
total derived tensor products.  Combining the assumptions of those two
sections, we assume that \V\ is a \csmhc\ equipped with a strong
monoidal adjunction $\sS\rightleftarrows\V$.  We work with a general
\tvva, in which the categories \M, \N, and \p\ will be \V-homotopical
categories.  Moreover, clearly the deformation retracts in
\autoref{def:pre-good} should be taken to be those with which
\V-homotopical categories are equipped.  Thus we modify
\autoref{def:pre-good} to the following.

\begin{defn}\label{def:good}
  Let $(\oast,\homl,\homr)\maps \M\ten\N\to\p$ be a \V-deformable
  \tvva\ between \V-homotopical categories, and let \D\ be a small
  \V-category.  We say that \D\ is \emph{good} for $\oast$ if the
  following conditions hold:
  \begin{itemize}
  \item $B(-,\D,-)$ is homotopical on $\M_Q^{\D\op} \times \N_Q^\D$;
  \item If $F\in \N_Q^\D$ and $G\in \M_Q^{\D\op}$, then
    \begin{itemize}
    \item $B(\D,\D,F)\in \N_Q^\D$,
    \item $B(G,\D,\D)\in \M_Q^{\D\op}$, and
    \item $B(G,\D,F)\in \p_Q$.
    \end{itemize}
  \end{itemize}
\end{defn}

\begin{defn}\label{def:hoten-gen}
  Let \D\ be a small \V-category which is good for $\oast$ and let
  $F\in \N^\D$ and $G\in \M^{\D\op}$.  Assume moreover that
  $\M_Q^{\D\op}$ and $\N_Q^\D$ are left deformation retracts of
  $\M^{\D\op}$ and $\N^\D$, respectively, with corresponding left
  deformations $Q_\M^{\D\op}$ and $Q_\N^\D$.  Define the
  \emph{corrected homotopy tensor product} of $G$ and $F$ to be
  \begin{equation}\label{eq:ho-ten-all}
    G \hooast_\D F = B(Q_\M^{\D\op} G,\D,Q_\N^\D F).
  \end{equation}
\end{defn}

The following generalization of \autoref{thm:hoten} is immediate.

\begin{thm}\label{thm:hoten-gen}
  If \D\ is good for $\oast$, simplicial homotopy equivalences in \M\
  and \N\ are weak equivalences, and $\M_Q^{\D\op}$ and $\N_Q^\D$ are
  left deformation retracts of $\M^{\D\op}$ and $\N^\D$, respectively,
  then the corrected homotopy tensor product $\hooast_\D$ is a derived
  functor of $\oast_\D$.
\end{thm}

We also state here the dual versions.  The lemmas and proofs are also
dual, so we omit them and give only the definitions and final results.

\begin{defn}\label{def:co-good}
  Let $(\oast,\homl,\homr)\maps \M\ten\N\to\p$ be a \tvva\ between
  \V-homotopical categories and let \D\ be a small \V-category.  We
  say that \D\ is \emph{good} for $\homr$ if the following conditions
  hold:
  \begin{itemize}
  \item $C(-,\D,-)$ is homotopical on $\N_Q^{\D} \times \p_R^\D$
  \item If $G\in \N_Q^\D$ and $F\in \p_R^{\D}$, then
    \begin{itemize}
    \item $B(\D,\D,G)\in \N_Q^\D$,
    \item $C(\D,\D,F)\in \p_R^{\D}$, and
    \item $C(G,\D,F)\in \M_R$.
    \end{itemize}
  \end{itemize}
\end{defn}

Note that goodness is still a \emph{co}fibrancy condition on \D.  We
will have a little more to say about this in \S\ref{sec:cofibrancy}.

\begin{defn}\label{def:hocoten}
  Let \D\ be a small \V-category which is good for $\homr$, and assume
  that $\N_Q^\D$ and $\p_R^\D$ are, respectively, a left and a right
  deformation retract of $\N^\D$ and $\p^\D$.  Then the
  \emph{(corrected) homotopy cotensor product} of $G\in \N^\D$ and
  $F\in \p^\D$ is
  \begin{equation*}
    \hohomr^\D(G,F) = C(Q_\N^\D G,\D,R_\p^\D F).
  \end{equation*}
\end{defn}

\begin{thm}\label{thm:hocoten}
  If \D\ is good for $\homr$, simplicial homotopy equivalences in \N\
  and \p\ are weak equivalences, and $\N_Q^\D$ and $\p_R^\D$ are,
  respectively, a left and a right deformation retract of $\N^\D$ and
  $\p^\D$, then the corrected homotopy cotensor product $\hohomr^\D$
  is a derived functor of $\homr^\D$.
\end{thm}

Now we want to apply \autoref{tvva->ho} to produce $\Ho\V$-enrichments
for the total derived tensor and cotensor products, but to do this we
first need a \V-homotopical structure on diagram categories.  The
wonderful thing is that we can do both of these things with the very
same construction, since the enriched-hom between diagrams is a
cotensor product construction.  This fulfills the promise made in
\S\ref{sec:coherent-transf} by showing that we can use ``coherent
transformations'', as defined there, to invert the objectwise weak
equivalences and thereby construct the homotopy category of a diagram
category.

Recall that we write $\M_{BQ}^\D$ for the full image of $B(\D,\D,-)$
on $\M_Q^\D$.  Dually, we write $\M_{CR}^\D$ for the full image of
$C(\D,\D,-)$ on $\M_R^\D$.

\begin{thm}\label{dgm-hoenriched}
  Let \M\ be a \V-homotopical category in which simplicial homotopy
  equivalences are weak equivalences, let \D\ be good for the enriched
  hom-functor $\M(-,-)\maps \M\op\ten\M\to\V$, and assume that $\M_Q^\D$
  and $\M_R^\D$ are, respectively, a left and right deformation of
  $\M^\D$.  Then $\M^\D$ is a \V-homotopical category when equipped
  with the deformation retracts
  \begin{equation*}
    \M_{BQ}^\D \qquad\text{and}\qquad \M_R^\D
  \end{equation*}
  and also when equipped with the deformation retracts
  \begin{equation*}
    \M_{Q}^\D \qquad\text{and}\qquad \M_{CR}^\D.
  \end{equation*}
  We call the first the \emph{bar \V-homotopical structure} and the
  second the \emph{cobar \V-homotopical structure} on $\M^\D$.
\end{thm}
\begin{proof}
  By \autoref{thm:hocoten}, the given deformation retracts, together with
  $\V_Q$, define a deformation retract for the tensor-hom-cotensor
  \tva\ of $\M^\D$.  Tensors and cotensors in $\M^\D$ are objectwise,
  and the bar (resp.\ cobar) construction is a colimit (resp.\ limit)
  construction, so it is preserved by tensors (resp.\ cotensors); thus
  the remaining axioms of \autoref{def:vhc} for $\M^\D$ follow from those
  for \M.
\end{proof}

\begin{cor}\label{dgm-locsmall}
  Under the conditions of \autoref{dgm-hoenriched}, the homotopy category
  $\Ho(\M^\D)$ is enriched, tensored, and cotensored over $\Ho\V$, and
  if $\Ho\V_0$ has small hom-sets, then so does $\Ho(\M^\D)_0$.
\end{cor}

Finally, we can complete \autoref{thm:hoten-gen} as follows.

\begin{cor}\label{total-hoten}
  Let $(\oast,\homl,\homr)\maps \M\ten\N\to\p$ be a \V-deformable
  \tvva\ between \V-homotopical categories, let \D\ be good for
  $\oast$, and assume that the bar \V-homotopical structure exists on
  $\M^\D$ and $\N^\D$.  Then the \tvva\
  \[(\oast_\D, \homl, \homr)\maps  \M^{\D\op}\ten \N^\D\too\p\]
  is \V-deformable of the same type as $(\oast,\homl,\homr)$, and thus
  has a ``total derived'' \tvhva\
  \[(\thooast_\D, \thohoml, \thohomr)\maps  \Ho(\M^{\D\op})\thoten\Ho(\N^\D)\too\Ho\p.\]
\end{cor}
\begin{proof}
  We showed in \autoref{thm:hoten} that $Q_\M \times B(\D,\D,Q_\N-)$ is a
  deformation for $\oast_\D$.  Since \D\ is good,
  $B(\D\op,\D\op,Q_\M-)$ has image objectwise in $\M_Q$, from which it
  follows that $\oast_\D$ is homotopical on the image of
  $B(\D\op,\D\op,Q_\M-)\times B(\D,\D,Q_\N-)$.  Since $\homl$ and
  $\homr$ are applied objectwise, it follows that they are homotopical
  on the correct deformation retracts as well.

  If the given adjunction is $\homl$- or $\homr$-\V-deformable, then
  it follows directly that the tensor product one is also.  For the
  $\oast$-\V-deformable case, we invoke the final axiom of goodness to
  see that $\oast_\D$ maps the image of $B(\D\op,\D\op,Q_\M-)\times
  B(\D,\D,Q_\N-)$ into $\p_Q$.
\end{proof}

\begin{rem}\label{rmk:total-hoten-variants}
  It is straightforward to show that if one of the diagram categories
  is given the cobar \V-homotopical structure, the result is still
  true.  If both are given the cobar structure, the result is true for
  $\homl$- and $\homr$-\V-deformability but not necessarily for
  $\oast$-\V-deformability.
\end{rem}

Of course, there is a dual version:

\begin{cor}\label{total-hocoten}
  Let $(\oast,\homl,\homr)\maps  \M\ten\N\to\p$ be a \V-deformable \tvva\
  between \V-homotopical categories, let \D\ be good for $\homr$, and
  assume that simplicial homotopy equivalences in \M\ and \N\ are weak
  equivalences.  Assume that the bar \V-homotopical structure exists
  on $\N^{\D}$ and the cobar \V-homotopical structure exists on
  $\p^\D$.  Then the \tvva\
  \[(\oast, \homl, \homr^\D)\maps  \M\ten \N^\D\too\p^\D\]
  is \V-deformable of the same type as $(\oast,\homl,\homr)$, and thus
  has a ``total derived'' \tvhva\
  \[(\thooast, \thohoml, \thohomr^\D)\maps  \Ho\M\thoten\Ho(\N^\D)\too\Ho(\p^\D).\]
\end{cor}

Combined with \autoref{enriched->middlederived}, these results tell us
that in favorable cases, the bar and cobar constructions are the
``universal'' deformations for functors defined on diagram categories.
This provides our final clue to the surprising ubiquity of the bar
construction.  Unfortunately, to actually apply
\autoref{enriched->middlederived} and produce middle derived functors, we
would need the bar construction to preserve fibrant objects or the
fibrant replacement functor to preserve the bar construction.
However, we will see in the next section that many functors of
interest have genuine left or right derived functors which can also be
produced with the bar construction.  We still don't really understand
the meaning of the ``leftness'' or ``rightness'' of a derived functor,
but it is an undeniably useful and frequently encountered property.

\section{Weighted homotopy colimits, II}
\label{sec:wgt-htpy-lim}

We now apply the general notions of homotopy tensor and cotensor
products to the specific cases of weighted limits and colimits, which
were our original motivation.  Note how natural the proofs are in the
global context of derived functors; this is another value of having
the theorem that local and global notions agree.

In this section, we continue our standing assumption that \V\ is a
\csmhc\ with a strong monoidal adjunction $\sS\rightleftarrows\V$, and
we assume that simplicial homotopy equivalences are weak equivalences
in all \V-homotopical categories, including \V, and that all
homotopical categories are saturated.  All these conditions are
satisfied if \V\ is a simplicial monoidal model category and we
consider \V-model categories only.

We also assume that for all \V-homotopical categories \M\ and all
small \V-categories \D\ considered, the subcategories $\M_Q^\D$ and
$\M_R^\D$ are, respectively, a left and right deformation retract of
$\M^\D$, so that the bar and cobar \V-homotopical structures exist.
Unless otherwise specified, we give diagram categories such as $\M^\D$
the bar \V-homotopical structure.

\begin{thm}\label{cor:wgtcolim}
  If \D\ is good for the tensor $\odot$ of \M\ and $G\maps \D\op\to\V$ is
  objectwise in $\V_Q$, then the \V-adjunction
  \[G\odot_\D- \maps  \M^\D \rightleftarrows \M \spam \cten{G,-}\]
  is both left and right \V-deformable.  Therefore the functor
  \begin{equation*}
    G\hoodot_\D- = B(G,\D,Q-)
  \end{equation*}
  is a left derived functor of the weighted colimit $G\odot_\D -$, and
  we have a derived $\Ho\V$-adjunction
  \[G\thoodot_\D- \maps \Ho(\M^\D)\rightleftarrows \Ho\M \spam \thocten{G,-}\]
\end{thm}
\begin{proof}
  Straightforward from \autoref{total-hoten}.
\end{proof}

We also have a dual version using \autoref{thm:hocoten}. which constructs
homotopy weighted limits, using the cobar \V-homotopical structure, as
follows.

\begin{thm}\label{cor:wgtlim}
  If \D\ is good for the cotensor $\cten{}$ of \M\ and $G\maps \D\to\V$ is
  objectwise in $\V_Q$, then the \V-adjunction
  \[G\odot -\maps  \M \rightleftarrows \M^\D \spam \cten{G,-}^\D\]
  is both left and right \V-deformable.  Therefore the functor
  \begin{equation*}
    \hocten{G,-}^\D = C(G,\D,R-)
  \end{equation*}
  is a right derived functor of the weighted limit $\cten{G,-}^\D$, and
  we have a derived $\Ho\V$-adjunction
  \[G\thoodot- \maps \Ho\M\rightleftarrows \Ho(\M^\D) \spam \thocten{G,-}^\D\]
\end{thm}

This dual result is especially important because, as remarked earlier,
while projective model structures often exist on the \V-category
$\M^\D$ of enriched diagrams, so that homotopy colimits of enriched
diagrams can often be defined (if not computed) in that way, there is
hardly ever an injective model structure on enriched diagram
categories.  Thus cobar constructions provide the \emph{only} viable
way to compute homotopy weighted limits.

An especially interesting sort of weighted colimit is a coend.  Recall
(from~\cite[\S3.10]{kelly}, for example) that the coend of a functor
$H\maps \D\op\ten\D\to \M$ is defined to be the weighted colimit
\[\int^\D H = \D(-,-)\odot_{\D\op\ten\D} H.\]
Thus, \autoref{cor:wgtcolim} tells us that the \emph{homotopy coend},
defined as
\begin{equation}
  \begin{split}
    \mathbb{L}\!\!\int^{d\in\D} H(d,d) &= \D\hoodot_{\D\op\ten\D} H\\
    &= B(\D,\D\op\ten\D,QH)
  \end{split}\label{eq:hocoend}
\end{equation}
is a derived functor of the usual coend functor, at least if \D\ is
good and also all its hom-objects $\D(d,d')$ are in $\V_Q$.  For
example, in~\cite[\S{}IV.1]{ekmm}, the topological Hochschild homology
(THH) of an $E_\infty$ ring spectrum $R$ with coefficients in a
bimodule $M$ is defined as the derived smash product
$R\hosmash_{R\wedge R\op} M$, which is a homotopy version of the coend
$R\wedge_{R\wedge R\op} M\iso\int^R M$.

However, there is a more economical construction of homotopy coends
which is also sometimes useful.  To motivate this, notice that most
coends which arise in practice are tensor products of functors.  The
tensor product which we have been writing $G\oast_\D F$ can be shown
to be isomorphic to the coend
\[\int^\D G \xoast F\]
where the ``external tensor product'' is defined by $(G \xoast
F)(d,d') = Gd\oast Fd'$.  Applying~(\ref{eq:hocoend}) to this naively,
we would obtain the homotopy coend
\[\mathbb{L}\!\!\int^\D G \xoast F = B(\D,\D\ten\D\op,Q(G\xoast F)).
\]
But we know that the homotopy tensor product of functors may be given
instead by the simpler expression $B(QG,\D,QF)$.  In fact, there is an
analogous simplification that works for all coends, using a
generalized notion of bar construction described in~\cite{meyer_i}
and~\cite{meyer_ii}, for example.  This bar construction associates to
a bifunctor $H\maps \D\op\ten\D\to\M$ the realization $B(\D,H)$ of the
simplicial object
\begin{equation*}
  B_n(\D,H) =  \coprod_{\alpha\maps [n]\to\D}
  \big(\D(\alpha_{n-1},\alpha_n)\ten \dots \ten \D(\alpha_{0},\alpha_1)\big)
  \odot H(\alpha_n,\alpha_0).
\end{equation*}
One can then modify the proof of \autoref{thm:hoten} to show that under
suitable conditions on \D, this construction also defines a homotopy
coend.  The ``cyclic bar constructions'' frequently used to compute
Hochschild homology are special cases of this latter type of bar
construction.

We end this section by considering how homotopy limits and colimits
behave under the action of functors which change the target category
\M.  First we need to know that \V-deformable functors can be applied
objectwise to diagrams and remain deformable.  The restriction to
adjunctions in the following proposition, rather than more general
\V-functors, is merely for simplicity.

\begin{prop}\label{induced-hoadj}
  Let $\Phi\maps  \M\rightleftarrows\N\spam \Psi$ be a left (resp.\ right)
  \V-deformable \V-adjunction.  Then the induced \V-adjunction
  \[\Phi^\D\maps  \M^\D \rightleftarrows \N^\D \spam \Psi^\D\]
  is also left (resp.\ right) \V-deformable and hence gives rise to a
  derived $\Ho\V$-adjunction
  \[\mathbf{L}\Phi^\D\maps  \Ho(\M^\D)\rightleftarrows\Ho(\N^\D)\spam \mathbf{R}\Psi^\D.\]
  If the original adjunction $\Phi\adj\Psi$ was a \vheovhc, then so is
  the induced one $\Phi^\D\adj\Psi^\D$.
\end{prop}
\begin{proof}
  Since $\Phi$ is homotopical on $\M_Q$, $\Phi^\D$ is homotopical on
  $\M_Q^\D$ and therefore on the image of $B(\D,\D,Q-)$.  Similarly,
  $\Psi$ is homotopical on $\N_R$, so $\Psi^\D$ is homotopical on
  $\N_R^\D$.  Suppose that the given adjunction is left \V-deformable;
  then $\Phi$ maps $\M_Q$ to $\N_Q$, and since it is a left adjoint it
  preserves the bar construction, so it maps the image of
  $B(\D,\D,Q-)$ in $\M^\D$ to the corresponding image in $\N^\D$.
  Thus the induced adjunction is also left \V-deformable.  The case
  when it is right deformable is even easier.

  Finally, since $\Phi^\D$ and $\Psi^\D$ are simply $\Phi$ and $\Psi$
  applied objectwise, it is clear from condition~(\ref{item:1}) of
  \autoref{vadj-equiv} that if $\Phi\adj\Psi$ is a \vheovhc, so is
  $\Phi^\D\adj\Psi^\D$.
\end{proof}

\begin{prop}\label{compose-induced-hoadj}
  Let $\Phi_1\adj \Psi_1$ and $\Phi_2\adj \Psi_2$ be composable
  \V-adjunctions which are either both left \V-deformable or both
  right \V-deformable.  Then we have $\Ho\V$-natural isomorphisms
  \begin{align*}
    \mathbf{L}\Phi_2^\D \circ \mathbf{L}\Phi_1^\D &\iso
    \mathbf{L}(\Phi_2\circ \Phi_1)^\D\\
    \mathbf{R}\Psi_1^\D \circ \mathbf{R}\Psi_2^\D &\iso
    \mathbf{R}(\Psi_1\circ \Psi_2)^\D
  \end{align*}
\end{prop}
\begin{proof}
  This follows directly from \autoref{compose-vadj}.
\end{proof}

With this result in hand, it it is easy to show that left derived
functors of left adjoints preserve homotopy colimits, just as left
adjoints preserve ordinary colimits, and dually.

\begin{prop}\label{preserve-hocolim}
  Let $\Phi\maps \M\rightleftarrows \N\spam \Psi$ be a left or right
  \V-deformable \V-adjunction, let \D\ be good for the tensors of \M\
  and \N, and let $G\maps \D\op\to\V$ be objectwise in $\V_Q$.  Then
  there is a $\Ho\V$-natural isomorphism
  \[
  G\thoodot_\D \mathbf{L}\Phi^\D(F) \iso \mathbf{L}\Phi(G\thoodot_\D F).
  \]
  Similarly, if \E\ is good for the cotensors of \M\ and \N\ and
  $G\maps \E\to\V$ is objectwise in $\V_Q$, then there is a
  $\Ho\V$-natural isomorphism
  \[
  \thocten{G, \mathbf{R}\Psi^\E(F)}^\E \iso \mathbf{R}\Psi\left(\thocten{G,F}^\E\right).
  \]
\end{prop}
\begin{proof}
  This follows from Theorems~\ref{cor:wgtcolim} and~\ref{cor:wgtlim}
  and Propositions~\ref{induced-hoadj} and~\ref{compose-vadj}.
\end{proof}

\section{Homotopy theory of enriched diagrams}
\label{sec:homot-theory-enrich}

In this section, we consider the effect of functors that change the
shape category \D.  These give rise to Kan extension functors on
diagram categories.  Namely, given $K\maps \D\to\E$, we can restrict
along $K$ giving $K^*\maps \M^\E\to\M^\D$, and this functor has left
and right adjoints $\lan_K$ and $\ran_K$.  The homotopical behavior of
these functors is important in many different contexts, such as the
comparison in~\cite{mmss} of various types of diagram spectra,
including symmetric spectra and orthogonal spectra.  The
``prolongation'' functors in that paper are left Kan extensions for
which the enrichment, as considered here, is an essential aspect.

We continue our standing assumptions from the last section.  Recall
from~\cite[\S4.1]{kelly} that (enriched) left Kan extensions can be
computed as tensor products of functors, $\lan_K F(d) \iso
\D(K-,d)\odot_\D F$.  Thus, in addition to \autoref{def:good} and
\autoref{def:co-good}, for this section we make the following definition.

\begin{defn}\label{def:very-good}
  A small \V-category \D\ is \emph{very good for $\odot$} (resp.\
  \emph{very good for $\cten{}$}) if it is good for $\odot$ (resp.\
  good for $\cten{}$) and moreover all its hom-objects are in $\V_Q$.
\end{defn}

This condition may initially seem unreasonably strong, especially in
topological contexts, but the key is that in such contexts $\V_Q$ does
not have to be the actual cofibrant objects in a model structure;
often something much weaker suffices.  We will say more about these
situations in \S\ref{sec:cofibrancy}.

\begin{prop}\label{derived-lan}
  If $K\maps \D\to\E$ is a \V-functor and \D\ is very good for $\odot$,
  then the adjunction
  \[\lan_K\maps \M^\D \rightleftarrows \M^\E\spam K^*\]
  is right \V-deformable.  It therefore has a derived
  $\Ho\V$-adjunction
  \[\tholan_K\maps  \Ho(\M^\D) \rightleftarrows \Ho(\M^\E) \spam\mathbf{R}K^*.\]
  Moreover, for any other \V-functor $H\maps \C\to\D$, where \C\ is
  also very good for $\odot$, we have a $\Ho\V$-natural isomorphism
  \[\tholan_K\circ \tholan_H \iso \tholan_{KH}.\]
\end{prop}
\begin{proof}
  It follows as in \autoref{cor:wgtcolim} that $\lan_K$ is homotopical on
  the image of $B(\D,\D,Q-)$, and $K^*$ is homotopical everywhere
  since weak equivalences are objectwise.  Moreover, for the same
  reason $K^*$ maps $\M_R^\E$ into $\M_R^\D$, so the adjunction is
  right \V-deformable.  The last statement follows from
  \autoref{compose-vadj}.
\end{proof}

\begin{prop}\label{derived-ran}
  If $K\maps \D\to\E$ is a \V-functor and \D\ is very good for $\cten{}$,
  then the adjunction
  \[K^*\maps \M^\E \rightleftarrows \M^\D\spam \ran_K\]
  is left \V-deformable when the diagram categories are given the
  cobar \V-homotopical structure.  It therefore has a derived
  $\Ho\V$-adjunction
  \[\mathbf{L}K^*\maps  \Ho(\M^\E) \rightleftarrows \Ho(\M^\D) \spam\thoran_K.\]
  Moreover, for any other \V-functor $H\maps \C\to\D$, where \C\ is
  also very good for $\cten{}$, we have
  \[\thoran_K\circ \thoran_H \iso \thoran_{KH}.\]
\end{prop}

\begin{rem}\label{lan-ran-chain}
  Here we see a relatively common phenomenon.  We have a functor $K^*$
  which has both left and right adjoints, and both left and right
  derived functors (when \D\ is very good for both $\odot$ and
  $\cten{}$).  A natural question is whether the two agree.  In this
  case, it is easy to see that they do.  Since $K^*$ is already
  homotopical, it descends to homotopy categories directly, and the
  result is equivalent to any other derived functor of it; thus
  $\mathbf{L}K^* \iso \mathbf{R}K^*$.  In the terminology of
  \autoref{submiddle-unique}, we could say that the identity functor is a
  middle deformation for $K^*$.

  So far, however, this only gives agreement as unenriched functors,
  and we would like to know that the $\Ho\V$-enrichments also
  coincide.  It is straightforward, if a bit tedious, to check that
  the hypothesis of \autoref{mod-opmod-agree} is satisfied for $K^*$, so
  that the two enriched functors are also the same; thus we have a
  derived string of $\Ho\V$-adjunctions:
  \[\tholan_K \adj \mathbf{R}K^* \iso \mathbf{L}K^* \adj \thoran_K.\]
\end{rem}

\begin{prop}\label{lan-equiv}
  Let $K\maps \D\to\E$ be a \V-functor such that
  \begin{itemize}
  \item \D\ and \E\ are both very good for the tensor $\odot$ of \M;
  \item Each map $K\maps \D(d,d')\to \E(Kd,Kd')$ is a weak equivalence in
    \V; and
  \item $K_0$ is essentially surjective.
  \end{itemize}
  Then the \V-adjunction
  \[\lan_K\maps \M^\D \rightleftarrows \M^\E\spam K^*\]
  is a \vheovhc.
\end{prop}
\begin{proof}
  We use condition~(\ref{item:2}) of \autoref{vadj-equiv}.  Let $F\in
  \M_{Q}^\D$, so that $B(\D,\D,F)\in\M_{BQ}^\D$; we must show first
  that for each $d\in D$,
  \[B(\D(-,d),\D,F) \too B(\E(K-,Kd),\D,F) \too R B(\E(K-,Kd),\D,F)\]
  is a weak equivalence.  But since $F\in\M_Q^\D$, \D\ and \E\ are
  both very good, and each $\D(d,-)\to \E(Kd,K-)$ is a weak
  equivalence by assumption, the first map is a weak equivalence, and
  the second is trivially a weak equivalence.

  Now let $G\in \N_R^\D$; we must show that for each $e\in\E$,
  \begin{equation}
    B(\E(K-,e),\D,GK) \too \E(K-,e)\odot_\D GK \too Ge\label{eq:lan-equiv-right}
  \end{equation}
  is a weak equivalence.  Since $K$ is essentially surjective, we may
  assume that $e=Kd$ for some $d\in\D$.  Consider the composite
  \[B(\D(-,d),\D,GK)\too B(\E(K-,Kd),\D,GK)\too GKd.\] Since $K\maps
  \D(-,d)\to\E(K-,Kd)$ is a weak equivalence and \D\ and \E\ are very
  good, the left map is a weak equivalence.  The composite is just the
  map \ep\ from \autoref{bar-repl}, so it is also a weak equivalence.
  Thus, by the 2-out-of-3 property, the desired map is a weak
  equivalence.
\end{proof}

In the case when \D\ and \E\ have one object, this result implies that
a weak equivalence between ``cofibrant'' monoids in a monoidal
homotopical category induces an equivalence between their categories
of modules.

We would like to weaken the condition of essential surjectivity of $K$
in \autoref{lan-equiv} to a ``homotopical'' one.  A natural set of
definitions is the following.

\begin{defn}
  Let $K\maps \D\to\E$ be a \V-functor.
  \begin{itemize}
  \item Say $K$ is \emph{homotopically fully faithful} if each map
    $K\maps \D(d,d')\to \E(Kd,Kd')$ is a weak equivalence in \V.
  \item Say $K$ is \emph{homotopically essentially surjective} if
    every object $e$ is connected to some object $Kd$ in the image of
    $K$ by a zigzag of \V-equivalences (as in \autoref{def:vequiv}) in
    $\E_0$.
  \item If $K$ is both homotopically fully faithful and homotopically
    essentially surjective, say it is a \emph{homotopical equivalence
      of \V-categories}.
  \end{itemize}
\end{defn}

\begin{rem}
  If $K$ is homotopically fully faithful, then the induced
  $\Ho\V$-functor $hK\maps h\D\to h\E$ is fully faithful.  The converse is
  true if \V\ is saturated.  Similarly, if $K$ is homotopically
  essentially surjective, then $hK$ is essentially surjective.  The
  converse is true if all hom-objects $\E(e,e')$ are in $\V_R$, since
  then by \autoref{vhc-with-vequiv} $h\E_0$ is the localization of $\E_0$
  at the \V-equivalences, so any isomorphism in $h\E_0$ must be the
  image of a zigzag of \V-equivalences in $\E_0$.

  In particular, this shows that in the case $\V=\sS$, a homotopical
  equivalence of \sS-categories is in particular a
  \emph{DK-equivalence}, and that the converse is true for
  simplicially enriched categories which are locally Kan.  In our
  terminology, a DK-equivalence is a \sS-functor $K$ which is
  homotopically fully faithful and such that $(hK)_0$ is an
  equivalence.  The notion of DK-equivalence was first defined
  in~\cite{dk-function-complexes}, where it was called a ``weak
  equivalence of simplicial categories;'' see also~\cite{bergner}.
\end{rem}

We can now prove the following generalization of \autoref{lan-equiv}.

\begin{prop}\label{lan-equiv-gen}
  Let $K\maps \D\to\E$ be a homotopical equivalence of \V-categories such
  that \D\ and \E\ are very good for the tensor $\odot$ of \M, and
  assume that \V\ and \M\ are saturated.  Then the \V-adjunction
  \[\lan_K\maps \M^\D \rightleftarrows \M^\E\spam K^*\]
  is a \vheovhc.
\end{prop}
\begin{proof}
  The first half of the proof of \autoref{lan-equiv} still works, and the
  second half shows that~(\ref{eq:lan-equiv-right}) is a weak
  equivalence for all $e=Kd$.  Thus it remains only to prove that
  if~(\ref{eq:lan-equiv-right}) is a weak equivalence for $e$ and
  $e\uto e'$ is a \V-equivalence (in either direction),
  then~(\ref{eq:lan-equiv-right}) is a weak equivalence for $e'$.  In
  this case we have the following commuting square:
  \begin{equation*}
    \xymatrix{
      B(\E(K-,e),\D,GK) \ar[r]^(0.7)\sim \ar@{-}[d] & Ge\ar@{-}[d]\\
      B(\E(K-,e'),\D,GK) \ar[r]  & Ge'
    }
  \end{equation*}
  Since $e\uto e'$ is a \V-equivalence in \E, the induced map
  $\E(K-,e)\uto\E(K-,e')$ is a \V-equivalence in $\V$, hence a weak
  equivalence since \V\ is saturated.  Since \D\ and \E\ are very
  good, this makes the left vertical map a weak equivalence.
  Similarly, since $G$ is a \V-functor, by \autoref{pres-vequiv} it
  preserves all \V-equivalences, so since \M\ is saturated, the right
  vertical map is also a weak equivalence.  By the 2-out-of-3
  property, the bottom horizontal map is also a weak equivalence, as
  desired.
\end{proof}

There is, of course, a dual result for right Kan extensions.  By
\autoref{lan-ran-chain} and the uniqueness of adjoint equivalences, it
follows that if $K$ satisfies the hypotheses of \autoref{lan-equiv} or
\autoref{lan-equiv-gen}, then $\tholan_K \iso \thoran_K$.

We now consider the associativity of the derived tensor product of
functors.


\begin{prop}
  Let \D\ and \E\ be good for the product $\ten$ of \V\ and consider
  \V-functors $H\maps \D\op\to\V$, $G\maps \E\op\ten\D\to\V$, and
  $F\maps \E\to\V$.  Then we have a $\Ho\V$-natural isomorphism
  \begin{equation*}
    H \thoten_\D (G\thoten_\E F)
    \iso
    (H \thoten_\D G)\thoten_\E F.
  \end{equation*}
\end{prop}
\begin{proof}
  We use a two-variable version of \autoref{compose-vadj}.  Both
  $\ten_\D$ and $\ten_\E$ are \V-deformable, so it suffices to check
  that $\ten_\D$ maps $(\V^{\D\op})_Q\ten (\V^{\E\op\ten\D})_Q$ into
  $(\V^\E)_Q$, and similarly for $\ten_\E$.  By $(-)_Q$ we refer here
  to the bar \V-homotopical structure.  It is straightforward to check
  that
  \[B(\D\ten\E\op,\D\ten\E\op,QG) \iso B(\E\op,\E\op,B(\D,\D,QG)).\]
  Since $\ten_\D$ is applied pointwise with respect to \E, it commutes
  with $B(\E\op,\E\op,-)$.  Thus $\ten_\D$ maps the image of the
  deformation
  \[B(\D\op,\D\op,Q-)\times B(\D\ten\E\op,\D\ten\E\op,Q-)\]
  into the image of $B(\E\op,\E\op,Q-)$, as desired.  The case of
  $\ten_\E$ is similar.
\end{proof}

We note in passing that combining this result with methods similar to
those used in \autoref{mhc->ho} produces the following result.
\begin{prop}
  Let $\B_\V$ be the \V-enriched bicategory whose 0-cells are very
  good small \V-categories, whose 1-cells from $\D$ to $\E$ are
  distributors (i.e.\ \V-functors $\E\op\ten\D\to\V$), and whose
  2-cells are \V-natural transformations.  Then $\B_\V$ has a homotopy
  bicategory $\Ho(\B_\V)$ which is $\Ho\V$-enriched.
\end{prop}
This ``homotopical bicategory'' has recently arisen in connection with
unpublished work by Kate Ponto involving fixed point theory.  The left
and right homotopy Kan extensions considered above equip this homotopy
bicategory with ``base change functors'' similar to those considered
in~\cite{pht}.


\section{How to prove goodness}
\label{sec:cofibrancy}

In this section, we make the case that being good or very good is
basically a cofibrancy condition on the small category \D.  The most
appropriate notion of ``cofibration,'' however, can vary widely.  In
particular, even when a model structure exists, the cofibrations of
the model structure are not necessarily the best cofibrations to use
for this purpose.  In this section, we use the term ``cofibration''
without prejudice as to meaning, adopting from~\cite{ekmm}
and~\cite{pht} the terminology \emph{$q$-cofibration} for the model
structure cofibrations.  Another common type of cofibration is an
\emph{$h$-cofibration} (for ``Hurewicz''), generally defined by a
homotopy lifting property.

Our goal is to justify the following meta-statements.

\begin{metadefn}\label{metadef:coft-cat}
  Suppose that \V\ comes equipped with a notion of ``cofibration'' and
  a collection of ``good objects''.  Say a small \V-category \D\ is
  \emph{cofibrant} if each hom-object $\D(d,d')$ is good and each unit
  inclusion $E\to\D(d,d)$ is a cofibration.
\end{metadefn}

\begin{metathm}\label{metathm:coft-good}
  In the above situation, if $(\oast,\homl,\homr)\maps \M\ten\N\to\p$
  is a \tvva\ and \M, \N, and \p\ also come equipped with notions of
  ``cofibration'' and collections of ``good objects,'' then any
  cofibrant small \V-category \D\ is good for $\oast$.
\end{metathm}

One situation in which these meta-statements apply is when \V, \M, \N,
and \p\ are model categories, and the cofibrations and cofibrant
objects are the $q$-cofibrations and $q$-cofibrant objects.  In this
case, we will prove a version of \autoref{metathm:coft-good} as
\autoref{model-coft-good}.  However, there are other situations in which
the metatheorem applies, but the cofibrations in question do not arise
from a model structure.

For example, it applies to the $\mathbb{L}$-spectra and $R$-modules
of~\cite{ekmm} when the cofibrations are the Hurewicz cofibrations
(maps with the homotopy extension property) and every spectrum is
``good.''  It applies to Lewis-May spectra with the Hurewicz
cofibrations, taking the ``good'' spectra to be the \emph{tame}
spectra.  It also applies to the ``well-grounded topological
categories'' of~\cite[ch.~5]{pht}, when the cofibrations are the
$cyl$-cofibrations and the ``good objects'' are the well-grounded
objects.  These situations are all topological, but we suspect that
similar results are true in other contexts, such as for chain
complexes.

A proof of goodness generally has three steps.
\begin{enumerate}
\item Show that if \D\ is cofibrant and $F$ and $G$ are objectwise
  cofibrant, then $\sB(G,\D,F)$ is Reedy cofibrant.\label{good-coft}
\item Show that weak equivalences between functors give rise to weak
  equivalences between simplicial bar constructions.\label{good-we}
\item Show that geometric realization takes weak equivalences between
  Reedy cofibrant objects to weak equivalences between cofibrant
  objects.\label{good-htpical}
\end{enumerate}

The notion of \emph{Reedy cofibrant} is an obvious generalization of
the model-category definition.  When given a notion of cofibration, we
say a simplicial object $\sX$ is Reedy cofibrant if each latching map
$L_nX\to X_n$ is a cofibration.  Note that in the topological
literature, when the cofibrations are the Hurewicz cofibrations, or
$h$-cofibrations, Reedy $h$-cofibrant simplicial objects have usually
been called \emph{proper} simplicial objects.  In this case, the
latching objects $L_nX$ are often written as $sX_n$.

To clarify the essential points, we first give a more general
definition of bar constructions.  This is the most general definition
of bar constructions we have seen; it includes the two-sided bar
construction we have been using, as well as both the monadic bar
constructions used in~\cite{may:goils} and the generalized bar
constructions of~\cite{meyer_i, meyer_ii}.

\begin{defn}\label{def:functor-bimodule}
  Let $(\C,\ten,E)$ be a monoidal category and $D$ a monoid in it.  A
  functor $K\maps \C\op\times\C\to\p$ is called a \emph{$D$-bimodule}
  if it is equipped with two natural transformations
  \begin{align*}
    K(D\ten X) &\too K(X)\\
    K(X\ten D) &\too K(X)
  \end{align*}
  each satisfying evident associativity and unit conditions and
  related by the bimodule commutative square:
  \begin{equation*}
    \xymatrix{
      K(D\ten X \ten D) \ar[r]\ar[d] & K(D\ten X)\ar[d]\\
      K(X\ten D)\ar[r] & K(X).
    }
  \end{equation*}
\end{defn}

The reader should keep the following example in mind.

\begin{exmp}\label{eg:twosided-as-bimodule}
  Let \V\ be a monoidal category, let \Oscr\ be a set, and let \C\ be
  the category of \emph{\Oscr-graphs in \V}.  This is just the
  category $\V^{\Oscr\times\Oscr}$.  We equip \C\ with the following
  monoidal product:
  \begin{equation*}
    (\D\ten\E)(a,b) = \coprod_{c\in\Oscr} \E(c,b)\ten\D(a,c).
  \end{equation*}
  Then a monoid in \C\ is precisely a small \V-category with object
  set \Oscr.

  Let \D\ be such a small \V-category with object set \Oscr, and let
  $F\maps \D\to\N$, $G\maps \D\op\to\M$, and $\oast\maps
  \M\ten\N\to\p$ be \V-functors.  We define a functor $K\maps \C\to\p$ as
  follows.
  \begin{equation*}
    K(\E) = \coprod_{a,b\in\Oscr} \E(a,b)\odot\big(G(b)\oast F(a)\big).
  \end{equation*}
  Then the action of \D\ on $F$ and $G$ makes $K$ into a \D-bimodule
  as above.
\end{exmp}

\begin{defn}
  Given a $D$-bimodule $K$ as in \autoref{def:functor-bimodule}, where
  \p\ is simplicially enriched, we define the \emph{bar construction}
  $B(D,K)$ to be the geometric realization of the following simplicial
  bar construction.
  \begin{equation*}
    B_n(D,K) = K(\overbrace{D\ten\dots\ten D}^{n})
  \end{equation*}
  Here the degeneracies are induced by the unit map of $D$, the inner
  faces by the multiplication of $D$, and the outer faces by the
  action of $D$ on $K$.
\end{defn}

In the situation of \autoref{eg:twosided-as-bimodule}, this gives
precisely the enriched two-sided bar construction $B(G,\D,F)$ defined
in \S\ref{sec:enriched-bar}.

We now begin our proofs of goodness.  We have the following general
result about cofibrancy of generalized bar constructions.

\begin{prop}\label{gen-bar-coft}
  In the situation of \autoref{def:functor-bimodule}, suppose that \C\
  and \p\ have subcategories of ``cofibrations'' satisfying the
  following properties.
  \begin{itemize}
  \item The unit map $E\to D$ is a cofibration in \C;
  \item The pushout product of $E\to D$ with any other cofibration in
    \C\ is a cofibration in \C;
  \item The functor $D\ten -$ preserves pushouts (for example, when
    \C\ is closed monoidal);
  \item $K$ preserves pushouts and cofibrations; and
  \item The map $\emptyset\to K(E)$ is a cofibration in \p.
  \end{itemize}
  Then the simplicial bar construction $\sB(D,K)$ is Reedy cofibrant
  in \p.
\end{prop}
\begin{proof}
  We must show that each map
  \[L_n\sB(D,K)\too B_n(D,K)\]
  is a cofibration in \p.  When $n=0$, we have
  $L_0\sB(D,K)=\emptyset$, and $B_0(D,K)=K(E)$, so the final condition
  above gives the desired result in this case.  For $n\ge 1$, we
  observe that since $K$ preserves pushouts, the desired map is $K$
  applied to a map
  \begin{equation}
    L_nD \to D^n\label{eq:latching-of-monoid}
  \end{equation}
  in \C.  Since $K$ preserves cofibrations, it suffices to show that
  the maps~(\ref{eq:latching-of-monoid}) are cofibrations in \C.  We
  do this by induction on $n$, starting with $n=1$.  We could have
  started our induction with $n=0$, and omitted the last condition
  above, if we were willing to assume that $K$ preserves initial
  objects and that the map $\emptyset\to E$ is a cofibration in \C, as
  are frequently the case.

  When $n=1$, the map~(\ref{eq:latching-of-monoid}) is simply $E\to
  D$, which was assumed to be a cofibration.  So suppose that $L_nD\to
  D^n$ is a cofibration; then we have the following pushout square.
  \begin{equation*}
    \xymatrix{
      L_nD \ar[r]\ar[d] & D^n\ar[d]\ar[ddr]\\
      D\ten L_nD \ar[r]\ar[rrd] & L_{n+1}D \ar[dr]\\
      &&D^{n+1}
    }
  \end{equation*}
  The assumption that $D\ten -$ preserves pushouts ensures that the
  pushout here is, in fact, $L_{n+1}D$.  The induced map from the
  pushout is what we want to show is a cofibration, but it is
  evidently a pushout product of the given cofibration $L_nD\to D^n$
  with $E\to D$, hence a cofibration by assumption.  Thus the
  maps~(\ref{eq:latching-of-monoid}) are cofibrations for all $n\ge 1$
  and hence $\sB(D,K)$ is Reedy cofibrant, as desired.
\end{proof}

\begin{exmp}\label{eg:twosided-coft}
  Consider the situation of \autoref{eg:twosided-as-bimodule}.  Suppose
  that \V\ and \p\ have subcategories of ``cofibrations'' and
  collections of ``good objects'' satisfying the following conditions.
  \begin{itemize}
  \item The tensor product $\ten$ of \V\ preserves colimits;
  \item $\ten$ preserves good objects;
  \item The tensor $\odot$ of \p\ preserves good objects and
    cofibrations between good objects;
  \item The coproduct in \p\ of good objects is good;
  \item The pushout product of cofibrations in \V\ is a cofibration;
    and
  \item The coproduct of cofibrations in \V\ and in \p\ is a
    cofibration.
  \end{itemize}
  Suppose also that $\D$ is a small \V-category and that $G\maps
  \D\op\to\M$, $F\maps \D\to\N$, and $\oast\maps \M\ten\N\to\p$ are
  \V-functors satisfying the following conditions.
  \begin{itemize}
  \item Each map $\emptyset\to\D(d,d')$ and each unit inclusion
    $E\to\D(d,d)$ is a cofibration in \V;
  \item Each object $\D(d,d')$ is good;
  \item $G$ and $F$ are objectwise good and objectwise cofibrant
    (i.e.\ the maps $\emptyset\to Gd$ and $\emptyset\to Fd$ are
    cofibrations); and
  \item $\oast$ preserves good objects and cofibrant objects.
  \end{itemize}
  Then it is straightforward to check that the conditions of
  \autoref{gen-bar-coft} are satisfied, so that the two-sided simplicial
  bar construction $\sB(G,\D,F)$ is Reedy cofibrant in \p.  Moreover,
  in this case it is objectwise good; this will be needed below.
\end{exmp}

This completes step~(\ref{good-coft}).  Step~(\ref{good-we}), that
objectwise weak equivalences in $G$ and $F$ give rise to weak
equivalences between simplicial bar constructions, is easily handled
by the following observation.

\begin{prop}\label{twosided-good-we}
  In the situation of \autoref{eg:twosided-as-bimodule}, suppose that \V,
  \M, \N, and \p\ are homotopical and have collections of ``good''
  objects satisfying the following properties.
  \begin{itemize}
  \item The bifunctor $\oast$, and the tensor $\odot$ of \p, preserve
    good objects and weak equivalences between good objects;
  \item The coproduct of weak equivalences between good objects in \p\
    is a weak equivalence; and
  \item Each object $\D(d,d')$ is good.
  \end{itemize}
  Then objectwise weak equivalences $G\we G'$ and $F\we F'$ between
  functors which are objectwise good give rise to an objectwise weak
  equivalence
  \[\sB(G,\D,F)\we \sB(G',\D,F').\]
\end{prop}
\begin{proof}
  Straightforward from \autoref{def:enr-simp-bar}.
\end{proof}

We now turn to step~(\ref{good-htpical}).  The following proposition
is simply a categorical reformulation of the classical proofs, such as
that in~\cite[X.2.4]{ekmm}.

\begin{prop}\label{realiz-htpical}
  Suppose that \p\ is a homotopical category, simplicially enriched,
  and equipped with a subcategory of ``cofibrations'' and a collection
  of ``good objects'' satisfying the following properties.
  \begin{itemize}
  \item The simplicial tensor $\odot$ of \p\ preserves good objects,
    cofibrations between good objects, and weak equivalences between
    good objects in \p;
  \item  The pushout product under $\odot$ of a cofibration of
    simplicial sets (i.e.\ a monomorphism) and a cofibration in \p\ is
    a cofibration in \p;
  \item Pushouts of cofibrations in \p\ are cofibrations;
  \item (Gluing Lemma) If in the following diagram:
    \begin{equation*}
      \xymatrix{
        \cdot\ar[d]^{\sim} &
        \cdot\ar[l]\ar[d]^{\sim}\ar@{ >->}[r] &
        \cdot\ar[d]^{\sim}\\
        \cdot & \cdot\ar[l]\ar@{ >->}[r] & \cdot
      }
    \end{equation*}
    all objects are good, the maps displayed as $\cof$ are
    cofibrations, and the vertical maps are weak equivalences, then
    the pushouts are good and the induced map of pushouts is a weak
    equivalence.  In particular, weak equivalences between good
    objects are preserved by pushouts along cofibrations.
  \item (Colimit Lemma) If $M$ and $N$ are the colimits of sequences
    of cofibrations $M_k\cof M_{k+1}$ and $N_k\cof N_{k+1}$ between
    good objects, then $M$ and $N$ are good.  Moreover, a compatible
    sequence of weak equivalences $f_k\maps M_k\we N_k$ induces a weak
    equivalence $f\maps M\we N$.
  \end{itemize}
  Then if $\sX\to\sY$ is an objectwise weak equivalence between
  simplicial objects of \p\ which are Reedy cofibrant and objectwise
  good, its realization $|X|\to |Y|$ is a weak equivalence.
\end{prop}
\begin{proof}
  Let $\DD_{<n}$ be the full subcategory of $\DD$ spanned by the
  objects $m$ such that $m<n$, and for a simplicial object
  $\sX\maps \DD\op\to\M$ write $X_{<n}$ for its restriction to $\DD_{<n}$.
  Write
  \begin{equation*}
    |X|_n = \Delta^{<n}\odot_{\DD_{<n}\op} X_{<n}.
  \end{equation*}
  for the ``partial realization'' of \sX.  It is straightforward to
  check that we have
  \begin{equation*}
    |X| = \colim_{n\to\infty} |X|_n.
  \end{equation*}
  In other words, the geometric realization is ``filtered by
  simplicial degree.''  Therefore, by the colimit lemma, if we can
  show that each induced map $|X|_n \to |Y|_n$ is a weak equivalence
  between good objects, the desired result will follow.

  We prove this by induction on $n$.  The case $n=0$ follows directly
  because $|X|_0 = X_0$.  So assume that $|X|_{n-1}\to |Y|_{n-1}$ is a
  weak equivalence between good objects.  Consider the following
  sequence of pushouts, which it is straightforward to check produces
  $|X|_n$ from $|X|_{n-1}$.
  \begin{equation*}
    \xymatrix{
      \partial\Delta^n \odot L_nX \ar[r] \ar[d] &
      \partial\Delta^n \odot X_n \ar[d]\\
      \Delta^n \odot L_n X \ar[r] &
      P_nX \ar[d]\ar[r] &
      \Delta^n \odot X_n \ar[d]\\
      & |X|_{n-1} \ar[r] &
      |X|_n.
    }
  \end{equation*}
  Since $\odot$ preserves good objects and cofibrations between good
  objects, the first pushout is a pushout along a cofibration between
  good objects, and since $\odot$ preserves weak equivalences between
  good objects in \p, the gluing lemma implies that the induced map
  $P_nX\to P_nY$ is a weak equivalence between good objects.

  Now, the map $P_nX \to \Delta^n\odot X_n$ is the pushout product of
  $\partial\Delta^n\to \Delta^n$ and $L_nX\to X_n$, which are both
  cofibrations, hence is a cofibration.  The fact that $\odot$
  preserves weak equivalences between good objects, combined with the
  induction assumption, implies that the conditions of the gluing
  lemma are again satisfied for the second pushout.  Thus,
  $|X|_n\to|Y|_n$ is a weak equivalence between good objects, as
  desired.
\end{proof}

It follows that if all the conditions in the above propositions are
satisfied, then $B(-,\D,-)$ is homotopical on objectwise good
diagrams.  The other two conditions in \autoref{def:good} also follow, as
do those in \autoref{def:very-good}, if the ``good'' objects in the above
propositions are the left deformation retracts $\V_Q$, $\M_Q$, $\N_Q$,
and $\p_Q$ and the left adjoint $\sS\to\V$ lands in $\V_Q$.  In
particular, we have the following special case.

\begin{defn}\label{def:coft-cat}
  Let \V\ be a monoidal model category.  A small \V-category \D\ is
  \emph{$q$-cofibrant} if each hom-object $\D(d,d')$ is $q$-cofibrant
  in \V\ and all the unit inclusions $E\to\D(d,d)$ are
  $q$-cofibrations in \V.
\end{defn}

\begin{thm}\label{model-coft-good}
  Let \V\ be a simplicial monoidal model category, \M, \N, and \p\ be
  \V-model categories, and $(\oast,\homl,\homr)\maps  \M\ten\N\to\p$ be a
  Quillen \tvva.  Then a $q$-cofibrant \V-category \D\ is good for any
  Quillen \tvva, and very good for the tensor and cotensor of any
  \V-model category.
\end{thm}
\begin{proof}
  Taking the cofibrations to be the $q$-cofibrations and the good
  objects to be the $q$-cofibrant objects, it is well-known that all
  the conditions of \autoref{eg:twosided-coft}, \autoref{twosided-good-we},
  and \autoref{realiz-htpical} are satisfied.
\end{proof}

This completes our justification of the meta-statements at the
beginning of this section.

In the dual situation, to prove goodness for cotensor products, the
arguments are analogous.  The cofibrations in \V\ remain cofibrations,
but in \p\ we must consider fibrations instead, with properties dual
to those enumerated above.  Similarly, realization of simplicial
objects becomes totalization of cosimplicial objects, expressed as an
inverse limit of partial totalizations.


Note also that when we consider goodness for the hom-functor
$\M(-,-)\maps \M\op\ten\M\to\V$ of a \V-homotopical category, \M\
itself no longer needs any notion of cofibration or fibration.
Instead, \V\ must have both cofibrations and fibrations which interact
in a suitable way.  Thus, if \V\ is nice enough, \emph{all}
\V-homotopical categories have a well-behaved homotopy theory of
diagrams.  This shows the essential role of the enrichment, which puts
the real ``homotopy'' in homotopy theory.

\section{Objectwise good replacements}
\label{sec:objectwise-good}

In this section, we investigate under what conditions $\M_Q^\D$
is a left deformation retract of $\M^\D$, so that every diagram can be
replaced by one which is ``objectwise good.''  Our first observation
is trivial.

\begin{prop}\label{objw-good-unenriched}
  If \D\ is not enriched, then $\M_Q^\D$ is a left deformation retract
  of $\M^\D$.
\end{prop}
\begin{proof}
  Apply the functor $Q$ objectwise.
\end{proof}

The problem in the enriched case is that frequently in the case of
model categories, the functor $Q$ is produced by a small object
argument.  While the small object argument excels at producing
functorial replacements, it generally fails to produce \emph{enriched}
functorial replacements, and when $Q$ is not an enriched functor, we
cannot compose it with an enriched functor $F\maps \D\to\M$ to produce
a new enriched functor $QF\maps \D\to\M$.

There are some cases, however, in which $Q$ can be chosen to be an
enriched functor.  One trivial case is when $\M_Q=\M$ (every object is
cofibrant), so that $Q$ can be chosen to be the identity.  Dually, if
$\M_R=\M$, then $R$ can be chosen to be the identity.  This latter
case arises frequently in topology, which makes homotopy limits easier
to deal with; this is fortunate because, as we have mentioned
frequently, our results are generally more necessary for limits than
for colimits.

Another case in which $Q$ can be chosen enriched is when every object
of \V\ is cofibrant.  This allows us to actually perform a
small-object-argument factorization in an enriched-functorial way.
The only commonly occurring \V{}s for which this is true are $\V=\sS$
and $\V=\mathbf{Cat}$, but in these cases it simplifies things
greatly.  We briefly sketch a proof below; it is essentially the same
as \cite[4.3.8]{hirschhorn}, although the latter is stated
specifically for localizations of simplicial model categories.


\begin{prop}\label{objw-good-allcoft}
  Let \V\ be a monoidal model category, let \M\ be a \V-model category
  which is cofibrantly generated, and assume that for all $K\in \V$,
  the tensor $K \odot-$ preserves cofibrations in \M.  For example,
  this is true if all objects of \V\ are cofibrant.  Then there exists
  a \V-functorial cofibrant replacement functor $Q$ on \M.  Dually, if
  $K\odot -$ preserves trivial cofibrations for all $K$, there is a
  \V-functorial fibrant replacement functor $R$ on \M.
\end{prop}
\begin{proof}
  The construction is essentially the same as the usual small object
  argument.  Given $X$, we construct $QX$ as the colimit of a
  transfinite sequence $\{X_\beta\}$ starting with $X_0=\emptyset$ and
  using the following pushouts.
  \begin{equation}\label{eq:allcoft-smallobj}
    \xymatrix{
      \smash{\displaystyle\coprod_{A\to B}} \left(\M(A,X_\beta)\times_{\M(A,X)}\M(B,X)\right)\odot A \ar[r]\ar[d] &
      X_\beta \ar[d]\\
      \smash{\displaystyle\coprod_{A\to B}} \left(\M(A,X_\beta)\times_{\M(A,X)}\M(B,X)\right)\odot B \ar[r] &
      X_{\beta+1}
    }
  \end{equation}\ 

  \noindent Here the coproducts are over all maps $A\to B$ in a set of
  generating cofibrations.  The crucial difference from the usual
  small object argument is the use of the enriched `objects of
  commuting squares' $\M(A,X_\beta)\times_{\M(A,X)}\M(B,X)$ instead of
  a coproduct over all such commuting squares in $\M_0$.  It is
  straightforward to check by induction that this change makes
  $X\mapsto QX$ \V-functorial, and the usual arguments work to show
  that $QX\to X$ is a trivial fibration.  The assumption that $\odot$
  preserves cofibrations is needed to ensure that the left-hand maps
  are all still cofibrations, so that $QX$ is in fact cofibrant.  The
  construction of $R$ is dual.
\end{proof}




Unfortunately, having to use a small-object-argument factorization to
correct the bar construction makes it much less explicit, thus
partially nullifying one of its advantages.  However, as noted after
\autoref{def:hoten}, if the diagrams we begin with are already objectwise
cofibrant, there is no need for a correction at all.  The mere
\emph{existence} of the deformation $Q_\M^\D$ tells us that the
results we get this way agree with the global derived functor.

There are other situations in which $Q$ can be made \V-functorial,
notably those in which a weaker notion than model-category-theoretic
cofibrancy suffices.  In these cases, the replacement functors are
often given by explicit constructions, which can usually be made
\V-functorial.  Frequently, in fact, these explicit constructions are
themselves a sort of bar construction.


There is also a situation in which $\M_Q^\D$ is a deformation retract
of $\M^\D$ even though $Q$ is \emph{not} a \V-functor, and that is
when $\M^\D$ admits a model structure of its own in which the
cofibrant diagrams are, in particular, objectwise cofibrant.  In this
case, a cofibrant replacement functor for $\M^\D$ serves as a
deformation into $\M_Q^\D$.

As remarked in \S\ref{sec:enriched}, the only type of model structure
we can reasonably expect to have on an enriched diagram category is a
projective one, but these exist in reasonable generality.  There is,
however, something a little strange about the idea of applying a
projective-cofibrant replacement functor before the bar construction,
since projective cofibrancy is sometimes a good enough deformation for
homotopy colimits on its own.  However, there are a couple of reasons
why this is still a useful thing to consider.

Firstly, projective-cofibrant replacements are generally quite hard to
compute with.  But as in the situation of \autoref{objw-good-allcoft}, if
the diagram we started with is already objectwise cofibrant, we don't
need to make it projective-cofibrant before applying the bar
construction.  The value of the existence of the projective model
structure for us, then, is that it tells us that the bar construction
actually defines a derived functor on the whole diagram category,
rather than just the objectwise cofibrant diagrams.


Secondly, projective model structures are \emph{not} directly useful
for computing homotopy \emph{limits}, because the limit functor is not
right Quillen for them.  Moreover, as noted in \S\ref{sec:enriched},
injective model structures do not generally exist in the enriched
situation.  However, by combining a projective model structure with
the cobar construction, we can compute homotopy limits: since
projective-fibrant diagrams are, by definition, objectwise fibrant, a
projective-fibrant replacement works perfectly well to correct the
cobar construction.

With this as prelude, we now state a theorem about the existence of
projective model structures on enriched diagram categories.  This is a
straightforward application of the principle of transfer of model
structures along adjunctions, which was first stated in~\cite{crans}.
Various special cases of this result are scattered throughout the
literature, but we have been unable to find a statement of it in full
generality.

For example, essentially the same result for the case $\V=\M$ can be
found in~\cite[6.1]{ss:equiv-mmc}.  In the case when \D\ has one
object and $\V=\M$, it reduces to a result like those
of~\cite{hovey:mmc} and~\cite{ss:alg-mod-mmc} for modules over a
monoid in a monoidal model category.  It can also be viewed as a
special case of the model structures considered
in~\cite{brgmoer:coloured} for algebras over colored operads.

\begin{thm}\label{proj-model-str}
  Let \V\ be a monoidal model category and let \M\ be a cofibrantly
  generated \V-model category with generating cofibrations $I$ and
  generating trivial cofibrations $J$, whose domains are assumed small
  with respect to all of \M.  Let \D\ be a small \V-category such that
  $\D(d,d')\odot -$ preserves trivial cofibrations (for example, this
  occurs if each $\D(d,d')$ is cofibrant in \V, or if $\M=\V$ and \V
  satisfies the ``monoid axiom'').  Then $\M^\D$ has a cofibrantly
  generated \V-model structure in which the weak equivalences and
  fibrations are objectwise.  If $\D(d,d')\odot -$ also preserves
  cofibrations (as also occurs if each $\D(d,d')$ is cofibrant) then
  cofibrant diagrams are in particular objectwise cofibrant.
\end{thm}
\begin{proof}
  We use~\cite[Theorem 3.3]{crans}, or more precisely its restatement
  in~\cite[2.5]{brgmoer:operads}.  The forgetful functor $U\maps \M^\D
  \to \M^{\ob\D}$ has both left and right adjoints, given by left and
  right Kan extension.  Write $F$ for the left adjoint; we want to use
  $F$ to induce a model structure on $\M^\D$.

  First we observe that $\M^{\ob\D}$ has a model structure which is
  just a power of the model structure on \M.  This model structure is
  cofibrantly generated; sets of generating cofibrations and trivial
  cofibrations are given by
  \begin{align*}
    I' &= \{L_d f| f\in I,d\in\ob\D\}\\
    J' &= \{L_df| f\in J,d\in\ob\D\}
  \end{align*}
  where $L_df$ consists of $f$ at spot $d$ and $\emptyset$ elsewhere.

  Since $U$ has a right adjoint, it preserves colimits, and therefore
  its left adjoint $F$ preserves smallness of objects.  Therefore it
  remains to check that relative $FJ'$-cell complexes are weak
  equivalences in $\M^\D$.  However, $U$ preserves colimits and
  reflects weak equivalences, so it suffices to check that $UFJ'$-cell
  complexes are weak equivalences in $\M^{\ob\D}$.  Now by definition
  of $F$ as a left Kan extension, we have
  \[UF(L_df)_{d'} = \D(d,d')\odot f\]
  and so our assumption ensures that the maps in $UFJ'$ are trivial
  cofibrations, hence relative $UFJ'$-cell complexes are also trivial
  cofibrations in $\M^{\ob\D}$ and thus weak equivalences.

  This shows that the model structure exists.  To see that it is a
  \V-model structure, recall from~\cite[4.2.2]{hovey} that it suffices
  to check the lifting extension properties for the cotensor.  But
  since the fibrations, trivial fibrations, and cotensor are all
  objectwise, this follows from the corresponding property for \M.
  Similarly, the unit condition, when expressed using the cotensor,
  follows from the unit condition for \M.

  For the last statement, we use an adjointness argument to show more
  generally that projective-cofibrations are objectwise cofibrations.
  Consider $U$ and its \emph{right} adjoint which we denote $G$; we
  want to show that $U$ preserves cofibrations, which is equivalent to
  showing that $G$ preserves trivial fibrations.  Now any trivial
  fibration in $\M^{\ob\D}$ is a product of maps of the form $R_df$
  for trivial fibrations $f$ in \M, where $R_df$ consists of $f$ at
  spot $d$ and $1$ elsewhere, and by definition of $G$ as a right Kan
  extension, we have
  \[G(R_df)_{d'} = \cten{\D(d',d), f}.\]
  But since $\D(d',d)\odot -$ preserves cofibrations, by adjointness
  $\cten{\D(d',d),-}$ preserves trivial fibrations, so $G$ must
  preserve trivial fibrations, as desired.
\end{proof}

It follows, for instance, that when \D\ is cofibrant, $\M_Q^\D$ and
$\M_R^\D$ are, respectively, a left and a right deformation retract of
$\M^\D$, and so (in addition to the projective \V-model structure) the
bar and cobar \V-homotopical structures exist.  As remarked above, the
projective model structure can sometimes be used to compute homotopy
colimits, but not homotopy limits; for this we need to use the cobar
construction.

\appendix

\section{Proof of Theorem \ref*{dhks-comp}}
\label{sec:ugly-proof}

Here we will complete the proof of \autoref{dhks-comp}.  The proof we
give is slightly modified from the proof given in an early online
draft of~\cite{dhks}.  First, we recall the statement of the theorem.

\bigskip\noindent\usebox{\dhkstheorem}\bigskip

We will need some technical lemmas.  Let $H\maps \E\to\D$ be a functor
between small categories.  Applying the nerve functor $N$, we obtain a
map $NH\maps  N\E\to N\D$ of simplicial sets.  Given a simplex $\alpha$ in
$N\D$, let $\Lambda_H(\alpha)$ be the fiber of $NH$ over $\alpha$, so
that the following square is a pullback:
\begin{equation}\label{eq:lambda-pb}
  \xymatrix{
    \Lambda_H(\alpha) \ar[r]\ar[d] & N\E \ar[d]^{NH}\\
    \Delta^n \ar[r]_{\alpha} & N\D
  }
\end{equation}
This defines a functor $\Lambda_H\maps \Delta\D\to\sS$.
\begin{lem}\label{lan-comma}
  In the above situation, we have $\lan_{S\op}\Lambda_H \iso N(-\dn H)$.
\end{lem}
\begin{proof}
  As in the proof of \autoref{bar=may}, left Kan extension can be
  computed as a colimit over a comma-category.  Here, the formula
  is:
  \begin{equation}
    (\lan_{S\op}\Lambda_H)(d) \iso \colim_{(S\op\dn d)} \Lambda_H.
  \end{equation}
  Note moreover that we have $(S\op\dn d) \iso \Delta(d\dn\D)$.  Now
  the colimit~(\ref{eq:lan-fmla}) in \sS\ can also be viewed as a
  colimit in $(\sS\dn N\D)$, and pullback along $NH$, viewed as a
  functor $(\sS\dn N\D) \to (\sS\dn N\E)$, preserves colimits because
  it has a right adjoint.  Thus, applying colimits to the left
  vertical arrow in~(\ref{eq:lambda-pb}), and noting that
  $\colim_{\Delta\C}\Delta^n \iso N\C$ for any category \C, we get the
  following pullback square.
  \begin{equation*}
    \xymatrix{
      \lan_{S\op}\Lambda_H \ar[d]\ar[r] & N\E\ar[d]^{NH} \\
      N(c\dn\D)\ar[r] & N\D
    }
  \end{equation*}
  It is straightforward to check that there is also a similar pullback
  square with $\lan_{S\op}\Lambda_H$ replaced by $N(c\dn H)$.
\end{proof}

A particular case of this lemma is when $\E=\D$ and $H = \Id_\D$.
Then we have $\Lambda_H = \Delta\Pi\maps \Delta\D \to \sS$ (where
$\Delta\maps \DD\to\sS$ is the canonical cosimplicial simplicial set,
and $\Pi$ is the forgetful functor from~(\ref{eq:simpl-diag})), and
the result says that
\begin{equation}\label{eq:lan-under}
  \lan_{S\op}\Delta\Pi \iso N(-\dn \D).
\end{equation}

Now, given two simplices $\alpha$ and $\beta$ in $\D$, let
$\Lambda(\alpha, \beta)$ be their pullback over $N\D$, so that the
following diagram is a pullback:
\begin{equation*}
  \xymatrix{
    \Lambda(\alpha,\beta) \ar[r] \ar[d] &
    \Delta^m \ar[d]^\beta\\
    \Delta^n \ar[r]^\alpha & N\D.
  }
\end{equation*}
This defines a functor $\Lambda\maps \Delta\D\times\Delta\D\to\sS$.

\begin{lem}\label{lan-comma-2}
  \begin{equation*}
    \colim_{\Delta\D} \Lambda \iso \Delta\Pi\maps \Delta\D\to\sS
  \end{equation*}
\end{lem}
\begin{proof}
  Similar to the proof of \autoref{lan-comma}.
\end{proof}

\begin{lem}\label{lan-adjt}
  For any functors $K\maps \D\to\E$, $G\maps \D\op\to\sS$, and $F\maps
  \E\to\M$,
  \begin{equation*}
    \left(\lan_{K\op} G\right) \odot_\D F \iso G \odot_\E \left(K^* F\right).
  \end{equation*}
\end{lem}
\begin{proof}
  Note that because the tensor $\odot$ has adjoints on both sides (the
  cotensor and the enriched hom-functor), it preserves colimits in
  both variables.  The result follows from this, the fact that
  colimits commute with colimits, and the fact that $\E(-,Kd)\odot_\E
  F \iso F(Kd)$, as remarked in the proof of
  \autoref{one-sided-bar-pullout}.
\end{proof}

We now attack the proof of \autoref{dhks-comp}.  Recall that in
\S\ref{sec:comp-dhks}, we reduced the proof to
showing~(\ref{eq:dhks-repl}) and that $\Q F$ is Reedy cofibrant.  We
combine the above lemmas to manipulate $\hocolim$ into a form looking more
like $\Q F$ as follows.
\begin{xalignat*}{2}
    \hocolim F
    &= N(-\dn\D) \odot_\D QF
& \text{(by definition)}\\
    &\iso \lan_{S\op} \Delta\Pi \odot_{\D} QF
& \text{(by~(\ref{eq:lan-under}))}\\
    &\iso \Delta\Pi \odot_{\Delta\D} S^* QF
& \text{(by \autoref{lan-adjt})}\\
    &\iso \Big(\colim_{\alpha\in\Delta\D}\Lambda(\alpha,-)\Big)
    \odot_{\Delta\D} S^* QF 
& \text{(by \autoref{lan-comma-2})}\\
    &\iso \colim_{\alpha\in\Delta\D}
    \Big(\Lambda(\alpha,-)\odot_{\Delta\D} S^* QF\Big)
\end{xalignat*}
because colimits commute with tensor products of functors.  Thus $\hocolim
F$ is the colimit of the $\Delta\D$-diagram
\begin{equation}
  \label{eq:calq-def}
  \Lambda(\alpha,-)\odot_{\Delta\D} S^* QF
\end{equation}
and it remains to
identify that diagram with $\Q F$.  To begin with, we have
\begin{xalignat*}{2}
  \Lambda(\alpha,-)\odot_{\Delta\D} S^* QF
  &\iso \lan_{S\op}\Lambda(\alpha,-) \odot_\D QF
& \text{(by \autoref{lan-adjt})}\\
    &\iso N(-\dn\alpha)  \odot_\D  QF
& \text{(by \autoref{lan-comma}).}
\end{xalignat*}
Now note that for any functor $H\maps \E\to\D$, we have
\begin{equation*}
  \lan_{H\op} N(-\dn\E) \iso N(-\dn H).
\end{equation*}
This is easily verified using the formula $N(-\dn H) \iso B(*,\D,H)$.
Therefore for $H=\alpha$, we have:
\begin{xalignat*}{2}
  N(-\dn\alpha)  \odot_\D  QF
    &\iso \Big(\lan_{\alpha\op} N(-\dn[n]) \Big)  \odot_{[n]} QF \\
    &\iso N(-\dn[n])  \odot_{[n]} \alpha^* QF
& \text{(by \autoref{lan-adjt})}\\
    &\iso \Q F,
\end{xalignat*}
as desired.


Finally, we want to know that $\Q F$ is Reedy cofibrant.  Here we must
use some more technical model category theory.  Recall the following
result from~\cite[18.4.11]{hirschhorn}:
\begin{lem}\label{simplicial-reedy-2va}
  Let \R\ be a Reedy category and \M\ a simplicial model category.  If
  $G_1\to G_2$ is a Reedy cofibration of $\R\op$-diagrams of
  simplicial sets, and $F_1\to F_2$ is a Reedy cofibration of
  \R-diagrams in \M, then the pushout corner map:
  \begin{equation*}
    \big(G_1\odot_\R F_2\big)
    \bigsqcup_{G_1\odot_\R F_1}
    \big(G_2\odot_\R F_1\big)
    \too G_2\odot_\R F_2
  \end{equation*}
  is a cofibration in \M\ that is a weak equivalence if either given
  map is a Reedy weak equivalence.
\end{lem}

Note that since $\Delta$ is a Reedy cofibrant $\DD$-diagram in \sS,
\autoref{simplicial-reedy-2va} implies \autoref{realiz-quillen}.  The
methods used in the proof of \autoref{simplicial-reedy-2va} can be used
to show the following result.  Recall from~\cite[15.5.2]{hirschhorn}
that if $\R_1$ and $\R_2$ are Reedy categories, then $\R_1\times\R_2$
has a canonical Reedy structure, and the Reedy model structure on
$\M^{\R_1\times\R_2}$ is the same whether we regard it as
$\R_1\times\R_2$-diagrams in \M, $\R_1$-diagrams in $\M^{\R_2}$, or
$\R_2$-diagrams in $\M^{\R_1}$.

\begin{lem}\label{simplicial-double-reedy-2va}
  Let $\R_1$ and $\R_2$ be Reedy categories and let \M\ be a
  simplicial model category.  If $G_1\to G_2$ is a Reedy cofibration
  of $\R_1\op$-diagrams of simplicial sets, and $F_1\to F_2$ is a
  Reedy cofibration of $\R_1\times\R_2$-diagrams in \M, then the
  pushout corner map:
  \begin{equation*}
    \big(G_1\odot_\R F_2\big)
    \bigsqcup_{G_1\odot_\R F_1}
    \big(G_2\odot_\R F_1\big)
    \too G_2\odot_\R F_2
  \end{equation*}
  is a Reedy cofibration of $\R_2$-diagrams in \M\ that is a Reedy weak
  equivalence if either given map is a Reedy weak equivalence.
\end{lem}

These lemmas are both special cases of a more general result about
\tva{}s which the reader is invited to state and prove.

Taking $G_1$ and $F_1$ to be diagrams constant at the initial object,
we see that the tensor product of two Reedy cofibrant diagrams is
Reedy cofibrant.  Now recall from~(\ref{eq:calq-def}) that we have
\begin{equation*}
  \Q F(\alpha) = \Lambda(\alpha,-)\odot_{\Delta\D} S^* QF.
\end{equation*}
By \autoref{objw-coft->reedy}, we know that $S^* QF$ is a Reedy cofibrant
$\Delta\op\D$-diagram in \M.  Thus by
\autoref{simplicial-double-reedy-2va}, to conclude that $\Q F$ is Reedy
cofibrant in $\M^{\Delta\D}$, it suffices to observe that
$\Lambda(-,-)$ is a Reedy cofibrant $(\Delta\D\times\Delta\D)$-diagram
in \sS.

This completes the proof of \autoref{dhks-comp}.  We have given the
proof only for a simplicial model category, for ease of exposition.
However, a careful examination shows that it goes through without a
hitch for general model categories, as long as the simplicial tensor
is interpreted appropriately, using the by-now standard techniques of
framings or resolutions, as described in~\cite[ch.~5]{hovey}
or~\cite[ch.~16]{hirschhorn}.

\nocite{hovey:errata}
\bibliographystyle{alpha}
\raggedright
\bibliography{htpylim}

\end{document}